# 树结构和图同态

**龙旸靖**

| | |
|---|---|
| 合 作 导 师: | 吴耀琨教授 |
| 专 业 名 称: | 数 学 |
| 完 成 日 期: | 2015 年3月–2017 年9 月 |
| 提 交 日 期: | 2017 年9 月 |

**上海交通大学 数学科学学院**

**2017年9月**

# Phylogenetic trees and homomorphisms

**Yangjing Long**

**School of Mathematical Sciences**
**Shanghai Jiaotong University**
**Shanghai, China**

**Semptember 2017**

# 摘　要

第一章中我们完全刻画了全同态序中形成间隔的有限图对。通过这一刻画可以给出推广对偶对的存在性的一个简单证明。这一结果不仅对于无向图成立，并且对于有向图和双边的关系结构也成立。

第二章中我们给出了一个图同态序的普适性的新的简单证明，并且讨论了这一证明的一些应用。

第三章中我们用第二章中的简单证明，证明了图同态序的分形这一美妙的性质。

第四章中我们从组合学的角度分析了系统发生信息。考虑任意两个叶子节点之间的路径上的标号是孤一关系还是非孤一关系。我们证明了表示这种双边关系的图表示一定是树，并且我们完整刻画了这种树和其对应的系统发生树。

第五章中我们给出了符号标记无根树和符合三元度量的一一对应。

**关键词:** 图同态; 同态序; 稠密; 普适性; 全同态; 间隔; 分形性质; 系统发生组合学; 稀有事件; 双边关系; 符合三元度量; 中点; 无根系统发生树

# Abstract

In Chapter 1 we fully characterise pairs of finite graphs which form a gap in the full homomorphism order. This leads to a simple proof of the existence of generalised duality pairs. We also discuss how such results can be carried to relational structures with unary and binary relations.

In Chapter 2 we show a very simple and versatile argument based on divisibility which immediately yields the universality of the homomorphism order of directed graphs and discuss three applications.

In chapter 3, we show that every interval in the homomorphism order of finite undirected graphs is either universal or a gap. Together with density and universality this “fractal” property contributes to the spectacular properties of the homomorphism order.

In Chapter 4 we analyze the phylogenetic information content from a combinatorial point of view by considering the binary relation on the set of taxa defined by the existence of a single event separating two taxa. We show that the graph-representation of this relation must be a tree. Moreover, we characterize completely the relationship between the tree of such relations and the underlying phylogenetic tree.

In 1998, Böcker and Dress gave a 1-to-1 correspondence between symbolically dated rooted trees and symbolic ultrametrics. In Chapter 5 we consider the corresponding problem for unrooted trees. More precisely, given a tree $T$ with leaf set $X$ and a proper vertex colouring of its interior vertices, we can map every triple of three different leaves to the colour of its median vertex. We characterise all ternary maps that can be obtained in this way in terms of 4- and 5-point conditions, and we show that the corresponding tree and its colouring can be reconstructed from a ternary map that satisfies those conditions. Further, we give an additional condition that characterises whether the tree is binary, and we describe an algorithm that reconstructs general trees in a bottom-up fashion.

**Key Words:**

graph homomorphism ; homomorphism order; density; universality; full homomorpism; gaps; fractal property; Phylogenetic Combinatorics; Rare events ; Binary relations; symbolic ternary metric; median vertex ; unrooted phylogenetic tree

# Contents

# Chapter 1

# A brief introduction

This thesis is a combination of several published or submitted papers I wrote together with my coauthors during my postdoc period at Shanghai Jiao Tong University. The papers are listed below and each chapter is based on each paper with the same order number.

Due to the limit of the length of this thesis I decide to skip the introduction and preliminary Chapter. The necessary preliminary and introduction can be found in the beginning of each chapter. One can also find my detailed introduction and preliminary in the Introduction section of my PhD thesis [55].

1. Gaps in full homomorphism order, joint with J. Fiala, J. Hubička, *Electronic Notes in Discrete Mathematics*, (2017), **61**, 429–435

2. An universality argument for graph homomorphisms, joint with J. Fiala, J. Hubička, *Electronic Notes in Discrete Mathematics*, (2015), **49**, 643–649.

3. Fractal property of homomorphism order, joint with J. Fiala, J. Hubička, J. Nešetřil, *European Journal of Combinatorics*, in press, *http://www.sciencedirect.com/science/*

4. Inferring Phylogenetic Trees from the Knowledge of Rare Evolutionary Events, joint with M. Hellmuth, M. Hernandez-Rosales, P. F. Stadler, to appear in *Journal of Mathematical Biology*, available at *https://arxiv.org/abs/1612.09093*

5. Reconstructing unrooted phylogenetic trees from symbolic ternary metrics, joint with S. Grünewald, Y. Wu, minor revision in *Bulletin of Mathematical Biology*, available at *https://arxiv.org/abs/1702.00190*

# Chapter 2

# Gaps in full homomorphism order

## 2.1 Introduction

For given graphs $G = (V_G, E_G)$ and $H = (V_H, E_H)$ a *homomorphism* $f : G \to H$ is a mapping $f : V_G \to V_H$ such that $\{u, v\} \in E_G$ implies $\{f(u), f(v)\} \in E_H$. (Thus it is an edge preserving mapping.) The existence of a homomorphism $f : G \to H$ is traditionally denoted by $G\,\mathrm{Hom}_R(H,)$. This allows us to consider the existence of a homomorphism, $\mathrm{Hom}_R(,,)$ to be a (binary) relation on the class of graphs. A homomorphism $f$ is *full* if $\{u, v\} \notin E_G$ implies $\{f(u), f(v)\} \notin E_H$. (Thus it is an edge and non-edge preserving mapping). Similarly we will denote by $G \stackrel{F}{\to} H$ the existence of a full homomorphism $f : G \to H$.

As it is well known, the relations $\to$ and $\stackrel{F}{\to}$ are reflexive (the identity is a homomorphism) and transitive (a composition of two homomorphisms is still a homomorphism). Thus the existence of a homomorphism as well as the existence of full homomorphisms induces a quasi-order on the class of all finite graphs. We denote the quasi-order induced by the existence of homomorphisms and the existence of full homomorphism on finite graphs by $(\mathsf{Graphs}, \leq)$ and $(\mathsf{Graphs}, \leq^F)$ respectively. (Thus when speaking of orders, we use $G \leq H$ in the same sense as $G\,\mathrm{Hom}_R(H,)$ and $G \leq^F H$ in the sense $G \stackrel{F}{\to} H$.)

These quasi-orders can be easily transformed into partial orders by choosing a particular representative for each equivalence class. In the case of graph homomorphism such representative is up to isomorphism unique vertex minimal element of each class, the *(graph) core*. In the case of full homomorphisms we will speak of *F-core*.

The study of homomorphism order is a well established discipline and one of main topics of nowadays classical monograph of Hell and Nešetřil [34]. The order $(\mathsf{Graphs}, \leq^F)$ is a topic of several publications [83, 22, 33, 3, 25] which are primarily concerned about the full homomorphism equivalent of the homomorphism duality [62].

In this work we further contribute to this line of research by characterising *F-gaps* in $(\mathsf{Graphs}, \leq^F)$. That is pairs of non-isomorphic F-cores $G \leq^F H$ such that every F-core $H'$, $G \leq^F H' \leq^F H$, is isomorphic either to $G$ or $H$. We will show:

**Theorem 2.1.1.** *If $G$ and $H$ are F-cores and $(G, H)$ is an F-gap, then $G$ can be obtained from $H$ by removal of one vertex.*

First we show a known fact that F-cores correspond to point-determining graphs which have been studied in 70's by Sumner [77] (c.f. Feder and Hell [22]). We also show that there is a full homomorphism between two F-cores if and only if there is an embedding from one to another (see [22, Section 3]). These two observations shed a lot of light into the nature of full homomorphism order and makes the characterisation of F-gaps look particularly innocent (clearly gaps in embedding order are characterised by an equivalent of Theorem 2.1.1). The arguments in this area are however surprisingly subtle. This becomes even more apparent when one generalise the question to classes of graphs as done by Hell and Hernández-Cruz [33] where both results of Sumner [77] and Feder and Hell [22] are given for digraphs by new arguments using what one could consider to be surprisingly elaborate (and interesting) machinery needed to carry out the analysis.

We focus on minimising arguments about the actual structure of graphs and use approach which generalises easily to digraphs and binary relational structures in general (see Section 2.5). In Section 2.2 we outline the connection of point determining graphs and F-cores. In Section 2.3 we show proof of the main result. In Section 2.4 we show how the existence of gaps leads to a particularly easy proof of the existence of generalised dualities (main results of [83, 22, 33, 3]).

## 2.2 F-cores are point-determining

In a graph $G$, the *neighbourhood* of a vertex $v \in V_G$, denoted by $N_G(v)$, is the set of all vertices $v'$ of $G$ such that $v$ is adjacent to $v'$ in $G$. *Point-determining graphs* are graphs in which no two vertices have the same neighbourhoods. If we start with any graph $G$, and gradually merge vertices with the same neighbourhoods, we obtain a point-determining graph, denoted by $G_{\mathrm{pd}}$.

We write $G \sim^F H$ for any pairs of graphs such that $G \xrightarrow{F} H$ and $H \xrightarrow{F} G$. It is easy to observe that $G_{\mathrm{pd}}$ is always an induced subgraph of $G$. Moreover, for every graph $G$ it holds that $G_{\mathrm{pd}} \xrightarrow{F} G \xrightarrow{F} G_{\mathrm{pd}}$ and thus $G \sim^F G_{\mathrm{pd}}$. This motivates the following proposition:

**Proposition 2.2.1** ([22])**.** *A finite graph $G$ is an F-core if and only if it is point-determining.*

*Proof.* Recall that $G$ is an F-core if it is minimal (in the number of vertices) within its equivalence class of $\sim^F$. If $G$ is an F-core, $G_{\text{pd}}$ can not be smaller than $G$ and thus $G = G_{\text{pd}}$.

It remains to show that every point-determining graph is an F-core. Consider two point-determining graphs $G \sim^F H$ that are not isomorphic. There are full homomorphisms $f : G \xrightarrow{F} H$ and $g : H \xrightarrow{F} G$. Because injective full homomorphisms are embeddings, it follows that either $f$ or $g$ is not injective. Without loss of generality, assume that $f$ is not injective. Consider $u, v \in V_G$, $u \neq v$, such that $f(u) = f(v)$. Because full homomorphisms preserve both edges and non-edges, the preimage of any edge is a complete bipartite graph. If we apply this fact on edges incident with $f(u)$, we derive that $N_G(u) = N_G(v)$. □

**Proposition 2.2.2** ([22, 35]). *For F-cores $G$ and $H$ we have $G \xrightarrow{F} H$ if and only if $G$ is an induced subgraph of $H$.*

*Proof.* Embedding is a special case of a full homomorphisms. In the opposite direction consider a full homomorphism $f : G \xrightarrow{F} H$. By the same argument as in the proof of Proposition 2.2.1 we get that $f$ is injective, as otherwise $G$ would not be point-determining. □

## 2.3 Main result: characterisation of F-gaps

Given a graph $G$ and a vertex $v \in V_G$ we denote by $G \setminus v$ the graph created from $G$ by removing vertex $v$. We say that vertex $v$ *determines* a pair of vertices $u$ and $u'$ if $N_{G \setminus v}(u) = N_{G \setminus v}(u')$. This relation (pioneered in [22] and used in [22, 33, 83])will play key role in our analysis. We make use of the following Lemma:

**Lemma 2.3.1.** *Given a graph $G$ and a subset $A$ of the set of vertices of $G$ denote by $L$ a graph on the vertices of $G$, where $u$ and $u'$ are adjacent if and only if there is $v \in A$ that determines $u$ and $u'$. Let $S$ be any spanning tree of $L$. Denote by $B \subseteq A$ the set of vertices that determine some pair of vertices connected by an edge of $L$ and by $C \subseteq B$ set of vertices that determine some pair of vertices connected by an edge of $S$. Then $B = C$.*

*Proof.* Because for every pair of vertices there is at most one vertex determining them clearly $C \subseteq B \subseteq A$.

Assume to the contrary that there is vertex $v \in B \setminus C$ and thus every pair determined by $v$ is an edge of $L$ but not an edge of $S$. Denote by $\{u, u'\}$ some such edge of $L$ determined by $v \in B$. Adding this edge to $S$ closes a cycle. Denote by $u = v_1, v_2, \dots v_n = u'$ the vertices of $G$ such that every consecutive pair is an edge of $S$. Without loss of generality, we can assume that $v \in N_G(v_1)$ and $v \notin N_G(v_n)$. Because $v \in N_G(v_i)$ implies $v \in N_G(v_{i+1})$ unless $v$ determines pair $\{v_i, v_{i+1}\}$ we also know that there is $1 \leq i < n$ such that $v$ determines $v_i$ and $v_{i+1}$. A contradiction with the fact that $v_i, v_{i+1}$ forms an edge of $S$. □

As a warmup we show the following theorem which also follows by [77] (also shown as Corollary 3.2 in [22] for graphs and [33] for digraphs):

**Theorem 2.3.2** ([77, 22, 33])**.** *Every F-core $G$ with at least 2 vertices contains an F-core with $|V_G| - 1$ vertices as an induced subgraph.*

*Proof.* Denote by $n$ number of vertices of $G$. If there is a vertex $v$ of $G$ such that the graph $G \setminus v$ is point-determining, it is the desired F-core. Consider graph $S$ as in Lemma 2.3.1 where $A$ is the vertex set of $G$. Because $S$ has at most $n-1$ edges and every edge of $S$ is determined by at most one vertex, we know that there is vertex $v$ which does not determine any pair of vertices and thus $G \setminus v$ is point-determining. □

In fact both [77, 33] shows that every F-core $G$ with at least 2 vertices contains vertices $v_1 \neq v_2$ such that both $G \setminus v_1$ and $G \setminus v_2$ are F-cores. This follows by our argument, too but needs bit more detailed analysis. The main idea of the following proof of Theorem 2.1.1 can also be adapted to show this.

*Proof.* (of Theorem 2.1.1) Assume to the contrary that there are F-cores $G$ and $H$ such that $(G, H)$ is an F-gap, but $G$ differs from $H$ by more than one vertex. By induction we construct two infinite sequences of vertices of $H$ denoted by $u_0, u_1, \dots$ and $v_0, v_1, \dots$ along with two infinite sequences of induced subgraphs of $H$ denoted by $G_0, G_1, \dots$ and $G'_0, G'_1, \dots$ such that for every $i \geq 0$ it holds that:

1. $G_i$ and $G'_i$ are isomorphic to $G$,
2. $G_i$ does not contain $u_i$ and $v_i$,
3. $G'_i$ does not contain $u_i$ and $v_{i+1}$,
4. $u_i$ and $u_{i+1}$ is determined by $v_i$, and,
5. $v_i$ and $v_{i+1}$ is determined by $u_i$.

Put $G_0 = G$ and $A = V_H \setminus V_G$. Consider the spanning tree $S$ given by Lemma 2.3.1. Because no vertex of $A$ can be removed to obtain an induced point-determining subgraph, it follows that every vertex must have a corresponding edge in $S$. Consequently the number of edges of $S$ is at least $|A|$. Because $G$ itself is point-determining, it follows that every edge of $S$ must contain at least one vertex of $A$. These two conditions yields to the pair of vertices $v_0 \in A = V_H \setminus V_G$ and $v_1 \in V_G$ connected by an edge in $S$ and consequently we have a vertex $u_0 \in A$ which determines them. We have obtained $G_0, u_0, v_0, v_1$ with the desired properties. This finishes the initial step of the induction.

At the induction step assume we have constructed $G_i, u_i, v_i, v_{i+1}$. We show the construction of $G'_i$ and $u_{i+1}$. We consider two cases. If $v_{i+1} \notin V_{G_i}$ we put $G'_i = G_i$. If $v_{i+1} \in V_{G_i}$ we let $G'_i$ to be the graph induced by $H$ on $(V_{G_i} \setminus \{v_{i+1}\}) \cup \{v_i\}$. Because the neighbourhood of $v_i$ and $v_{i+1}$ differs only by a vertex $u_i \notin G_i$ which determines them we know that $G'_i$ is isomorphic to $G_i$ (and thus also to $G$) and moreover that $u_i$ is not a vertex of $G'_i$ (because $u_i \notin V_{G_i}$ can not determine itself and thus $u_i \neq v_i$). If $H$ was point-determining after removal of $v_{i+1}$ we would obtain a contradiction similarly as before. We can thus assume that $v_{i+1}$ determines at least one pair of vertices. Because neighbourhood $v_{i+1}$ and $v_i$ differs only by $u_i$ we know that one vertex of this pair is $u_i$. Denote by $u_{i+1}$ the second vertex.

Given $G'_i, u_i, u_{i+1}, v_{i+1}$ we proceed similarly. If $u_{i+1} \notin V_{G'_i}$ we put $G_{i+1} = G'_i$. If $u_{i+1} \in V_{G'_i}$ we let $G_{i+1}$ to be the graph induced by $H$ on $(V_{G'_i} \setminus \{u_{i+1}\}) \cup \{u_i\}$. Again $G_{i+1}$ is isomorphic to $G$ and does not contain $u_{i+1}$ nor $v_{i+1}$. Denote by $v_{i+2}$ a vertex determined by $u_{i+1}$ from $v_{i+1}$ (which again must exist by our assumption) and we have obtained $G_{i+1}, u_{i+1}, v_{i+1}, v_{i+2}$ with the desired properties. This finishes the inductive step of the construction.

Because $H$ is finite, we know that both sequences $u_0, u_1, \dots$ and $v_0, v_1, \dots$ contains repeated vertices. Without loss of generality we can assume that repeated vertex with lowest index $j$ appears in the first sequence. We thus have $u_j = u_i$ for some $i < j$. By minimality of $j$ we can assume that $v_i, v_{i+1}, \dots v_{j-1}$ are all unique. Assume that $v_i$ is in the neighbourhood of $u_i$, then $v_i$ is not in the neighbourhood of $u_{i+1}$ (because it determines this pair) and consequently also $u_{i+1}, u_{i+2}, \dots, u_j$. A contradiction with $u_j = u_i$. If $v_i$ is not in the neighbourhood of $u_i$ we proceed analogously. □

## 2.4 Generalised dualities always exist

To demonstrate the usefulness of Theorem 2.1.1 and Propositions 2.2.1 and 2.2.2 give a simple proof of the existence of generalised dualities in full homomorphism order. For two finite sets of graphs $\mathcal{F}$ and $\mathcal{D}$ we say that $(\mathcal{F}, \mathcal{D})$ is a *generalised finite F-duality pair* (sometimes also *$\mathcal{D}$-obstruction*) if for any graph $G$ there exists $F \in \mathcal{F}$ such that $F \xrightarrow{F} G$ if and only if $G \xrightarrow{F} D$ for no $D \in \mathcal{D}$.

Existence of (generalised) dualities have several consequences. To mention one, it implies that the decision problem "given graph $G$ is there $D \in \mathcal{D}$ and full homomorphism $G \to D$?" is polynomial time solvable for every fixed finite family $\mathcal{D}$ of finite graphs. In the graph homomorphism order the dualities (characterised in [62]) are rare. In the case of full homomorphisms they are however always guaranteed to exist.

**Theorem 2.4.1** ([83, 22, 33, 3])**.** *For every finite set of graphs $\mathcal{D}$ there is a finite set of graphs $\mathcal{F}$ such that $(\mathcal{F}, \mathcal{D})$ is a generalised finite F-duality pair.*

*Proof.* Without loss of generality assume that $\mathcal{D}$ is a non-empty set of F-cores. Consider set $\mathcal{X}$ of all F-cores $G$ such that there is $D \in \mathcal{D}$, $G \to D$. Because, by Proposition 2.2.2, the number of vertices of every such $G$ is bounded from above by the number of vertices of $D$ and because $\mathcal{D}$ is finite, we know that $\mathcal{X}$ is finite.

Now denote by $\mathcal{F}$ the set of all F-cores $H$ such that $H \notin \mathcal{X}$ and there is $G \in \mathcal{X}$ such that $(G, H)$ is a gap. By Theorem 2.1.1 this set is finite. We show that $(\mathcal{F}, \mathcal{D})$ is a duality pair.

Consider an F-core $G$, either $G \in \mathcal{X}$ and thus there is $D \in \mathcal{D}$, $G \to D$ or $G \notin \mathcal{X}$ and then consider a sequence of $F$-cores $G_1, G_2, \ldots, G_{|G|} = G$ such that $G_1 \in \mathcal{X}$ consists of single vertex, $G_{i+1}$ is created from $G_i$ by adding a single vertex for every $1 \leq i < |G|$ (such sequence exists by Theorem 2.3.2). Clearly there is $1 \leq j < |G|$ such that $G_j \in \mathcal{X}$ and $G_{j+1} \notin \mathcal{X}$. Because $(G_j, G_{j+1})$ forms a gap, we know that $G_{j+1} \in \mathcal{F}$. □

**Remark 2.4.1.** A stronger result is shown by Feder and Hell [22, Theorem 3.1] who shows that if $\mathcal{D}$ consists of single graph $G$ with $k$ vertices, then $\mathcal{F}$ can be chosen in a way so it contains graphs with at most $k + 1$ vertices and there are at most two graphs having precisely $k + 1$ vertices. While, by Theorem 2.1.1, we can also give the same upper bound on number of vertices of graphs in $\mathcal{F}$, it does not really follow that there are at most two graphs needed. It appears that the full machinery of [22] is necessary to prove this result.

In the opposite direction it does not seem to be possible to derive Theorem 2.1.1 from this characterisation of dualities, because given pair of non-isomorphic F-cores $G \xrightarrow{F} H$ and $\mathcal{D}$ a full homomorphism dual of $\{G\}$ it does not hold that for a graph $F \in \mathcal{D}$ such that $D \xrightarrow{F} H$ there is also full homomorphism $G \xrightarrow{F} H$.

## 2.5 Full homomorphisms of relational structures

To the date, the full homomorphism order has been analysed in the context of graphs and digraphs only. Let us introduce generalised setting of relational structures:

A language $L$ is a set of relational symbols $R \in L$, each associated with natural number $a(R)$ called *arity.* A *(relational) $L$-structure* $\mathbf{A}$ is a pair $(A, (R_\mathbf{A}; R \in L))$ where $R_\mathbf{A} \subseteq A^{a(R)}$ (i.e. $R_\mathbf{A}$ is a $a(R)$-ary relation on $A$). The set $A$ is called the *vertex set* of $\mathbf{A}$ and elements of $A$ are *vertices.* The language is usually fixed and understood from the context. If the set $A$ is finite we call $\mathbf{A}$ *finite structure.* The class of all finite relational $L$-structures will be denoted by $\mathrm{Rel}(L)$.

A *homomorphism* $f : \mathbf{A} \to \mathbf{B} = (B, (R_\mathbf{B}; R \in L))$ is a mapping $f : A \to B$ satisfying for every $R \in L$ the implication $(x_1, x_2, \ldots, x_{a(R)}) \in R_\mathbf{A} \implies (f(x_1), f(x_2), \ldots, f(x_{a(R)})) \in R_\mathbf{B}$. A homomorphism is *full* if the above implication is equivalence, i.e. if for every $R \in L$ we have $(x_1, x_2, \ldots, x_{a(R)}) \in R_\mathbf{A} \iff (f(x_1), f(x_2), \ldots, f(x_{a(R)})) \in R_\mathbf{B}$.

Given structure $\mathbf{A}$ its vertex $v$ is contained in a *loop* if there exists $(v, v, \ldots, v) \in R_\mathbf{A}$ for some $R \in L$ of arity at least 2. Given relation $R_\mathbf{A}$ we denote by $\overline{R}_\mathbf{A}$ its complement, that is the set of all $a(R)$-tuples $\vec{t}$ of vertices of $A$ that are not in $R_\mathbf{A}$.

When considering full homomorphism order in this context, the first problem is what should be considered to be the neighbourhood of a vertex. This can be described as follows: Given $L$-structure $\mathbf{A}$, relation $R \in L$ and vertex $v \in A$ such that $(v, v, \ldots, v) \notin R_\mathbf{A}$ the *$R$-neighbourhood* of $v$ in $\mathbf{A}$, denoted by $N^R_\mathbf{A}(v)$ is the set of all tuples $\vec{t} \setminus v$ created from $\vec{t} \in R_\mathbf{A}$ containing $v$. Here by $\vec{t} \setminus v$ we denote tuple created from $\vec{t}$ by replacing all occurrences of vertex $v$ by a special symbol $\bullet$ which is not part of any vertex set. If $(v, v \ldots, v) \in R_\mathbf{A}$ then the $R$-neighbourhood $N^R_\mathbf{A}(v)$ is the set of all tuples $\vec{t} \setminus v$ created from $\vec{t} \in \overline{R}_\mathbf{A} \cup \{(v, v, \ldots, v)\}$. The *neighbourhood* of $v$ in $\mathbf{A}$ is a function assigning every relational symbol its neighbourhood: $N_\mathbf{A}(v)(R) = N^R_\mathbf{A}(v)$.

We say that $L$-structure $\mathbf{A}$ is *point-determining* if there are no two vertices with same neighbourhood. With these definitions direct analogies of Proposition 2.2.1 and 2.2.2 for $\mathrm{Rel}(L)$ follows.

Analogies of Lemma 2.3.1, Theorem 2.3.2 and Theorem 2.1.1 do not follow for relational structures in general. Consider, for example, a relational structure with three vertices $\{a, b, c\}$ and a single ternary relation $R$ containing one tuple $(a, b, c)$. Such structure is point-determining, but the only point-determining substructures consist of single vertex. There is however deeper problem with carrying Lemma 2.3.1 to relational structures: if a pair of vertices $u, u'$ is determined by vertex $v$ their neighbourhood may differ by tuples containing additional vertices. Thus the basic argument about cycles can not be directly applied here. We consequently formulate results for relational language consisting of unary and binary relations only (and, as a special case, to digraphs):

**Theorem 2.5.1.** *Let $L$ be a language containing relational symbols of arity at most 2. If $\mathbf{A}$ and $\mathbf{B}$ are (relational) F-cores and $(\mathbf{A}, \mathbf{B})$ is an F-gap, then $\mathbf{A}$ can be obtained* *from* $\mathbf{B}$ *by removal of one vertex.*

The example above shows that the limit on arity of relational symbols is actually necessary. This may be seen as a surprise, because the results about digraph homomorphism orders tend to generalise naturally to relational structures and we thus close this paper by an open problem of characterising gaps in full homomorphism order of relational structures in general.

# Chapter 3

# Universality of homomorphism order

It is a non-trivial result that every countable partial order can be found as a suborder of the homomorphism order of graphs. This has been first proved in the even stronger setting of category theory [65]. Subsequently, it has been shown that many restricted classes of graphs (such as oriented trees [49], oriented paths [48], partial orders and lattices [54]) admit this universality property.

We show a very simple and versatile argument based on divisibility which immediately yields the universality of the homomorphism order of directed graphs and discuss three applications.

## 3.1 Universal partial orders

In this section we give a construction of a universal partial order. Let us first review some basic concept and notations.

In the whole paper we consider only finite and countable partial orders. An *embedding* of a partial order $(Q, \leq_Q)$ in $(P, \leq_P)$ is a mapping $e : P \to Q$ satisfying $x \leq_P y$ if and only if $e(x) \leq_Q e(y)$. In such a case we also say that $(Q, \leq_Q)$ is a *suborder* of $(P, \leq_P)$.

For a given partial order $(P, \leq)$, the *down-set* $\downarrow x$ is $\{y \in P \mid y \leq x\}$.

Any finite partial order $(P, \leq_P)$ can be represented by finite sets ordered by the inclusion, e.g. when $x$ is represented by $\downarrow x$. This is a valid embedding, because $\downarrow x \subseteq \downarrow y$ if and only if $x \leq_P y$.

Without loss of generality we may assume that $P$ a subset of some fixed countable set $A$, e.g. $\mathbb{N}$. Consequently, the partial order formed by the system $P_{\mathrm{fin}}(A)$ of all finite subsets of $A$ ordered by the inclusion contains any finite partial order as a suborder.

Such orders are called are *finite-universal.* We reserve the term *universal* for orders that contain every countable partial order as a suborder.

Finite-universal and universal orders may be viewed as countable orders of rich structure — they are of infinite dimension, and that they contain finite chains, antichains and decreasing chains of arbitrary length. While finite-universal partial orders are rather easy to construct, e.g., as the disjoint union of all finite partial orders, the existence of a universal partial order can be seen as a counter-intuitive fact: there are uncountably many countable partial orders, yet all of them can be "packed" into a single countable structure.

The universal partial order can be build in two steps. For these we need further terminology: An order is *past-finite*, if every down-set is finite. An order is *past-finite-universal* if it contains every past-finite order. Analogously, *future-finite* and *future-finite-universal* orders are defined w.r.t. finiteness of up-sets.

**1.** Observe that the mapping $e(x) = \downarrow x$ is also an embedding $e : (P, \leq) \to (P_{\text{fin}}(A), \subseteq)$ in the case when $(P, \leq)$ is past-finite and $P \subseteq A$. Since a past-finite partial order turns to be future-finite when the direction of inequalities is reversed, we get:

**Proposition 3.1.1.** *For any countably infinite set $A$ it holds that*

*(i) the order $(P_{\text{fin}}(A), \subseteq)$ is past-finite-universal, and*

*(ii) the order $(P_{\text{fin}}(A), \supseteq)$ is future-finite-universal.*

**2.** For a given partial order $(Q, \leq)$ we construct the *subset order*, $(P_{\text{fin}}(Q), \leq_Q^{\text{dom}})$, where

$$X \leq_Q^{\text{dom}} Y \iff \text{for every } x \in X \text{ there exists } y \in Y \text{ such that } x \leq y.$$

We show that the subset order is universal:

**Theorem 3.1.2.** *For every future-finite-universal partial order $(F, \leq_F)$ it holds that $(P_{\text{fin}}(F), \leq_F^{\text{dom}})$ is universal.*

*Proof.* Proof (sketch). It is easy to check that $(P_{\text{fin}}(F), \leq_F^{\text{dom}})$ is indeed partial order. We sketch the way to embed any given partial order in $(P_{\text{fin}}(F), \leq_F^{\text{dom}})$. Let be given any countable partial order $(P, \leq_P)$. Without loss of generality we may assume that $P \subseteq \mathbb{N}$. This way we enforce a linear order $\leq$ on the elements of $P$. The order $\leq$ is unrelated to the partial order $\leq_P$. We decomposed $(P, \leq_P)$ into:

1. The *forward order* $\leq_f$, where $x \leq_f y$ if and only if $x \leq_P y$ and $x \leq y$, and

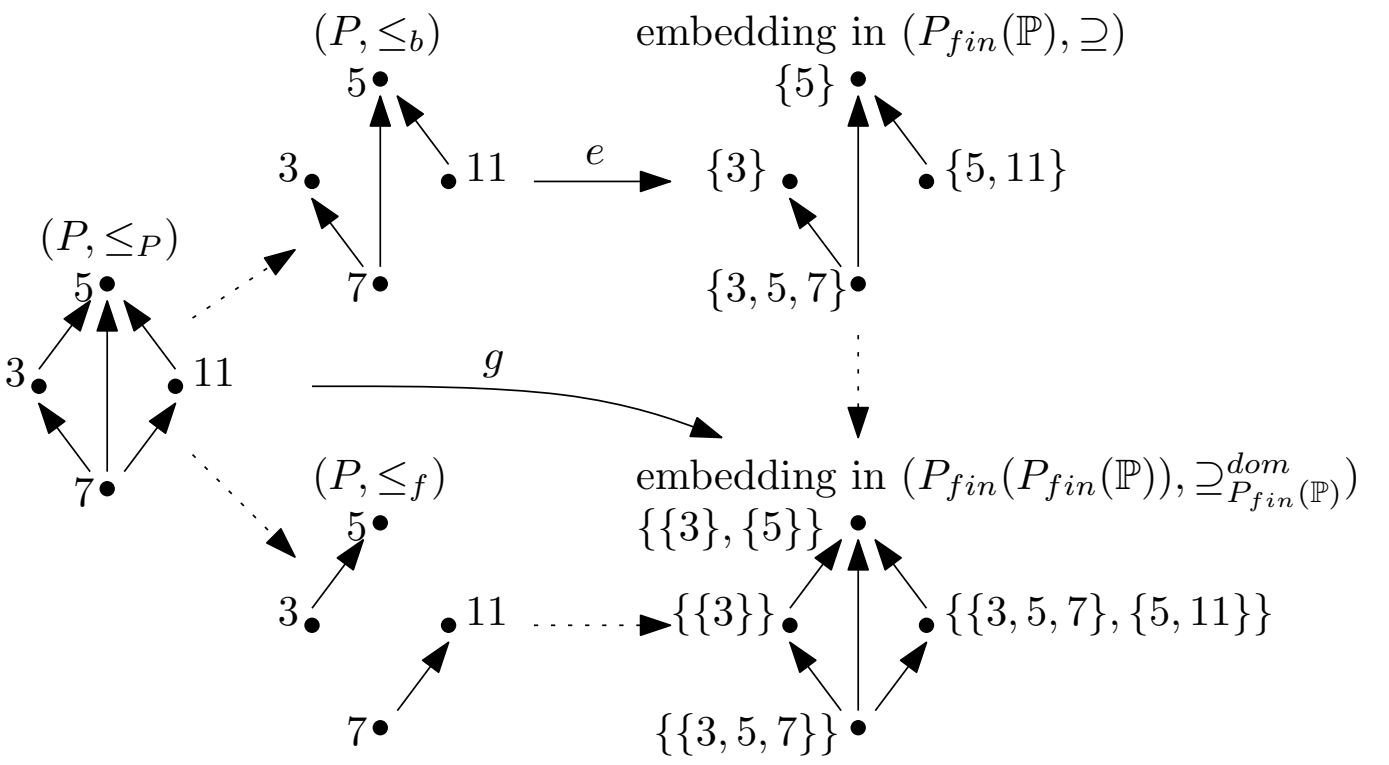

Figure 3.1: A representation of $(P, \leq_P)$ according to Theorem 3.1.2

2. the *backward order* $\leq_b$, where $x \leq_b y$ if and only if $x \leq_P y$ and $x \geq y$.

For every $x \in P$ both sets $\{y \mid y \leq_f x\}$ and $\{y \mid x \leq_b y\}$ are finite. In other words $(P, \leq_f)$ is past-finite and $(P, \leq_b)$ is future-finite.

Since $(F, \leq_F)$ is future-finite-universal, there is an embedding $e : (P, \leq_b) \to (F, \leq_F)$. For every $x \in P$ we now define: $g(x) = \{e(y) \mid y \leq_f x\}$.

□

An example of this construction is depicted in Figure 3.1. We chose $F$ to be set of prime numbers for reasons that will become clear shortly.

We remark that the embedding $g$ constructed in the proof of Theorem 3.1.2 has the property that $g(x)$ depends only on elements $y < x$. Such embeddings are known as on-line embeddings because they can be constructed inductively without a-priori knowledge of the whole partial order. See also [48, 49, 46] for similar constructions.

By Proposition 3.1.1 we see that a particular example of a past-finite-universal order is $(P_{\text{fin}}(\mathbb{P}), \subseteq)$, where $\mathbb{P}$ is the class of all odd prime numbers. It follows that $(P_{\text{fin}}(\mathbb{P}), \supseteq)$ is future-finite-universal. As for $X, Y \in P_{\text{fin}}(\mathbb{P})$ holds that $X \subseteq Y$ if and only if $\prod X$ divides $\prod Y$, we immediately obtain a special embeddings of the subset orders by divisibility as:

**Proposition 3.1.3.**

*a) The divisibility order* $(\mathbb{N}, |)$ *is past-finite-universal,*

*b) the reversed divisibility order* $(\mathbb{N}, \overleftarrow{|})$ *is future-finite-universal,*

*c) the* subset reverse divisibility order $(P_{\text{fin}}(\mathbb{N}), \overleftarrow{|}^{dom}_{\mathbb{N}})$ *is universal.*

In the following we show that the subset reverse divisibility order can be directly represented in the homomorphism order.

## 3.2 The homomorphism order

For given directed graphs $G$ and $H$ a *homomorphism* $f : G \to H$ is a mapping $f : V_G \to V_H$ such that $(u, v) \in V_G$ implies $(f(u), f(v)) \in V_H$. He existence of homomorphism $f : G \to H$ is traditionally denoted by $G\,\mathrm{Hom}_R(H,)$. This allows us to consider the existence of a homomorphism, $\mathrm{Hom}_R(,,)$ to be a (binary) relation on the class of directed graphs.

The relation $\to$ is reflexive (the identity is a homomorphism) and transitive (a composition of two homomorphisms is still a homomorphism). Thus the existence of a homomorphism induces a quasi-order on the class of all finite directed graphs. We denote the quasi-order induced by the existence of homomorphisms on directed graph by $(\mathsf{DiGraphs}, \leq)$ and on undirected graphs by $(\mathsf{Graphs}, \leq)$. When speaking of orders, we use $G \leq H$ in the same sense as $G\,\mathrm{Hom}_R(H,)$. These quasi-orders can be easily transformed into a partial order by choosing a particular representative for each equivalence class. In the case of graph homomorphism such representative is up to isomorphism unique vertex minimal element of each class, the *graph core*.

Both homomorphism orders $(\mathsf{DiGraphs}, \leq)$ and $(\mathsf{Graphs}, \leq)$ have been extensively studied and proved to be fruitful areas of research, see [34].

The original argument for universality of partial order [65] used complex graphs and ad-hoc constructions. It thus came as a surprise that the homomorphism order is universal even on the class of oriented paths [46]. While oriented paths is a very simple class of graphs, the universality argument remained rather complex. We can show show the universality of another restricted class easily.

Let $\overrightarrow{C}_k$ stand for the directed cycle on $k$ vertices with edges oriented in the same direction; $\mathsf{DiCycle}$ is the class of directed graphs formed by all $\overrightarrow{C}_k$, $k \geq 3$; and $\mathsf{DiCycles}$ is the class of directed graphs formed by disjoint union of finitely many graphs in $\mathsf{DiCycle}$.

**Theorem 3.2.1.** *The partial order* $(\mathsf{DiCycles}, \leq)$ *is universal.*

*Proof.* As $\overrightarrow{C}_k \to \overrightarrow{C}_l$ if and only if $k \overleftarrow{|} l$, we get the conclusion directly from Proposition 3.1.3. □

## 3.3 Applications

### 3.3.1 The fractal property of the homomorphism order

As a strengtening of the universality of homomorphism order we can show that every non-trivial interval in the order is in universal. This property under name of fractal property was first shown by Nešetřil [60] but the proof was difficult and never published. Easier proof based on the divisibility argument can be found in Chapter 3 and [47].

### 3.3.2 Universality of order induced by locally injective homomorphisms

Graph homomorphisms are just one of many mappings between graphs that induce a partial order. Monomorphisms, embeddings or full homomorphisms also induce partial orders. The homomorphism order however stands out as especially interesting and the universality result is one of unique properties of it. Other orders fails to be universal for rather trivial reasons, such as lack of infinite increasing or decreasing chains. A notable exception is the graph minor order, that is known to not be universal as a consequence of celebrated result of Robertson and Seymour [67]. We consider the following order:

A homomorphism $f : G \to H$ is *locally injective*, if for every vertex $v$ the restriction of the mapping $f$ to the domain $N_G(v)$ and range $N_H(f(v))$ is injective. (Here $N_G(v)$ denote the open neighborhood of a vertex). This order was first studied by Fiala, Paulusma and Telle in [26] where the degree refinement matrices are used to describe several interesting properties. We can further show:

**Theorem 3.3.1.** *The class of all finite connected graphs ordered by the existence of locally injective homomorphisms is universal.*

The proof of this theorem is based on a simple observation that every homomorphism between directed cycles is also locally injective homomorphism. The universality of locally injective homomorphism order on DiCycles thus follows from Theorem 3.2.1. This is a key difference between Theorem 3.2.1 and the universality of oriented paths: homomorphisms between oriented paths require flipping that can not be easily interpreted by locally injective homomorphisms.

In the second part of proof of Theorem 3.3.1 the cycles need to be connected together into a single connected graph in a way preserving all homomorphisms intended. This argument is technical and will appear in [25].

### 3.3.3 Universality of homomorphism order of line graphs

We close the paper by yet another application answering question of Roberson [66] asking about the universality of homomorphism order on the class of linegraphs of graphs with a vertices of degree at most $d$. We were able to give an affirmative answer.

**Theorem 3.3.2** ([23]). *The homomorphism order of line graphs of regular graphs with maximal degree $d$ is universal for every $d \geq 3$.*

This result may seem counter-intuitive with respect to the Vizing theorem. Vizing class 1 contains the graphs whose chromatic index is the same as the maximal degree of a vertex, while Vizing class 2 contains the remaining graphs. Because the Vizing class 1 is trivial it may seem that the homomorphism order on the Vizing class 2 should be simple, too. The converse is true.

# Chapter 4

# Fractal property of the graph homomorphism order

## 4.1 Introduction

In this note we consider finite simple graphs and countable partial orders. On these graphs we consider all homomorphisms between them. Recall that for graphs $G = (V_G, E_G)$ and $H = (V_H, E_H)$ a *homomorphism* $f : G \to H$ is an edge preserving mapping $f : V_G \to V_H$, that is:

$$\{x, y\} \in E_G \implies \{f(x), f(y)\} \in E_H.$$

If there exists a homomorphism from graph $G$ to $H$, we write $G \to H$.

Denote by $\mathscr{C}$ the class of all finite simple undirected graphs without loops and multiple edges, and by $\leq$ the following order:

$$G \to H \iff G \leq H.$$

$(\mathscr{C}, \leq)$ is called the *homomorphism order*.

The relation $\leq$ is clearly a quasiorder which becomes a partial order when factorized by homomorphism equivalent graphs. This homomorphism equivalence takes particularly simple form, when we represent each class by the so called core. Here, a *core* of a graph is its minimal homomorphism equivalent subgraph. It is well known that up to an isomorphism every equivalence class contains a unique core [34]. However, for our purposes it is irrelevant whether we consider $(\mathscr{C}, \leq)$ as a quasiorder or a partial order. For brevity we speak of the homomorphism order in both cases.

The homomorphism order has some special properties, two of which are expressed as follows:

**Theorem 4.1.1.** $(\mathscr{C}, \leq)$ *is (countably) universal.*

*Explicitly: For every countable partial order $P$ there exists an embedding of $P$ into $(\mathscr{C}, \leq)$.*

Here an *embedding* of partial order $(P, \leq)$ to partial order $(P', \leq')$ is an injective function $f : P \to P'$ such that for every $u, v \in P$, $u \leq v'$ if and only if $f(u) \leq' f(v)$.

**Theorem 4.1.2.** $(\mathscr{C}, \leq)$ *is dense.*

*Explicitly: For every pair of graphs $G_1 < G_2$ there exists $H$ such that $G_1 < H < G_2$. This holds with the single exception of $K_1 < K_2$, which forms the only gap of the homomorphism order of undirected graphs.*

As usual, $K_n$ denotes the complete graph with $n$ vertices. We follow the standard graph terminology as e.g. [34]. As the main result of this paper we complement these structural results by the following statement:

**Theorem 4.1.3.** $(\mathscr{C}, \leq)$ *has the fractal property.*

*Explicitly: For every pair $G_1 < G_2$, distinct from $K_1$ and $K_2$ (i.e. the pair $(G_1, G_2)$ is not a gap), there exists an order-preserving embedding $\Phi$ of $\mathscr{C}$ into the interval*

$$[G_1, G_2]_{\mathscr{C}} = \{H; G_1 < H < G_2\}.$$

Putting otherwise, every nonempty interval in $\mathscr{C}$ contains a copy of $\mathscr{C}$ itself.

Theorem 4.1.1 was proved first in [32, 65] and reproved in [48, 49]. Theorem 4.1.2 was proved in [81] and particularly simpler proof was given by Perles and Nešetřil [60], see also [62, 34].

Theorem 4.1.3 was formulated in [60] and remained unpublished since. The principal ingredient of the proof is the Sparse Incomparability Lemma [61]. In addition, we give yet another proof of Theorem 4.1.3. In fact, we prove all three Theorems 4.1.1, 4.1.2 and 4.1.3 (Theorem 4.1.2 is a corollary of Theorem 4.1.3).

First, to make the paper self-contained we also give in Section 4.2 a short and easy proof of universality of $(\mathscr{C}, \leq)$ which was developed in [23] and sketched in [24]. Then, in Section 4.3 we give first proof of Theorem 4.1.3 based on the Sparse Incomparability Lemma [61, 34]. Then in Section 4.4 we prove a strenghtening of Theorem 4.1.2 (stated as Lemma 4.4.1). This will be needed for our second proof of Theorem 4.1.3 which is flexible enough for applications. Thus this paper summarizes perhaps surprisingly easy proofs of theorems which originally had difficult proofs.

## 4.2 Construction of a universal order

### 4.2.1 Particular universal partial order

Let $(\mathcal{P}, \leq_{\mathcal{P}})$ be a partial order, where $\mathcal{P}$ consists of all finite sets of odd integers, and where for $A, B \in \mathcal{P}$ we put $A \leq_{\mathcal{P}} B$ if and only if for every $a \in A$ there is $b \in B$ such that $b$ divides $a$. We make use of the following:

**Theorem 4.2.1** ([23]). *The order $(\mathcal{P}, \leq_{\mathcal{P}})$ is a universal partial order.*

To make the paper self-contained we give a brief proof of this assertion. (See also [32, 48, 46] for related constructions of universal partial orders.) The proof of Theorem 4.2.1 follows from two simple lemmas.

We say that a countable partial order is *past-finite* if every down-set $x^{\downarrow} = \{y; y \leq x\}$ is finite. A countable partial order is *past-finite-universal*, if it contains every past-finite partial order as a suborder. *Future-finite* and *future-finite-universal* orders are defined analogously with respect to up-sets $x^{\uparrow} = \{y; y \geq x\}$.

Let $P_f(X)$ denote the set of all finite subsets of $X$. The following lemma extends a well known fact about representing finite partial orders by sets ordered by the subset relation.

**Lemma 4.2.2.** *For any countably infinite set $X$, the partial order $(P_f(X), \subseteq)$ is past-finite-universal.*

*Proof.* Consider an arbitrary past-finite order $(Q, \leq_Q)$. Without loss of generality we may assume that $Q \subseteq X$. Let $\Phi$ be the mapping that assigns every $x \in Q$ its down-set, i.e. $\Phi(x) = \{y \in Q; y \leq x\}$. It is easy to verify that $\Phi$ is indeed an embedding $(Q, \leq_Q) \to (P_f(X), \subseteq)$. □

By the *divisibility partial order*, denoted by $(\mathbb{N}, \leq_d)$, we mean the partial order on positive integers, where $n \leq_d m$ if and only if $n$ is divisible by $m$. Denote by $\mathbb{Z}_o$ the set of all odd integers $n$, $n \geq 3$.

**Lemma 4.2.3.** *The divisibility partial order $(\mathbb{Z}_o, \leq_d)$ is future-finite-universal.*

*Proof.* Denote by $\mathbb{P}$ the set of all odd prime numbers. Apply Lemma 4.2.2 for $X = \mathbb{P}$. Observe that $A \in P_f(\mathbb{P})$ is a subset of $B \in P_f(\mathbb{P})$ if and only if $\prod_{p \in A} p$ divides $\prod_{p \in B} p$. □

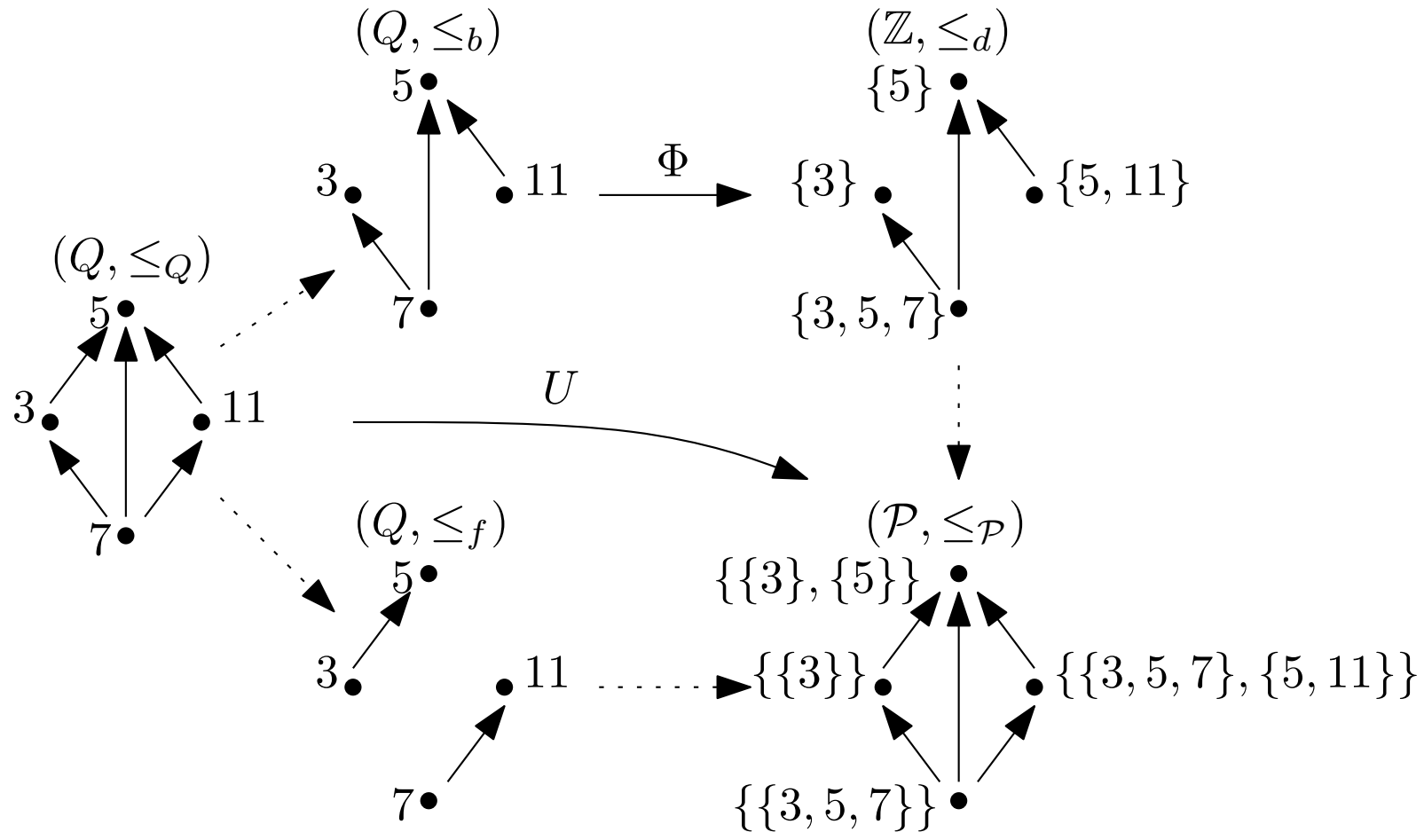

Figure 4.1: Example of the construction of an embedding $U : (Q, \leq_Q) \to (\mathcal{P}, \leq_\mathcal{P})$.

*Proof of Theorem 4.2.1.* Let $(Q, \leq_Q)$ be a given partial order. Without loss of generality we may assume that $Q$ is a subset of $\mathbb{P}$. This way we obtain also the usual linear order $\leq$ (i.e. the comparison by the size) on the elements of $Q$. In the following construction, the order $\leq$ determines in which sequence the elements of $Q$ are processed (it could be also interpreted as the time of creation of the elements of $Q$).

We define two new orders on $Q$: the *forwarding order* $\leq_f$ and the *backwarding order* $\leq_b$ as follows:

1. We put $x \leq_f y$ if and only if $x \leq_Q y$ and $x \leq y$.

2. We put $x \leq_b y$ if and only if $x \leq_Q y$ and $x \geq y$.

Thus the partial order $(Q, \leq_Q)$ has been decomposed into $(Q, \leq_f)$ and $(Q, \leq_b)$. For every vertex $x \in Q$ both sets $\{y; y \leq_f x\}$ and $\{y; x \leq_b y\}$ are finite. It follows that $(Q, \leq_f)$ is past-finite and that $(Q, \leq_b)$ is future-finite.

Since $(\mathbb{Z}_o, \leq_d)$ is future-finite-universal (Lemma 4.2.3), there is an embedding $\Phi : (Q, \leq_b) \to (\mathbb{Z}_o, \leq_d)$. The desired embedding $U : (Q, \leq_Q) \to (\mathcal{P}, \leq_\mathcal{P})$ is obtained by representing each $x \in Q$ by a set system $U(x)$ defined by (see Figure 4.1):

$$U(x) = \{\Phi(y); y \leq_f x\}.$$

To argue the correctness we first show that $U(x) \leq_\mathcal{P} U(y)$ implies $x \leq_Q y$. From the definition of $\leq_\mathcal{P}$ and the fact that $\Phi(x) \in U(x)$ follows that at least one $w \in Q$ exists, such that $\Phi(w) \in U(y)$ and $\Phi(x) \leq_d \Phi(w)$. By the definition of $U$, $\Phi(w) \in U(y)$ if and only if $w \leq_f y$. By the definition of $\Phi$, $\Phi(x) \leq_d \Phi(w)$ if and only if $x \leq_b w$. It follows that $x \leq_b w \leq_f y$ and thus also $x \leq_Q w \leq_Q y$ and consequently $x \leq_Q y$.

To show that $x \leq_Q y$ implies $U(x) \leq_\mathcal{P} U(y)$ we consider two cases.

1. When $x \le y$ then $U(x) \subseteq U(y)$ and thus also $U(x) \le_{\mathcal{P}} U(y)$.

2. Assume $x > y$ and take any $w \in Q$ such that $\Phi(w) \in U(x)$. From the definition of $U(x)$ we have $w \le_f x$. Since $x \le_Q y$ we have $w \le_Q y$. If $w \le y$, then $w \le_f y$ and we have $\Phi(w) \in U(y)$. In the other case if $w > y$ then $w \le_b y$ and thus $\Phi(w) \le_d \Phi(y)$. Because the choice of $w$ is arbitrary, it follows that $U(x) \le_{\mathcal{P}} U(y)$.

□

Clearly, as in e.g. [48] this can be interpreted as Alice-Bob game played on finite partial orders. Alice always wins.

### 4.2.2 Representing divisibility

Denote by $\overrightarrow{C_p}$ the directed cycle of length $p$, i.e. the graph $(\{0, 1, \ldots, p-1\}, \{(i, i+1); i = 0, 1, \ldots, p-1\})$, where addition is performed modulo $p$. Denote by $\mathcal{D}$ the class of disjoint unions of directed cycles.

**Theorem 4.2.4.** *The homomorphism order $(\mathcal{D}, \le)$ is universal.*

*Proof.* Observe first that a homomorphism $f : \overrightarrow{C_p} \to \overrightarrow{C_q}$ between two cycles $\overrightarrow{C_p}$ and $\overrightarrow{C_q}$ exists if and only if $q$ divides $p$.

Consequently, for two collections of disjoint cycles $\sum_{p \in A} \overrightarrow{C_p}$ and $\sum_{q \in B} \overrightarrow{C_q}$ a homomorphism

$$f : \sum_{p \in A} \overrightarrow{C_p} \to \sum_{q \in B} \overrightarrow{C_q}$$

exists if and only if

$$A \le_{\mathcal{P}} B,$$

with respect to the universal partial order $(\mathcal{P}, \le_{\mathcal{P}})$ of Theorem 4.2.1.

Since we have used odd primes in the proof of Lemma 4.2.3, the minimum of each set in $\mathcal{P}$ is at least three. Hence, each $A \in \mathcal{P}$ corresponds to a disjoint union of odd cycles. □

**Remark 4.2.1.** Denote by $\overrightarrow{\mathscr{C}}$ the class of all directed graphs. Theorem 4.2.4 yields immediately that the order $(\overrightarrow{\mathscr{C}}, \le)$ is also universal.

The class of disjoint union of directed odd cycles is probably the simplest class for which the homomorphism order is universal. However note that here the key property is that objects are not connected and contains odd cycles of unbounded length. If we want to obtain connected graphs with bounded cycles then we have to refer to [48, 46, 49] where it is proved that that the class of finite oriented trees $\mathcal{T}$ and even the class of finite oriented paths form universal partial orders. These strong notions are not needed in this paper. However note that from our results here it also follows that not only the class of planar graphs but also the class of outer-planar graphs form a universal partial order.

## 4.3 The fractal property

To prove Theorem 4.1.3, we use the following result proved by Nešetřil and Rödl by non-constructive methods [61]. Later, non-trivial constructions were given by Matoušek-Nešetřil [56] and Kun [51]:

**Theorem 4.3.1** (Sparse Incomparability Lemma [61], see e.g. Theorem 3.12 of [34]). *Let $l$ be positive integer. For any non-bipartite graphs $G_1$ and $G_2$ with $G_1 < G_2$, there exists a connected graph $F$ such that*

- *$F$ is (homomorphism) incomparable with $G_1$ (i.e. $F \not\leq G_1$ and $G_1 \not\leq F$);*
- *$F < G_2$; and*
- *$F$ has girth at least $l$, where the girth of a graph is the length of its shortest cycle.*

In the sequel we combine the Sparse Incomparability Lemma and the universality of $(\mathcal{D}, \leq)$, together with the standard indicator technique developed by Hedrlín and Pultr [31, 34].

The essential construction of this method takes an oriented graph $G$ and a graph $I$ with two distinguished vertices $a$, $b$ and creates a graph $G * (I, a, b)$, obtained by substituting every arc $(x, y)$ of $G$ by a copy of the graph $I$, where $x$ is identified with $a$ and $y$ is identified with $b$, see Figure 4.2 for an example.

*Proof of Theorem 4.1.3.* Let be given undirected graphs $G_1 < G_2$ not forming a gap. By our assumptions $G_2$ is a core distinct from $K_2$ as otherwise $G_1 = K_1$. (If $G_2 = K_2$ then $K_1 = K_1$ and we have a gap). We may also assume without loss of generality that $G_1$ is not bipartite since in such a case we may replace $G_1$ by a graph $G_1 < G_1' < G_2$ (given by Theorem 4.1.2), making the interval $[G_1', G_2]_{\mathscr{C}}$ even narrower than $[G_1, G_2]_{\mathscr{C}}$. Because all bipartite graphs are homomorphically equivalent it follows that $G_1'$ is non-bipartite.

Let $l \geq 5$ be any odd integer s.t. the longest odd cycle of $G_2$ has length at most $l$.

For the indicator we use the graph $I_l$, depicted in Figure 4.3. The graph $I_l$ can be viewed either as a subdivision of $K_4$, or as 3 cycles of length $l + 2$ amalgamated together. The indicator $I_l$ is rigid, i.e. the only homomorphism $I_l \to I_l$ is the identity [34, Proposition 4.6]. Note also that $I_l$ allows a homomorphism to the undirected cycle of length $l$, as is also depicted in Figure 4.3.

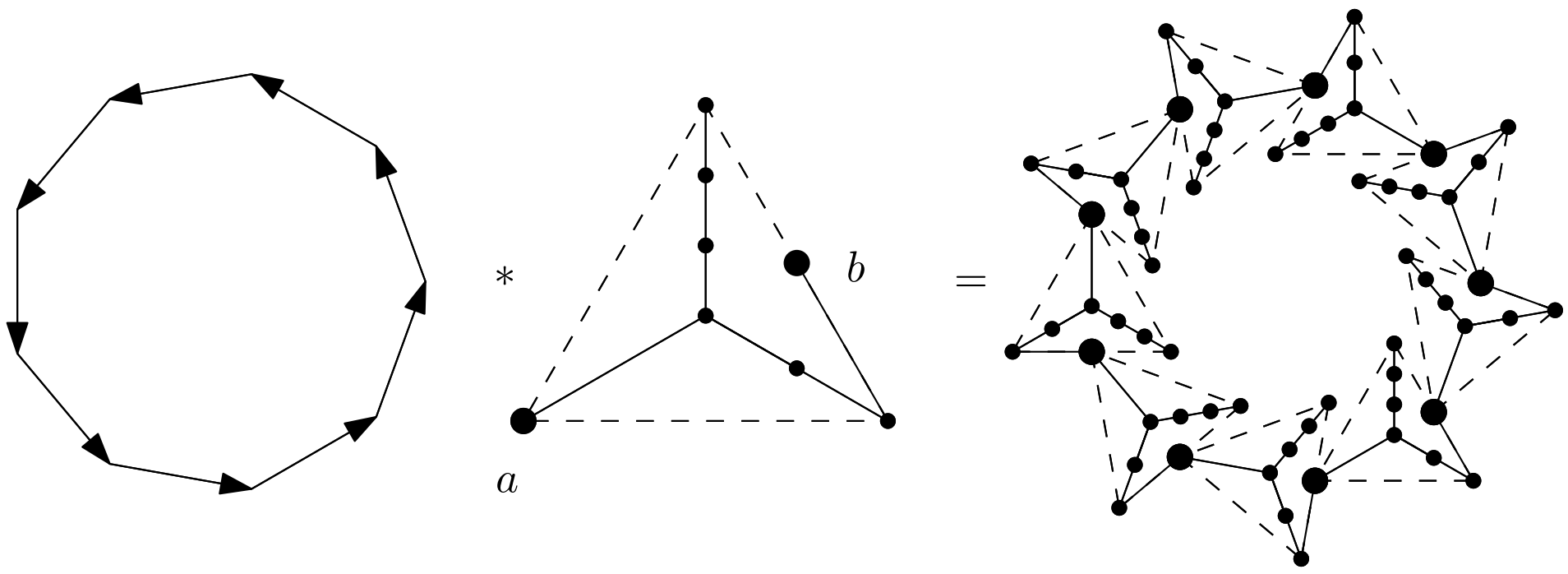

Figure 4.2: Construction of $\overrightarrow{C_p} * (I, a, b)$ .

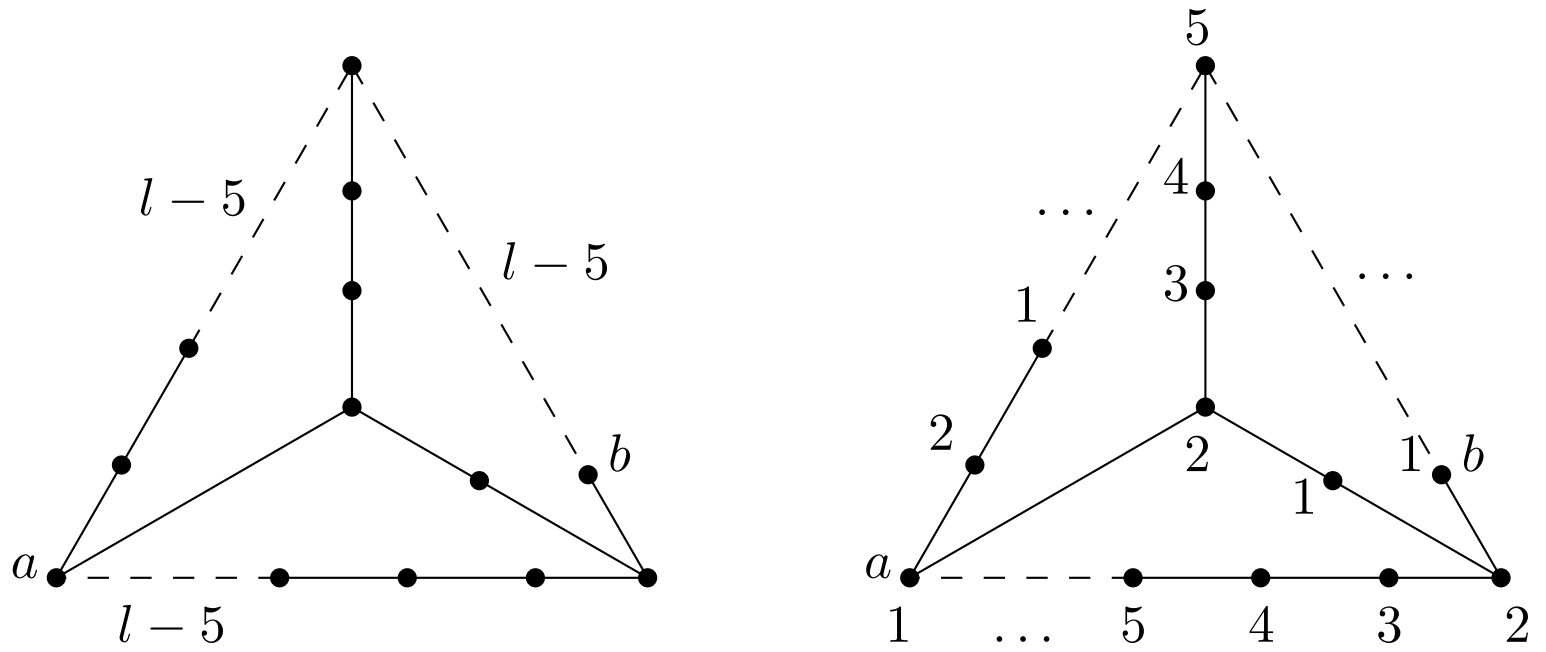

Figure 4.3: The indicator $I_l$ and its homomorphism to $C_l$. Dashed lines represent paths with $l-5$ internal vertices.

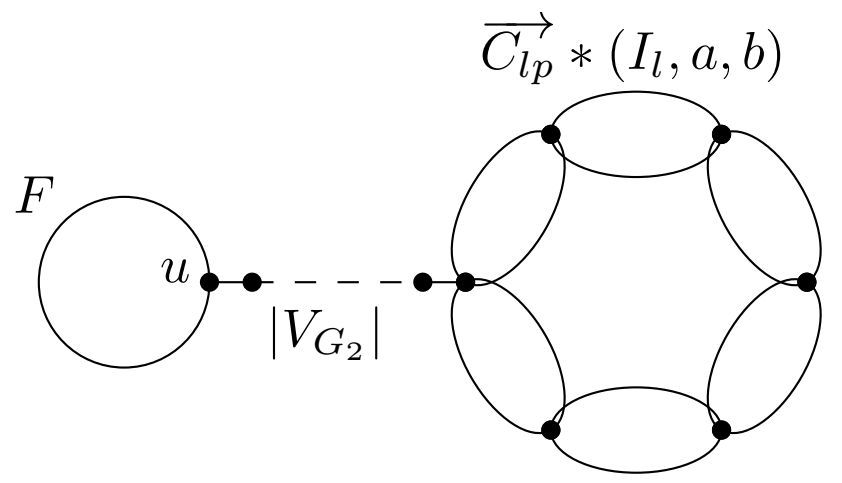

Figure 4.4: Graph $H_p$.

We continue with the construction of a graph $H_A$ from a set of odd integers $A \in \mathcal{P}$. Let $F$ be a connected graph satisfying the conclusions of the Sparse Incomparability Lemma. We fix an arbitrary vertex $u$ of $F$.

Then, given a positive integer $p \geq 3$, we apply the indicator $I_l, a, b$ on the directed cycle of length $l \cdot p$ to obtain $\overrightarrow{C_{lp}} * (I_l, a, b)$. (Observe that $\overrightarrow{C_{lp}} * (I_l, a, b) \to \overrightarrow{C_{lq}} * (I_l, a, b)$ if and only if $q$ divides $p$.) We then join any vertex of the original cycle $\overrightarrow{C_{lp}}$ to $u$ by a path of length $|V_{G_2}|$, see Figure 4.4.

Observe that the resulting graph $H_p$ allows a homomorphism to $G_2$, since:

1. There exists a homomorphism $f : F \to G_2$ by Theorem 4.3.1;

2. the indicator $I_l$ has a homomorphism to a cycle of length $l$, which can be simply transformed to a homomorphism $g$ to any odd cycle of length $l' \leq l$ in $G_2$ (by the choice of $l$);

3. the mapping $g$ could be chosen s.t. $g(a) = g(b)$, hence it can be extended to all vertices of $\overrightarrow{C_p}$;

4. the distance between the image of $u$ and the cycle of length $l'$ is at most $|V_{G_2}|$, therefore both homomorphisms $f$ and $g$ can be combined together and extended to the whole graph $H_p$ straightforwardly.

To complete the construction of $H_A$, we put

$$H_A = \sum_{p \in A} H_p + G_1.$$

The conclusion of Theorem 4.1.3 follows from the following three properties:

1. For every $A \in \mathcal{P} : G_1 < H_A$.

   The $\leq$ inequality is obvious, since $G_1$ is a component of each $H_A$.

   Since $F$ is a subgraph of each $H_p$, by Theorem 4.3.1 there is no homomorphism $H_A \to G_1$ whenever $A$ is nonempty.

2. For every $A \in \mathcal{P} : H_A < G_2$.

   The existence of homomorphisms $H_p \to G_2$ and $G_1 \to G_2$ yields a homomorphism $H_A \to G_2$.

   As $G_2 \not\leq F$, and as the shortest cycle in $\overrightarrow{C_{lp}} * (I_l, a, b)$ has length $l + 2$, which is by the choice of $l$ longer than the length of any odd cycle in $G_2$, there is no homomorphism $G_2 \to H_A$.

3. For every $A, B \in \mathcal{P} : H_A \to H_B$ if and only if $A \leq_{\mathcal{P}} B$.

   It is easy to see that $q$ divides $p$ iff $\overrightarrow{C_{lp}} * (I_l, a, b) \to \overrightarrow{C_{lq}} * (I_l, a, b)$. This is a standard argument. Note that the paths between $F$ and $\overrightarrow{C_{lp}}$ in $H_p$, and between $F$ and $\overrightarrow{C_{lq}}$ in $H_q$ have the same length and the

   vertex $u$ of attachment has been chosen in the same way in both cases. Therefore, $H_p \to H_q$ and consequently, $A \leq_{\mathcal{P}} B$ implies $H_A \to H_B$.

   Assume now that $H_A \to H_B$. We have already excluded $H_p \to G_1$, hence by the connectivity of each $H_p$ neccessarily follows that $\sum_{p \in A} H_p \to \sum_{q \in B} H_q$. This in turns leads to $\sum_{p \in A} \overrightarrow{C_{lp}} * (I_l, a, b) \to \sum_{q \in B} \overrightarrow{C_{lq}} * (I_l, a, b)$ which is equivalent to $A \leq_{\mathcal{P}} B$.

These three properties guarantee that gave a full embedding of $(\mathcal{D}, \leq)$ into $(\mathscr{C}, \leq)$, which maps every $\sum_{p \in I} \overrightarrow{C_p}$ into the interval $[G, G']_{\mathscr{C}}$ in $\mathscr{C}$. By Theorem 4.2.4 the universality of $[G, G']_{\mathscr{C}}$ follows. □

## 4.4 Alternative elementary proof

The sparse incomparability lemma holds for "dense" classes of graphs. Here we establish the fractal property of graphs by a different technique which allows us to prove the fractal property of some "sparse" classes of graphs as well. For example we can reduce the (stronger form of) density for planar graphs to the fractal property of the class of planar graphs. But first we formulate the proof for general graphs. We shall make use of the following two assertions.

**Lemma 4.4.1.** *Given graphs $G_1 < G_2$, $G_2$ non-bipartite, there exists connected graphs $H_1$ and $H_2$ with properties:*

1. *$H_1$ and $H_2$ are homomorphically incomparable, and,*
2. *$G_1 < H_i < G_2$ for $i = 1, 2$.*

*In other words, any non-gap interval in the homomorphism order contains two incomparable graphs.*

*Proof.* Proof is a variant of the Nešetřil-Perles' proof of density [34]. Put

$$H_1 = (G_2 \times H) + G_1,$$

where $H$ is a graph that we specify later, $+$ is the disjoint union and $\times$ denotes the direct product.

Then obviously $G_1 \leq H_1 \leq G_2$. If the odd girth of $G_2$ is larger than the odd girth of $G_1$ then $G_2 \not\to H_1$. If the chromatic number $\chi(H) > |V(G_1)|^{|V(G_2)|}$ then any homomorphism $G_2 \times H \to G_1$ induces a homomorphism $G_2 \to G_1$ which is absurd (see [34, Theorem 3.20]). Thus $G_2 \times H \not\to G_1$ and

$$G_1 < H_1 < G_2.$$

We choose $H$ to be a graph with large odd girth and chromatic number (known to exist [21]). This finishes construction of $H_1$. (Note that here we use the fact that the odd girth of the product is the maximum of the odd girths of its factors.) Now we repeat the same argument with the pair $G_1$ and $G_2 \times H$, put

$$H_2 = (G_2 \times H') + G_1.$$

If the odd girth of $H'$ is larger than the odd girth of $G_2 \times H$ then $G_2 \times H \not\to H_2$ (assuming $H$ and thus $G_2 \times H$ is connected). Thus in turn $H_1 \not\to H_2$. If $\chi(H') > \max(|V(G_2 \times H)|^{|V(G_2)|}, |V(G_2)|^{|V(G_1)|})$ then again $G_2 \times H' \to G_2 \times H$ implies $G_2 \to G_2 \times H$ which is absurd. Similarly $G_2 \times H' \to G_1$ implies $G_2 \to G_1$ and thus $G_2 \times H' \not\to G_2$

We may also assume that graphs $H_1$ and $H_2$ from Lemma 4.4.1 are connected as otherwise we can join components by a long enough path. Connectivity also follows from the following folklore fact:

**Claim 4.4.2.** *For every connected non-bipartite graph $H$ there exists an integer $l$ such that for any two vertices $x, y \in V(H)$ and any $l' \geq l$ there exists a homomorphism $f : P_{l'} \to H$ such that $f(0) = x$ and $f(l') = y$. ($P_{l'}$ is the path of length $l'$ with vertices $\{0, 1, \ldots, l'\}$).*

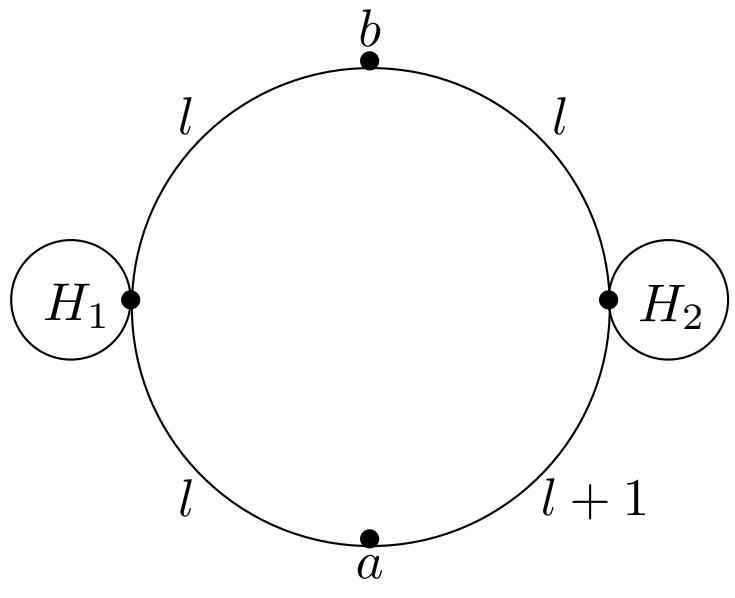

Figure 4.5: Gadget $I$.

This concludes the construction of $H_1$ and $H_2$. □

*Second proof of Theorem 4.1.3.* Let $G_1 < G_2$ be a non-gap, thus $G_2$ is non-bipartite. Assume without loss of generality that $G_2$ is connected. Since the homomorphism order is universal, we prove the universality of the interval $(G_1, G_2)$ by embedding the homomorphism order into it.

Let $H_1$, $H_2$ be two connected graphs given by Lemma 4.4.1. We may assume that both $H_1$ and $H_2$ are cores. Let $l$ be the number given by Claim 4.4.2 for graph $G_2$. We may assume $l > \max\{|V(H_1)|, |V(H_2)|, |V(G_2)|\}$. We construct the gadget $I$ consisting of graphs $H_1$, $H_2$ joined together by two paths of length $2l$ and $2l+1$. We choose two distinguished vertices $a$, $b$ to be the middle vertices of these two paths, see Figure 4.5.

We observe that any homomorphism $f : I \to I$ is surjective (because both $G_1$ and $G_2$ are cores), it is also an identity on vertices of $I \setminus (G_1 \cup G_2)$ and that there exists a homomorphism $f : I \to G_2$ such that $f(a) = f(b)$.

For every oriented graph $G$ define graph $\Phi(G)$ as $\Phi(G) = G * (I, a, b)$. We know $G_1 < \Phi(G) < G_2$ (as any homomorphisms $f_i : H_i \to G_2$, $i \in \{1, 2\}$, can be extended to a homomorphism $\Phi(G) \to G_2$).

We finish the proof by proving $\Phi(G) \to \Phi(G')$ if and only if $G \to G'$.

Assume first that there exists a homomorphism $f : G \to G'$. Consider the function $g$ defined as $f$ on vertices of $G$ and as the unique mapping which maps a copy of $(I, a, b)$ in $G$ corresponding to edge $(u, v)$ to the copy of $(I, a, b)$ in $G'$ corresponding to edge $(f(u), f(v))$. Hence $g$ is a homomorphism.

Let now $g$ be a homomorphism $\Phi(G) \to \Phi(G')$. By the girth assumption and connectivity of $H_1$ and $H_2$ we know that $g$ maps every copy of $H_1$ (or $H_2$)) in $\Phi(G)$ to a copy of $H_1$ (or $H_2$) in $\Phi(G')$. Again, by the girth argument it follows that every copy of the indicator $(I, a, b)$ in $G$ is mapped to a copy of the indicator $(I, a, b)$ in $G'$. But the only copies of $(I, a, b)$ in both $G$ and $G'$ are those corresponding to the edges of $G$ and $G'$.

Since $I$ is a core, it follows that any pair of vertices $(a, b)$ in a copy of $(I, a, b)$ has to be mapped to the vertices $(a, b)$ in any other copy of $(I, a, b)$. As copies of $(I, a, b)$ and hence also the pairs $(a, b)$ correspond to edges of $G'$, it follows that $g$ induces a mapping $f : V(G) \to V(G')$, which is a homomorphism $G \to G'$. This argument concludes the second proof of Theorem 4.1.3. $\square$

**Remark 4.4.1.** Note that in this second proof we have an one-to-one correspondence between homomorphisms $G \to G$ and $\Phi(G) \to \Phi(G')$.

## 4.5 Concluding remarks

**1.** Gaps on oriented graphs and on more general relational structures are more sophisticated. They were characterized by Nešetřil and Tardif [62]. In the same paper, a nice 1-1 correspondence between gaps and dualities has been shown. Consequently, the full discussion of fractal property of relational structures is more complicated and it will appear elsewhere.
**2.** The whole paper deals with finite graphs but there is no difficulty in generalizing our results to infinite graphs.
**3.** An interesting question (already considered in [62]) is: which intervals induce *isomorphic orders*. We provide neither a characterization nor a conjecture in this direction.
**4.** It seems that any natural universal class $\mathcal{K}$ of structures posesses the *gap-universal dichotomy*: An interval $[A, B]_{\mathcal{K}}$ in $\mathcal{K}$ either contains a gap or it contains a copy of $\mathcal{K}$ itself. While in general this fails, it remains to be seen, whether this is true for some general class of structures.
**5.** There is a great difference in treating universality, density and fractal property. Whereas the universality was established in many classes of partial orders (and categories), the density and fractal property was only established for just a few basic classes. Perhaps this should be investigated in greater depth. Apart from general relational structures (which we hope to treat in another paper) another interesting case is provided by structures with two equivalences (here the universality is established by [63]).

# Chapter 5

# Inferring Phylogenetic Trees from the Knowledge of Rare Evolutionary Events

Rare events have played an increasing role in molecular phylogenetics as potentially homoplasy-poor characters. In this contribution we analyze the phylogenetic information content from a combinatorial point of view by considering the binary relation on the set of taxa defined by the existence of a single event separating two taxa. We show that the graph-representation of this relation must be a tree. Moreover, we characterize completely the relationship between the tree of such relations and the underlying phylogenetic tree. With directed operations such as tandem-duplication-random-loss events in mind we demonstrate how non-symmetric information constrains the position of the root in the partially reconstructed phylogeny.

## 5.1 Introduction

Shared derived characters (synapomorphies or "Hennigian markers") that are unique to specific clades form the basis of classical cladistics [44]. In the context of molecular phylogenetics *rare genomic changes (RGCs)* can play the same important role [69, 7]. RGCs correspond to rare mutational events that are very unlikely to occur multiple times and thus are (almost) free of homoplasy. A wide variety of processes and associated markers have been proposed and investigated. Well-studied RGCs include presence/absence patterns of protein-coding genes [20] as well as microRNAs [71], retroposon integrations [74], insertions and deletions (indels) of introns [68], pairs of mutually exclusive introns (NIPs) [50], protein domains [15, 84], RNA secondary structures [59], protein fusions [76], changes in gene order [70, 8, 53], metabolic networks [27, 28, 57], transcription factor binding sites [64], insertions and deletions of arbitrary sequences [75, 2, 17], and variations of the genetic code [1]. RGCs clearly have proved to be phylogenetically informative and helped to resolve many of the phylogenetic questions where sequence data lead to conflicting or equivocal results.

Not all RGCs behave like cladistic characters, however. While presence/absence characters are naturally stored in character matrices whose columns can vary independently, this is not the case e.g. for gene order characters. From a mathematical point of view, character-based parsimony analysis requires that the mutations have a product structure in which characters are identified with factors and character states can vary independently of each other [79]. This assumption is violated whenever changes in the states of two distinct characters do not commute. Gene order is, of course, the prime example on non-commutative events.

Three strategies have been pursued in such cases: (i) Most importantly, the analog of the parsimony approach is considered for a particular non-commutative model. For the genome rearrangements an elaborated theory has been developed that considers various types of operations on (usually signed) permutations. Already the computation of editing distances is non-trivial. An added difficulty is that the interplay of different operations such as reversals, transpositions, and tandem-duplication-random-loss (TDRL) events is difficult to handle [5, 30]. (ii) An alternative is to focus on distance-based methods [80]. Since good rate models are usually unavailable, distance measures usually are not additive and thus fail to precisely satisfy the assumptions underlying the most widely used methods such as neighbor joining. (iii) Finally, the non-commutative data structure can be converted into a presence-absence structure, e.g., by using pairwise adjacencies [78] as a representation of permutations or using list alignments in which rearrangements appear as pairs of insertions and deletions [29]. While this yields character matrices that can be fed into parsimony algorithms, these can only result in approximate heuristics.

While it tends to be difficult to disentangle multiple, super-imposed complex changes such as genome rearrangements or tandem duplication, it is considerably simpler to recognize whether two genes or genomes differ by a single RGC operation. It make sense therefore to ask just how much phylogenetic information can be extracted from elementary RGC events. Of course, we cannot expect that a single RGC will allow us to (re)construct a detailed phylogeny. It can, however, provide us with solid, well-founded constraints. Furthermore, we can hope that the combination of such constraints can be utilized as a practicable method for phylogenetic inference. Recently, we have shown that orthology assignments in large gene families imply triples that must be displayed by the underlying species tree [45, 37]. In a phylogenomics setting a sufficient number of such triple constraints can be collected to yield fully resolved phylogenetic trees [43], see [36] for an overview.

A plausible application scenario for our setting is the rearrangement of mitogenomes [70]. Since mitogenomes are readily and cheaply available, the taxon sampling is sufficiently dense so that the gene orders often differ by only a single rearrangement or not at all. These cases are identifiable with near certainty [5]. Moreover, some RGC are inherently directional. Probably the best known example is the tandem duplication random loss (TDRL) operation [13]. We will therefore also consider a directed variant of the problem.

In this contribution, we ask how much phylogenetic information can be retrieved from single RGCs. More precisely, we consider a scenario in which we can, for every pair of taxa distinguish, for a given type of RGC, whether $x$ and $y$ have the same genomic state, whether $x$ and $y$ differ by exactly one elementary change, or whether their distance is larger than a single operation. We formalize this problem in the following way. Given a relation $\sim$, there is a phylogenetic tree $T$ with an edge labeling $\lambda$ (marking the elementary events) such that $x \sim y$ if and only if the edge labeling along the unique path $\mathbb{P}(x, y)$ from $x$ to $y$ in $T$ has a certain prescribed property $\Pi$. After defining the necessary notation and preliminaries, we give a more formal definition of the general problem in section 5.3.

The graphs defined by path relations on a tree are closely related to *pairwise compatibility graphs* (PCGs). A graph $G = (V, E)$ is a PCG if there is a tree $T$ with leaf set $V$, a positive edge-weight function $w : E(T) \to \mathbb{R}^+$, and two nonnegative real numbers $d_{\min} \le d_{\max}$ such that there is an edge $uv \in E(G)$ if and only if $d_{\min} \le d_{T,w}(x, y) \le d_{\max}$, where $d_{T,w}(x, y)$ is the sum of the weights of the edges on the unique path $\mathbb{P}(x, y)$ in $T$. One writes $G = \mathrm{PCG}(T, w, d_{\min}, d_{\max})$. In this contribution we will primarily be interested in the special case where $\Pi$ is "a single event along the path". Although PCGs have been studied extensively, see e.g., [12, 86, 85, 11, 58, 19], the questions are different from our motivation and, to our knowledge, no results have been obtained that would simplify the characterization of the PCGs corresponding to the "single-1-relation" in Section 5.4. Furthermore, PCGs are always treated as undirected graphs in the literature. We also consider an antisymmetric (Section 5.5) and a general directed (Section 5.6) versions of the single-1-relation motivated by RGCs with directional information.

The main result of this contribution can be summarized as follows: (i) The graph of a single-1-relation is always a forest. (ii) If the single-1-relation is connected, there is a unique minimally resolved tree that explains the relation. The same holds true for the connected components of an arbitrary relation. (iii) Analogous results hold for the anti-symmetric and the mixed variants of the single-1-relation. In this case not only the tree topology but also the position of the root can be determined or at least constrained. Together, these results in a sense characterize the phylogenetic information contained in rare events: if the single-1-relation graph is connected, it is a tree that through a bijection corresponds to a uniquely defined, but not necessarily fully resolved, phylogenetic tree. Otherwise, it is forest whose connected components determine subtrees for which the rare events provide at least some phylogenetically relevant information.

## 5.2 Preliminaries

### 5.2.1 Basic Notation

We largely follow the notation and terminology of the book by [72]. Throughout, $X$ denotes always a finite set of at least three taxa. We will consider both undirected and directed graphs $G = (V, E)$ with finite vertex set $V(G) := V$ and edge set or arc set $E(G) := E$. For a digraph $G$ we write $\underline{G}$ for its *underlying undirected graph* where $V(G) = V(\underline{G})$ and $\{x, y\} \in E(\underline{G})$ if $(x, y) \in E(G)$ or $(y, x) \in E(G)$. Thus, $\underline{G}$ is obtained from $G$ by ignoring the direction of edges. For simplicity, edges $\{x, y\} \in E(G)$ (in the undirected case) and arcs $(x, y) \in E(G)$ (in the directed case) will be both denoted by $xy$.

The representation $G(R) = (V, E)$ of a relation $R \subseteq V \times V$ has vertex set $V$ and edge set $E = \{xy \mid (x, y) \in R\}$. If $R$ is irreflexive, then $G$ has no loops. If $R$ is symmetric, we regard $G(R)$ as an undirected graph. A *clique* is a complete subgraph that is maximal w.r.t. inclusion. An equivalence relation is *discrete* if all its equivalence classes consist of single vertices.

A tree $T = (V, E)$ is a connected cycle-free undirected graph. The vertices of degree 1 in a tree are called leaves, all other vertices of $T$ are called *inner vertices.* An edge of $T$ is *interior* if both of its end vertices are inner vertices, otherwise the edge is called *terminal.* For technical reasons, we call a vertex $v$ an inner vertex and leaf if $T$ is a single vertex graph $(\{v\}, \emptyset)$. However, if $T$ is an edge $vw$ we refer to $v$ and $w$ as leaves but not as inner vertices. Hence, in this case the edge $vw$ is not an interior edge

A *star* $S_m$ with $m$ leaves is a tree that has at most one inner vertex. A *path* $P_n$ (on $n$ vertices) is a tree with two leaves and $n - 3$ interior edges. There is a unique path $\mathbb{P}(x, y)$ connecting any two vertices $x$ and $y$ in a tree $T$. We write $e \in \mathbb{P}(x, y)$ if the edge $e$ connects two adjacent vertices along $\mathbb{P}(x, y)$. We say that a directed graph is a tree if its underlying undirected graph is a tree. A directed path $P$ is a tree on vertices $x_1, \ldots, x_n$ s.t. $x_i x_{i+1} \in E(P)$, $1 \leq i \leq n-1$. A graph is a forest if all its connected components are trees.

A tree is *rooted* if there is a distinguished vertex $\rho \in V$ called the *root*. Throughout this contribution we assume that the root is an inner vertex. Given a rooted tree $T = (V, E)$, there is a partial order $\preceq$ on $V$ defined as $v \preceq u$ if $u$ lies on the path from $v$ to the root. Obviously, the root is the unique maximal element w.r.t $\preceq$. For a non-empty subset of $W \subseteq V$, we define $lca(W)$, or the *least common ancestor of* $W$, to be the unique $\preceq_T$-minimal vertex of $T$ that is an ancestor of every vertex in $W$. In case $W = \{x, y\}$, we put $lca(x, y) := lca(\{x, y\})$. If $T$ is rooted, then by definition $lca(x, y)$ is a uniquely defined inner vertex along $\mathbb{P}(x, y)$.

We write $L(v)$ for the set of leaves in the subtree below a fixed vertex $v$, i.e., $L(v)$ is the set of all leaves for which $v$ is located on the unique path from $x \in L(v)$ to the root of $T$. The *children* of an inner vertex $v$ are its direct descendants, i.e., vertices $w$ with $vw \in E(T)$ s.t. that $w$ is further away from the root than $v$. A rooted or unrooted tree that has no vertices of degree two (except possibly the root of $T$) and leaf set $X$ is called a *phylogenetic tree $T$ (on $X$)*.

Suppose $X' \subseteq X$. A phylogenetic tree $T$ on $X$ *displays* a phylogenetic tree $T'$ on $X'$ if $T'$ can be obtained from $T$ by a series of vertex deletions, edge deletions, and suppression of vertices of degree 2 other than possibly the root, i.e., the replacement of an inner vertex $u$ and its two incident edges $e'$ and $e''$ by a single edge $e$, cf. Def. 6.1.2 in the book by [72]. In the rooted case, only a vertex between two incident edges may be suppressed; furthermore, if $X'$ is contained in a single subtree, then the $lca(X')$ becomes the root of $T'$. It is useful to note that $T'$ is displayed by $T$ if and only if it can be obtained from $T$ step-wisely by removing an arbitrarily selected leaf $y \in X \setminus X'$, its incident edge $e = yu$, and suppression of $u$ provided $u$ has degree 2 after removal of $e$.

We say that a rooted tree $T$ *contains* or *displays* the triple $\mathsf{xy|z}$ if $x, y$, and $z$ are leaves of $T$ and the path from $x$ to $y$ does not intersect the path from $z$ to the root of $T$. A set of triples $\mathcal{R}$ is consistent if there is a rooted tree that contains all triples in $\mathcal{R}$. For a given leaf set $L$, a triple set $\mathcal{R}$ is said to be *strict dense* if for any three distinct vertices $x, y, z \in L$ we have $|\{\mathsf{xy|z}, \mathsf{xz|y}, \mathsf{yz|x}\} \cap R| = 1$. It is well-known that any consistent strict-dense triple set $\mathcal{R}$ has a unique representation as a binary tree [43, Suppl. Material]. For a consistent set $R$ of rooted triples we write $R \vdash (\mathsf{xy|z})$ if any phylogenetic tree that displays all triples of $R$ also displays $(\mathsf{xy|z})$. [10] extend and generalized results by [16] and showed under which conditions it is possible to infer triples by using only subsets $R' \subseteq R$, i.e., under which conditions $R \vdash (\mathsf{xy|z}) \implies R' \vdash (\mathsf{xy|z})$ for some $R' \subseteq R$. In particular, we will use the following inference rules:

$$\{(\mathsf{ab|c}), (\mathsf{ad|c})\} \vdash (\mathsf{bd|c}) \tag{i}$$

$$\{(\mathsf{ab|c}), (\mathsf{ad|b})\} \vdash (\mathsf{bd|c}), (\mathsf{ad|c}) \tag{ii}$$

$$\{(\mathsf{ab|c}), (\mathsf{cd|b})\} \vdash (\mathsf{ab|d}), (\mathsf{cd|a}). \tag{iii}$$

## 5.3 Path Relations and Phylogenetic Trees

Let $\Lambda$ be a non-empty set. Throughout this contribution we consider a phylogenetic tree $T = (V, E)$ with edge-labeling $\lambda\colon E \to \Lambda$. An edge $e$ with label $\lambda(e) = k$ will be called a *k-edge*. We interpret $(T, \lambda)$ so that a RGC occurs along edge $e$ if and only if $\lambda(e) = 1$. Let $\Pi$ be a subset of the set of $\Lambda$-labeled paths. We interpret $\Pi$ as a property of the path and its labeling. The tree $(T, \lambda)$ and the property $\Pi$ together define a binary relation $\sim_\Pi$ on $X$ by setting

$$x \sim_\Pi y \qquad \Longleftrightarrow \qquad (\mathbb{P}(x, y), \lambda) \in \Pi \tag{i}$$

The relation $\sim_\Pi$ has the graph representation $G(\sim_\Pi)$ with vertex set $X$ and edges $xy \in E(G(\sim_\Pi))$ if and only if $x \sim_\Pi y$.

**Definition 5.3.1.** Let $(T, \lambda)$ be a $\Lambda$-labeled phylogenetic tree with leaf set $L(T)$ and let $G$ be a graph with vertex set $L(T)$. We say that $(T, \lambda)$ *explains* $G$ *(w.r.t. to the path property* $\Pi$*)* if $G = G(\sim_\Pi)$.

For simplicity we also say "$(T, \lambda)$ explains $\sim$" for the binary relation $\sim$.

We consider in this contribution the conceptually "inverse problem": Given a definition of the predicate $\Pi$ as a function of edge labels along a path and a graph $G$, is there an edge-labeled tree $(T, \lambda)$ that explains $G$? Furthermore, we ask for a characterization of the class of graph that can be explained by edge-labeled trees and a given predicate $\Pi$.

A straightforward biological interpretation of an edge labeling $\lambda : E \to \{0, 1\}$ is that a certain type of evolutionary event has occurred along $e$ if and only if $\lambda(e) = 1$. This suggests that in particular the following path properties and their associated relations on $X$ are of practical interest:

$x \overset{0}{\sim} y$ if and only if all edges in $\mathbb{P}(x, y)$ are labeled 0; For convenience we set $x \overset{0}{\sim} x$ for all $x \in X$.

$x \overset{1}{\sim} y$ if and only if all but one edges along $\mathbb{P}(x, y)$ are labeled 0 and exactly one edge is labeled 1;

$x \overset{1}{\rightharpoondown} y$ if and only if all edges along $\mathbb{P}(u,x)$ are labeled 0 and exactly one edge along $\mathbb{P}(u,y)$ is labeled 1, where $u = lca(x,y)$.

$x \overset{\geq k}{\sim} y$ with $k \geq 1$ if and only if at least $k$ edges along $\mathbb{P}(x,y)$ are labeled 1;

$x \rightsquigarrow y$ if all edges along $\mathbb{P}(u,x)$ are labeled 0 and there are one or more edges along $\mathbb{P}(u,y)$ with a non-zero label, where $u = lca(x,y)$.

We will call the relation $\overset{1}{\sim}$ the *single-1-relation*. It will be studied in detail in the following section. Its directed variant $\overset{1}{\rightharpoondown}$ will be investigated in Section 5.5. The more general relations $\overset{\geq k}{\sim}$ and $\rightsquigarrow$ will be studied in future work.

As noted in the introduction there is close relationship between the graphs of path relations introduced above and PCGs. For instance, the single-1-relations correspond to a graph of the form $G = \text{PCG}(T, \lambda, 1, 1)$ for some tree $T$. The exact-$k$ leaf power graph $\text{PCG}(T, \lambda, k, k)$ arise when $\lambda(e) = 1$ for all $e \in E(T)$ [9]. The "weight function" $\lambda$, however, may be 0 in our setting. It is not difficult to transform our weight functions to strictly positive values albeit at the expense of using less "beautiful" values of $d_{\min}$ and $d_{\max}$. The literature on the PCG, to our knowledge, does not provide results that would simplify our discussion below. Furthermore, the applications that we have in mind for future work are more naturally phrased in terms of Boolean labels, such as the "at least one 1" relation, or even vector-valued structures. We therefore do not pursue the relationship with PCGs further.

The combinations of labeling systems and path properties of primary interest to us have *nice properties*:

(L1) The label set $\Lambda$ is endowed with a semigroup $\boxplus : \Lambda \times \Lambda \to \Lambda$.

(L2) There is a subset $\Lambda_\Pi \subseteq \Lambda$ of labels such that $(\mathbb{P}(x,y), \lambda) \in \Pi$ if and only if $\lambda(\mathbb{P}(x,y)) := \boxplus_{e \in \mathbb{P}(x,y)} \lambda(e) \in \Lambda_\Pi$ or $\lambda(\mathbb{P}(x,y)) := \boxplus_{e \in \mathbb{P}(lca(x,y),y)} \lambda(e) \in \Lambda_\Pi$.

For instance, we may set $\Lambda = \mathbb{N}$ and use the usual addition for $\boxplus$. Then $\overset{0}{\sim}$ corresponds to $\Lambda_\Pi = \{0\}$, $\overset{1}{\sim}$ corresponds to $\Lambda_\Pi = \{1\}$, etc. The bounds $d_{\min}$ and $d_{\max}$ in the definition of PCGs of course is just a special case of of the predicate $\Lambda_\Pi$.

We now extend the concept of a phylogenetic tree displaying another one to the $\Lambda$-labeled case.

**Definition 5.3.2.** Let $(T, \lambda)$ and $(T', \lambda')$ be two phylogenetic trees with $L(T') \subseteq L(T)$. Then $(T, \lambda)$ displays $(T', \lambda')$ w.r.t. a path property $\Pi$ if (i) $T$ displays $T'$ and (ii) $(\mathbb{P}_T(x,y), \lambda) \in \Lambda_\Pi$ if and only if $(\mathbb{P}_{T'}(x,y), \lambda') \in \Lambda_\Pi$ for all $x, y \in L(T')$.

The definition is designed to ensure that the following property is satisfied:

**Lemma 5.3.3.** *Suppose* $(T, \lambda)$ *displays* $(T', \lambda')$ *and* $(T, \lambda)$ *explains a graph* $G$. *Then* $(T', \lambda')$ *explains the induced subgraph* $G[L(T')]$.

**Lemma 5.3.4.** *Let* $(T, \lambda)$ *display* $(T', \lambda')$. *Assume that the labeling system satisfies (L1) and (L2) and suppose* $\lambda'(e) = \lambda(e') \boxplus \lambda(e'')$ *whenever* $e$ *is the edge resulting from suppressing the inner vertex between* $e'$ *and* $e''$. *If* $T'$ *is displayed by* $T$ *then* $(T', \lambda')$ *is displayed by* $(T, \lambda)$ *(w.r.t. any path property* $\Pi$*).*

*Proof.* Suppose $T'$ is obtained from $T$ by removing a single leaf $w$. By construction $T'$ is displayed by $T$ and $\lambda(\mathbb{P}(x, y))$ is preserved upon removal of $w$ and suppression of its neighbor. Thus (L2) implies that $(T, \lambda)$ displays $(T', \lambda')$. For an arbitrary $T'$ displayed by $T$ this argument can be repeated for each individual leaf removal on the editing path from $T$ to $T'$. □

We note in passing that this construction is also well behaved for PCGs: it preserves path length, and thus distances between leaves, by summing up the weights of edges whenever a vertex of degree 2 between them is omitted.

Let us now turn to the properties of the specific relations that are of interest in this contribution.

**Lemma 5.3.5.** *The relation* $\overset{0}{\sim}$ *is an equivalence relation.*

*Proof.* By construction, $\overset{0}{\sim}$ is symmetric and reflexive. To establish transitivity, suppose $x \overset{0}{\sim} y$ and $y \overset{0}{\sim} z$, i.e., $\lambda(e) = 0$ for all $e \in \mathbb{P}(x, y) \cup \mathbb{P}(y, z)$. By uniqueness of the path connecting vertices in a tree, $\mathbb{P}(x, z) \subseteq \mathbb{P}(x, y) \cup \mathbb{P}(y, z)$, i.e., $\lambda(e) = 0$ for all $e \in \mathbb{P}(x, z)$ and therefore $x \overset{0}{\sim} z$. □

Since $\overset{0}{\sim}$ is an equivalence relation, the graph $G(\overset{0}{\sim})$ is a disjoint union of complete graphs, or in other words, each connected component of $G(\overset{0}{\sim})$ is a clique.

We are interested here in characterizing the pairs of trees and labeling functions $(T, \lambda)$ that explain a given relation $\rho$ as its $\overset{0}{\sim}$, $\overset{1}{\sim}$ or $\overset{1}{\rightharpoonup}$ relation. More precisely, we are interested in the least resolved trees with this property.

**Definition 5.3.6.** Let $(T, \lambda)$ be an edge-labeled phylogenetic tree with leaf set $X = L(T)$. We say that $(T', \lambda')$ is *edge-contracted from* $(T, \lambda)$ if the following conditions hold: (i) $T' = T/e$ is the usual graph-theoretical edge contraction for some interior edge $e = \{u, v\}$ of $T$.

(ii) The labels satisfy $\lambda'(e') = \lambda(e')$ for all $e' \neq e$.

Note that we do not allow the contraction of terminal edges, i.e., of edges incident with leaves.

**Definition 5.3.7** (Least and Minimally Resolved Trees)**.** Let $R \in \{\overset{0}{\sim}, \overset{1}{\sim}, \overset{1}{\rightharpoonup}, \overset{1}{\sim} / \overset{0}{\sim}, \overset{1}{\rightharpoonup} / \overset{0}{\sim}\}$. A pair $(T, \lambda)$ is *least resolved* for a prescribed relation $R$ if no edge contraction leads to a tree $T', \lambda')$ of $(T, \lambda)$ that explains $R$. A pair $(T, \lambda)$ is *minimally resolved* for a prescribed relation $R$ if it has the fewest number of vertices among all trees that explain $R$.

Note that every minimally resolved tree is also least resolved, but not *vice versa.*

## 5.4 The single-1-relation

The single-1-relation does not convey any information on the location of the root and the corresponding partial order on the tree. We therefore regard $T$ as unrooted in this section.

**Lemma 5.4.1.** *Let $(T, \lambda)$ be an edge-labeled phylogenetic tree with leaf set $X$ and resulting relations $\overset{0}{\sim}$ and $\overset{1}{\sim}$ over $X$. Assume that $A, B$ are distinct cliques in $G(\overset{0}{\sim})$ and suppose $x \overset{1}{\sim} y$ where $x \in A$ and $y \in B$. Then $x' \overset{1}{\sim} y'$ holds for all $x' \in A$ and $y' \in B$.*

*Proof.* First, observe that $\mathbb{P}(x', y') \subseteq \mathbb{P}(x', x) \cup \mathbb{P}(x, y) \cup \mathbb{P}(y, y')$ in $T$. Moreover, $\mathbb{P}(x', x)$ and $\mathbb{P}(y, y')$ have only edges with label 0. As $\mathbb{P}(x, y)$ contains exactly one non-0-label, thus $\mathbb{P}(x', y')$ contains at most one non-0-label. If there was no non-0-label, then $\mathbb{P}(x, y) \subseteq \mathbb{P}(x, x') \cup \mathbb{P}(x', y') \cup \mathbb{P}(y', y)$ would imply that $\mathbb{P}(x, y)$ also has only 0-labels, a contradiction. Therefore $x' \overset{1}{\sim} y'$. $\square$

As a consequence it suffices to study the single-1-relation on the quotient graph $G(\overset{1}{\sim})/\overset{0}{\sim}$. To be more precise, $G(\overset{1}{\sim})/\overset{0}{\sim}$ has as vertex set the equivalence classes of $\overset{0}{\sim}$ and two vertices $c_i$ and $c_j$ are connected by an edge if there are vertices $x \in c_i$ and $y \in c_j$ with $x \overset{1}{\sim} y$. Analogously, the graph $G(\overset{1}{\rightharpoonup})/\overset{0}{\sim}$ is defined.

For a given $(T, \lambda)$ and its corresponding relation $\overset{0}{\sim}$ consider an arbitrary nontrivial equivalence class $c_i$ of $\overset{0}{\sim}$. Since $\overset{0}{\sim}$ is an equivalence relation, the induced subtree $T'$ with leaf set $c_i$ and inner vertices $lca(c)$ for any subset $c \subseteq c_i$ contains only 0-edges and is maximal w.r.t. this property. Hence, we could remove $T'$ from $T$ and identify the root $lca(c_i)$ of $T'$ in $T$ by a representative of $c_i$, while keeping the information of $\overset{0}{\sim}$ and $\overset{1}{\sim}$. Let us be a bit more explicit about this point. Consider trees $(T_Y, \lambda_Y)$ displayed by $(T, \lambda)$ with leaf sets $Y$ such that $Y$ contains exactly one (arbitrarily chosen) representative from each $\overset{0}{\sim}$ equivalence class of $(T, \lambda)$. For any such trees $(T_Y, \lambda_Y)$ and $(T'_Y, \lambda'_Y)$ with the latter property, there is an isomorphism $\alpha : T_Y \to T_{Y'}$ such that $\alpha(y) \overset{0}{\sim} y$ and $\lambda_{Y'}(\alpha(e)) = \lambda_{Y'}(e)$. Thus, all such $(T_Y, \lambda_Y)$ are isomorphic and differ basically only in the choice of the particular representatives of the equivalence classes of $\overset{0}{\sim}$. Furthermore, $T_Y$ is isomorphic to the quotient graph $T/\overset{0}{\sim}$ obtained from $T$ by replacing the (maximal) subtrees where all edges are labeled with 0 by a representative of the corresponding $\overset{0}{\sim}$-class. Suppose $(T, \lambda)$ explains $G$. Then $(T_Y, \lambda_Y)$ explains $G[Y]$ for a given $Y$. Since all $(T_Y, \lambda_Y)$ are isomorphic, all $G[Y]$ are also isomorphic, and thus $G[Y] = G/\overset{0}{\sim}$ for all $Y$.

To avoid unnecessarily clumsy language we will say that "$(T, \lambda)$ explains $G(\overset{1}{\sim})/\overset{0}{\sim}$" instead of the more accurate wording "$(T, \lambda)$ displays $(T_Y, \lambda_Y)$ where $Y$ contains exactly one representative of each $\overset{0}{\sim}$ equivalence class such that $(T_Y, \lambda_Y)$ explains $G(\overset{1}{\sim})/\overset{0}{\sim}$".

In contrast to $\overset{0}{\sim}$, the single-1-relation $\overset{1}{\sim}$ is not transitive. As an example, consider the star $S_3$ with leaf set $\{x, y, z\}$, inner vertex $v$, and edge labeling $\lambda(v, x) = \lambda(v, z) = 1 \neq \lambda(v, y) = 0$. Hence $x \overset{1}{\sim} y$, $y \overset{1}{\sim} z$ and $x \overset{1}{\not\sim} z$. In fact, a stronger property holds that forms the basis for understanding the single-1-relation:

**Lemma 5.4.2.** *If $x \overset{1}{\sim} y$ and $x \overset{1}{\sim} z$, then $y \overset{1}{\not\sim} z$.*

*Proof.* Uniqueness of paths in $T$ implies that there is a unique inner vertex $u$ in $T$ such that $\mathbb{P}(x, y) = \mathbb{P}(x, u) \cup \mathbb{P}(u, y)$, $\mathbb{P}(x, z) = \mathbb{P}(x, u) \cup \mathbb{P}(u, z)$, $\mathbb{P}(y, z) = \mathbb{P}(y, u) \cup \mathbb{P}(u, z)$. By assumption, each of the three sub-paths $\mathbb{P}(x, u)$, $\mathbb{P}(y, u)$, and $\mathbb{P}(z, u)$ contains at most one 1-label. There are only two cases: (i) There is a 1-edge in $\mathbb{P}(x, u)$. Then neither $\mathbb{P}(y, u)$ nor $\mathbb{P}(z, u)$ may have another 1-edge, and thus $y \overset{0}{\sim} z$, which implies that $y \overset{1}{\not\sim} z$. (ii) There is no 1-edge in $\mathbb{P}(x, u)$. Then both $\mathbb{P}(y, u)$ and $\mathbb{P}(z, u)$ must have exactly one 1-edge. Thus $\mathbb{P}(y, z)$ harbors exactly two 1-edges, whence $y \overset{1}{\not\sim} z$. □

Lemma 5.4.2 can be generalized as follows.

**Lemma 5.4.3.** *Let $x_1, \dots, x_n$ be vertices s.t. $x_i \overset{1}{\sim} x_{i+1}$, $1 \leq i \leq n-1$. Then, for all $i, j$, $x_i \overset{1}{\sim} x_j$ if and only if $|i - j| = 1$.*

*Proof.* For $n = 3$, we can apply Lemma 5.4.2. Assume the assumption is true for all $n < K$. Now let $n = K$. Hence, for all vertices $x_i, x_j$ along the paths from $x_1$ to $x_{K-1}$, as well as the paths from $x_2$ to $x_K$ it holds that $|i-j| = 1$ if and only if we have $x_i \overset{1}{\sim} x_j$. Thus, for the vertices $x_i, x_j$ we have $|i-j| > 1$ if and only if we have $x_i \overset{1}{\not\sim} x_j$. Therefore, it remains to show that $x_1 \overset{1}{\not\sim} x_n$.

Assume for contradiction, that $x_1 \overset{1}{\sim} x_n$. Uniqueness of paths on $T$ implies that there is a unique inner vertex $u$ in $T$ that lies on all three paths $\mathbb{P}(x_1, x_2)$, $\mathbb{P}(x_1, x_n)$, and $\mathbb{P}(x_2, x_n)$.

There are two cases, either there is a 1-edge in $\mathbb{P}(x_1, u)$ or $\mathbb{P}(x_1, u)$ contains only 0-edges.

If $\mathbb{P}(x_1, u)$ contains a 1-edge, then all edges along the path $\mathbb{P}(u, x_n)$ must be 0, and all the edge on path $\mathbb{P}(u, x_2)$ must be 0, However, this implies that $x_2 \overset{0}{\sim} x_n$, a contradiction, as we assumed that $\overset{0}{\sim}$ is discrete.

Thus, there is no 1-edge in $\mathbb{P}(x_1, u)$ and hence, both paths $\mathbb{P}(u, x_n)$ and $\mathbb{P}(u, x_2)$ contain each exactly one 1-edge.

Now consider the unique vertex $v$ that lies on all three paths $\mathbb{P}(x_1, x_2)$, $\mathbb{P}(x_1, x_3)$, and $\mathbb{P}(x_2, x_3)$.

Since $u, v \in \mathbb{P}(x_1, x_2)$, we have either (A) $v \in \mathbb{P}(x_1, u)$ where $u = v$ is possible, or (B) $u \in \mathbb{P}(x_1, v)$ and $u \neq v$. We consider the two cases separately.

Case (A): Since there is no 1-edge in $\mathbb{P}(x_1, u)$ and $x_1 \overset{1}{\sim} x_n$, resp., $x_1 \overset{1}{\sim} x_2$ there is exactly one 1-edge in $\mathbb{P}(u, x_n)$, resp., $\mathbb{P}(u, x_2)$. Moreover, since $x_2 \overset{1}{\sim} x_3$ the path $\mathbb{P}(v, x_3)$ contains only 0-edges, and thus $x_3 \overset{1}{\sim} x_n$, a contradiction.

Case (B): Since there is no 1-edge in $\mathbb{P}(x_1, u)$ and $x_1 \overset{1}{\sim} x_n$, the path $\mathbb{P}(u, x_n)$ contains exactly one 1-edge.

In the following, we consider paths between two vertices $x_i, x_{n-i} \in \{x_1, \dots, x_n\}$ step-by-step, starting with $x_1$ and $x_{n-1}$.

The induction hypothesis implies that $x_1 \overset{1}{\nsim} x_{n-1}$ and since $\overset{0}{\sim}$ is discrete, we can conclude that $x_1 \overset{0}{\nsim} x_{n-1}$. Let $\mathbb{P}(x_1, x_n) = \mathbb{P}(x_1, a) \cup ab \cup \mathbb{P}(b, x_n)$ where $e = ab$ is the 1-edge contained in $\mathbb{P}(x_1, x_n)$. Let $c_1$ be the unique vertex that lies on all three paths $\mathbb{P}(x_1, x_n), \mathbb{P}(x_1, x_{n-1})$, and $\mathbb{P}(x_{n-1}, x_n)$. If $c_1$ lies on the path $\mathbb{P}(x_1, a)$, then $\mathbb{P}(c_1, x_{n-1})$ contains only 0-edges, since $\mathbb{P}(x_{n-1}, x_n) = \mathbb{P}(x_{n-1}, c_1) \cup \mathbb{P}(c_1, a) \cup ab \cup \mathbb{P}(b, x_n)$ and $x_n \overset{1}{\sim} x_{n-1}$. However, in this case the path $\mathbb{P}(x_1, c_1) \cup \mathbb{P}(c_1, x_{n-1})$ contains only 0-edges, which implies that $x_1 \overset{0}{\sim} x_{n-1}$, a contradiction. Thus, the vertex $c_1$ must be contained in $\mathbb{P}(b, x_n)$. Since $x_1 \overset{1}{\sim} x_n$, the path $\mathbb{P}(c_1, x_n)$ contains only 0-edges. Hence, the path $\mathbb{P}(c_1, x_{n-1})$ contains exactly one 1-edge, because $x_n \overset{1}{\sim} x_{n-1}$. In particular, by construction we see that $\mathbb{P}(x_1, x_n) = \mathbb{P}(x_1, u) \cup \mathbb{P}(u, c_1) \cup \mathbb{P}(c_1, x_n)$.

Now consider the vertices $x_n$ and $x_{n-2}$. Let $a'b'$ be the 1-edge on the path $\mathbb{P}(c_1, x_{n-1}) = \mathbb{P}(c_1, a') \cup a'b' \cup \mathbb{P}(b', x_{n-1})$. Since $x_n \overset{1}{\nsim} x_{n-2}$ and $x_n \overset{0}{\nsim} x_{n-2}$ we can apply the same argument and conclude that there is a vertex $c_2 \in \mathbb{P}(b', x_{n-1})$ s.t. the path $\mathbb{P}(c_2, x_{n-2})$ contains exactly one 1-edge. In particular, by construction we see that $\mathbb{P}(x_1, x_{n-2}) = \mathbb{P}(x_1, u) \cup \mathbb{P}(u, c_1) \cup \mathbb{P}(c_1, c_2) \cup \mathbb{P}(c_2, x_2)$ s.t. the path $\mathbb{P}(c_1, c_2)$ contains exactly one 1-edge.

Repeating this argument, we arrive at vertices $x_2$ and $x_4$ and can conclude analogously that there is a path $\mathbb{P}(c_{n-2}, x_2)$ that contains exactly one 1-edge and in particular, that $\mathbb{P}(x_1,x_2) = \mathbb{P}(x_1,u) \cup \mathbb{P}(u,c_1) \cup \left(\bigcup_{1\le i\le n-3} \mathbb{P}(c_i,c_{i+1})\right) \cup \mathbb{P}(c_{n-2},x_{n-2})$, where each of the distinct paths $\mathbb{P}(c_i,c_{i+1})$, $1 \le i \le n-3$ contains exactly one 1-edge. However, this contradicts that $x_1 \overset{1}{\sim} x_2$. □

**Corollary 5.4.1.** The graph $G(\overset{1}{\sim})/\overset{0}{\sim}$ is a forest, and hence all paths in $G(\overset{1}{\sim})/\overset{0}{\sim}$ are induced paths.

Next we analyze the effect of edge contractions in $T$.

**Lemma 5.4.4.** *Let $(T,\lambda)$ explain $G(\overset{1}{\sim})/\overset{0}{\sim}$ and let $(T',\lambda')$ be the result of contracting an interior edge $e$ in $T$. If $\lambda(e) = 0$ then $(T',\lambda')$ explains $G(\overset{1}{\sim})/\overset{0}{\sim}$. If $G(\overset{1}{\sim})/\overset{0}{\sim}$ is connected and $\lambda(e) = 1$ then $(T',\lambda')$ does not explain $G(\overset{1}{\sim})/\overset{0}{\sim}$.*

*If $G(\overset{1}{\sim})/\overset{0}{\sim}$ is connected and $(T,\lambda)$ is a tree that explains $G(\overset{1}{\sim})/\overset{0}{\sim}$, then $(T,\lambda)$ is least resolved if and only if all 0-edge are incident to leaves and each inner vertex is incident to exactly one 0-edge.*

*If, in addition, $(T,\lambda)$ is minimally resolved, then all 0-edge are incident to leaves and each inner vertex is incident to exactly one 0-edge.*

*Proof.* Let $\overset{1}{\sim}_T$, $\overset{0}{\sim}_T$, $\overset{1}{\sim}_{T'}$, and $\overset{0}{\sim}_{T'}$ be the relations explained by $T$ and $T'$, respectively. Since $(T,\lambda)$ explains $G(\overset{1}{\sim})/\overset{0}{\sim}$ , we have $\overset{1}{\sim}_T = \overset{1}{\sim} / \overset{0}{\sim}$. Moreover, since $\overset{0}{\sim}_T$ is discrete, no two distinct leaves of $T$ are in relation $\overset{0}{\sim}_T$.

If $\lambda(e) = 0$, then contracting the interior edge $e$ clearly preserves the property of $\overset{0}{\sim}_{T'}$ being discrete. Since only interior edges are allowed to be contracted, we have $L(T) = L(T')$. Therefore, $\overset{0}{\sim}_T = \overset{0}{\sim}_{T'}$ and the 1-edges along any path from $x \in L(T) = L(T')$ to $y \in L(T) = L(T')$ remains unchanged, and thus $\overset{1}{\sim}_T = \overset{1}{\sim}_{T'}$. Hence, $(T',\lambda')$ explains $G(\overset{1}{\sim})/\overset{0}{\sim}$.

If $G(\overset{1}{\sim})/\overset{0}{\sim}$ is connected, then for every 1-edge $e$ there is a pair of leaves $x'$ and $x''$ such that $x' \overset{1}{\sim} x''$ and $e$ is the only 1-edge along the unique path connecting $x'$ and $x''$. Consequently, contracting $e$ would make $x'$ and $x''$ non-adjacent w.r.t. the resulting relation. Thus no 1-edge can be contracted in $T$ without changing $G(\overset{1}{\sim})/\overset{0}{\sim}$.

Now assume that $(T,\lambda)$ is a least resolved tree that explains the connected graph $G(\overset{1}{\sim})/\overset{0}{\sim}$. By the latter arguments, all interior edges of $(T,\lambda)$ must be 1-edges and thus any 0-edge must be incident to leaves. Assume for contradiction, that there is an inner vertex $w$ such that for all adjacent leaves $x', x''$ we have $\lambda(wx') = \lambda(wx'') = 1$. Thus, for any such leaves we have $x' \overset{1}{\nsim} x''$. In particular, any path from $x'$ to any other leaf $y$ (distinct from the leaves adjacent to $w$) contains an interior 1-edge. Thus, $x' \overset{1}{\nsim} y$ for any such leaf of $T$. However, this immediately implies that $x'$ is an isolated vertex in $G(\overset{1}{\sim})/\overset{0}{\sim}$; a contradiction to the connectedness of $G(\overset{1}{\sim})/\overset{0}{\sim}$. Furthermore, if there is an

**Algorithm 1** Compute $(T(Q), \lambda)$

**Require:** $Q$

**Ensure:** $T(Q)$

set $T(Q) \leftarrow Q$

Retain the labels of all leaves of $Q$ in $T(Q)$ and relabel all inner vertices $u$ of $Q$ as $u'$.

Label all edges of the copy of $Q$ by $\lambda(e) = 1$.

For each inner vertex $u'$ of $Q$ add a vertex $u$ to $T(Q)$ and insert the edge $uu'$.

Label all edges of the form $e = uu'$ with $\lambda(e) = 0$.

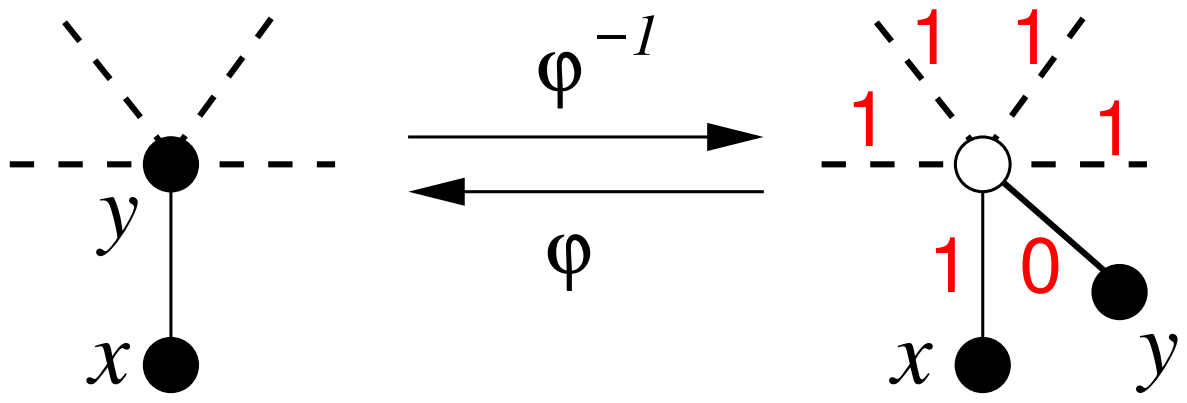

Figure 5.1: Illustration of the bijection $\varphi$. It contracts, at each inner vertex (white) of $(T, \lambda)$ the unique 0-edge and transfers $y \in X$ as vertex label at the inner vertex of $Q$. Leafs in $(T, \lambda)$ with incident to 1-edges remain unchanged. The inverse map $\varphi^{-1}$ is given by Algorithm 1.

inner vertex $w$ such that for adjacent leaves $x', x''$ it holds that $\lambda(wx') = \lambda(wx'') = 0$, then $x' \overset{0}{\sim} x''$; a contradiction to $\overset{0}{\sim}$ being discrete. Therefore, each inner vertex is incident to exactly one 0-edge.

If $G(\overset{1}{\sim})/ \overset{0}{\sim}$ is connected and $(T, \lambda)$ explains $G(\overset{1}{\sim})/ \overset{0}{\sim}$ such that all 0-edge are incident to leaves, then all interior edges are 1-edges. As shown, no interior 1-edge can be contracted in $T$ without changing the corresponding $\overset{1}{\sim}$ relation. Moreover, no leaf-edge can be contracted since $L(T) = V(G(\overset{1}{\sim})/ \overset{0}{\sim})$. Hence, $(T, \lambda)$ is least resolved.

Finally, since any minimally resolved tree is least resolved, the last assertions follows from the latter arguments. $\square$

Let $\mathbb{T}$ be the set of all trees with vertex set $X$ but no edge-labels and $\mathcal{T}$ denote the set of all edge-labeled 0-1-trees $(T, \lambda)$ with leaf set $X$ such that each inner vertex $w \in W$ has degree at least 3 and there is exactly one adjacent leaf $v$ to $w$ with $\lambda(wv) = 0$ while all other edges $e$ in $T$ have label $\lambda(e) = 1$.

**Lemma 5.4.5.** *The map $\varphi : \mathcal{T} \to \mathbb{T}$ with $\varphi : (T, \lambda) \mapsto Q$, $V(Q) = L(T)$, and $Q \simeq T^*$, where $T^*$ is the underlying unlabeled tree obtained from $(T, \lambda)$ by contracting all edges labeled 0, is a bijection.*

*Proof.* We show first that $\varphi$ and $\varphi^{-1}$ are maps. Clearly, $\varphi$ is a map, since the edge-contraction is well-defined and leads to exactly one tree in $\mathbb{T}$. For $\varphi^{-1}$ we construct $(T, \lambda)$ from $Q$ as in Algorithm 1. It is easy to see that $(T, \lambda) \in \mathcal{T}$. Now consider $T^*$ obtained from $(T, \lambda)$ by contracting all edges labeled 0. By construction, $T^* \simeq Q$ (see Fig. 5.1). Hence, $\varphi : \mathcal{T} \to \mathbb{T}$ is bijective. □

The bijection is illustrated in Fig. 5.1.

**Lemma 5.4.6.** *Let $(T, \lambda) \in \mathcal{T}$ and $Q = \varphi((T, \lambda)) \in \mathbb{T}$. The set $\mathcal{T}$ contains* all *least resolved trees that explain $Q$.*

*Moreover, if $Q$ is considered as a graph $G(\overset{1}{\sim})/\overset{0}{\sim}$ with vertex set $X$, then $(T, \lambda)$ is the unique least resolved tree $\mathcal{T}$ that explains $Q$ and therefore, the unique minimally resolved tree that explains $Q$.*

*Proof.* We start with showing that $(T, \lambda)$ explains $Q$. Note, since $Q \in \mathbb{T}$, the graph $Q$ must be connected. By construction and since $Q = \varphi((T, \lambda))$, $T^* \simeq Q$ where $T^*$ is the tree obtained from $T$ after contracting all 0-edges. Let $v, w \in X$ and assume that there is exactly a single 1 along the path from $v$ to $w$ in $(T, \lambda)$. Hence, after contracting all edges labeled 0 we see that $vw \in E(T^*)$ where $T^* \simeq Q$ and thus $v \overset{1}{\sim} w$. Note, no path between any two vertices in $(T, \lambda)$ can have only 0-edges (by construction). Thus, assume that there is more than a single 1-edge on the path between $v$ and $w$. Hence, after after contracting all edges labeled 0 we see that there is still a path in the tree $T^*$ from $v$ to $w$ with more than one 1-edge. Since $T^* \simeq Q$, we have $vw \notin E(Q)$ and therefore, $v \overset{1}{\not\sim} w$. Thus, $(T, \lambda)$ explains $Q$.

By construction of the trees in $\mathcal{T}$ all 0-edges are incident to a leaf. Thus, by Lemma 5.4.4, every least resolved tree that explains $Q$ is contained in $\mathcal{T}$.

It remains to show that the least resolved tree $(T, \lambda) \in \mathcal{T}$ with $T^* \simeq Q$ that explains $Q$ is minimally resolved. Assume there is another least resolved tree $(T', \lambda') \in \mathcal{T}$ with leaf set $V(Q)$ that explains $Q$. By Lemma 5.4.5, there is a bijection between those $Q$ and elements in $\mathcal{T}$ for which $T^* \simeq Q$. Thus, $T'^* \simeq Q' \not\simeq Q \simeq T^*$. However, this implies that $T' \not\simeq T$. However, since in this case $(T', \lambda')$ explains $Q'$ and $Q' \not\simeq Q$, the pair $(T', \lambda')$ cannot explain $Q$; a contradiction. □

As an immediate consequence of these considerations we obtain

**Theorem 5.4.2.** *Let $Q$ be a connected component in $G(\overset{1}{\sim})/\overset{0}{\sim}$ with vertex set $X$. Then the tree $(T, \lambda)$ constructed in Algorithm 1 is the unique minimally resolved tree that explains $Q$.*

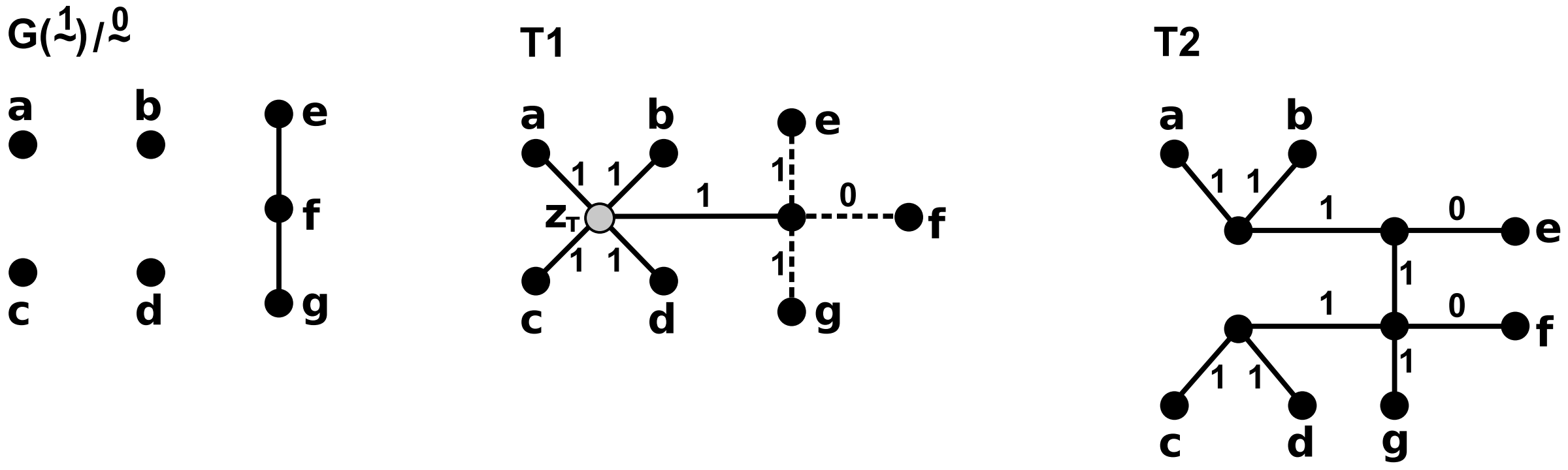

Figure 5.2: Least resolved trees explaining a given relation are unique whenever $G(\overset{1}{\sim})/\overset{0}{\sim}$ is connected (cf. Lemma 5.4.6) and are, therefore, also minimally resolved. Now, consider the disconnected graph $G(\overset{1}{\sim})/\overset{0}{\sim}$ shown on the left. Both pairs $(T_1, \lambda_1)$ (middle) and $(T_2, \lambda_2)$ (right) are least resolved trees that explain $G(\overset{1}{\sim})/\overset{0}{\sim}$. Thus, uniqueness of least resolved trees that explain $G(\overset{1}{\sim})/\overset{0}{\sim}$ is not always satisfied. The tree $(T_1, \lambda_1)$ is obtained with Alg. 2, by connecting the vertex $z_T$ via 1-edges to the (inner) vertices of the respective minimally resolved trees $T(Q)$ (the dashed subtree and the single vertex trees $a, b, c, d$) that explain the connected components $Q$ of $G(\overset{1}{\sim})/\overset{0}{\sim}$. In this example, $(T_1, \lambda_1)$ is the unique minimally resolved tree that explains $G(\overset{1}{\sim})/\overset{0}{\sim}$, since each tree $T(Q)$ has a unique interior vertex and due to Thm. 5.4.3.

*Moreover, for any pair $(T', \lambda')$ that explains $Q$, the tree $T$ is obtained from $T'$ by contracting all* interior *0-edges and putting $\lambda(e) = \lambda'(e)$ for all edges that are not contracted.*

*Proof.* The first statement follows from Lemma 5.4.4 and 5.4.6. To see the second statement, observe that Lemma 5.4.4 implies that no interior 1-edge but every 0-edge can be contracted. Hence, after contracting all 0-edges, no edge can be contracted and thus, the resulting tree is least resolved. By Lemma 5.4.4, we obtain the result. □

We emphasize that although the minimally resolved tree that explains $G(\overset{1}{\sim})/\overset{0}{\sim}$ is unique, this statement is in general not satisfied for least resolved trees, see Figure 5.2.

We are now in the position to demonstrate how to obtain a least resolved tree that explains $G(\overset{1}{\sim})/\overset{0}{\sim}$ also in the case that $G(\overset{1}{\sim})/\overset{0}{\sim}$ itself is not connected. To this end, denote by $Q_1, \dots Q_k$ the connected components of $G(\overset{1}{\sim})/\overset{0}{\sim}$. We can construct a phylogenetic tree $T(G(\overset{1}{\sim})/\overset{0}{\sim})$ with leaf set $X$ for $G(\overset{1}{\sim})/\overset{0}{\sim}$ using Alg. 2. It basically amounts to constructing a star $S_k$ with inner vertex $z$, where its leaves are identified with the trees $T(Q_i)$.

**Algorithm 2** Compute $(T(G(\overset{1}{\sim})/\overset{0}{\sim})), \lambda)$

**Require:** disconnected $G(\overset{1}{\sim})/\overset{0}{\sim})$

**Ensure:** $T(G(\overset{1}{\sim})/\overset{0}{\sim}))$

```
T(G(~1)/~0)) ← ({z_T}, ∅)
for For each connected component Q_i do
    construct (T(Q_i), λ_i) with Alg. 1 and add to T(G(~1)/~0)).
    if T(Q_i) is the single vertex graph ({x}, ∅) then
        add edge z_T x
    else if T(Q_i) is the edge v_i w_i then
        remove the edge v_i w_i from T(Q_i), insert a vertex x_i in T(Q_i) and the edges x_i v_i,
        x_i w_i.
        set either λ_i(x_i v_i) = 1 and λ(x_i w_i) = 0 or λ_i(x_i v_i) = 0 and λ(x_i w_i) = 1.
        add edge z_T x_i to T(G(~1)/~0)).
    else
        add edge z_T q'_i to T(G(~1)/~0)) for an arbitrary inner vertex q'_i of T(Q_i).
    end if
end for
Set λ(z_T v) = 1 for all edges z_T v and λ(e) = λ_i(e) for all edges e ∈ T(Q_i).
```

**Lemma 5.4.7.** *Let $G(\overset{1}{\sim})/\overset{0}{\sim}$ have connected components $Q_1, \dots Q_k$. Let $T'$ be a tree that explains $G(\overset{1}{\sim})/\overset{0}{\sim}$ and $T'_i$ be the subtree of $T'$ with leaf set $V(Q_i)$ that is minimal w.r.t. inclusion, $1 \leq i \leq k$. Then, $V(T'_i) \cap V(T'_j) = \emptyset$, $i \neq j$ and, in particular, any two vertices in $T'_i$ and $T'_j$, respectively, have distance at least two in $T'$.*

*Proof.* We start to show that two distinct subtrees $T'_i$, and $T'_j$ do not have a common vertex in $T'$. If one of $Q_i$ or $Q_j$ is a single vertex graph, then $T'_i$ or $T'_j$ consists of a single leaf only, and the statement holds trivially.

Hence, assume that both $Q_i$ and $Q_j$ have at least three vertices. Lemma 5.4.4 implies that each inner vertex of the minimally resolved trees $T(Q_i)$ and $T(Q_j)$ is incident to exactly one 0-edge as long as $Q_i$ is not an edge. Since $T(Q_i)$ can be obtained from $T'_i$ by the procedure above, for each inner vertex $v$ in $T'_i$ there is a leaf $x$ in $T'_i$ such that the unique path from $v$ to $x$ contains only 0-edges. The same arguments apply, if $Q_i$ is an edge $xy$. In this case, the tree $T'_i$ must have $x$ and $y$ as leaves, which implies that $T'_i$ has at least one inner vertex $v$ and that there is exactly one 1-edge along the path from $x$ to $y$. Thus, for each inner vertex in $T'_i$ there is a path to either $x$ or $y$ that contains only 0-edges.

Let $v$ and $w$ be arbitrary inner vertices of $T'_i$ and $T'_j$, respectively, and let $x$ and $y$ be leaves that are connected to $v$ and $w$, resp., by a path that contains only 0-edges. If $v = w$, then $x \overset{0}{\sim} y$, contradicting the property of $\overset{0}{\sim}$ being discrete. Thus, $T'_i$ and $T'_j$ cannot have a common vertex in $T'$. Moreover, there is no edge $vw$ in $T'$, since otherwise either $x \overset{0}{\sim} y$ (if $\lambda(vw) = 0$) or $x \overset{1}{\sim} y$ (if $\lambda(vw) = 1$). Hence, any two distinct vertices in $T'_i$ and $T'_j$ have distance at least two in $T'$. □

We note that Algorithm 2 produces a tree with a single vertex of degree 2, namely $z_T$, whenever $G(\overset{1}{\sim})/\overset{0}{\sim}$ consists of exactly two components. Although this strictly speaking violates the definition of phylogenetic trees, we tolerate this anomaly for the remainder of this section.

**Theorem 5.4.3.** *Let $Q_1, \dots Q_k$ be the connected components in $G(\overset{1}{\sim})/\overset{0}{\sim}$. Up to the choice of the vertices $q'_i$ in Line 11 of Alg. 2, the tree $T^* = T(G(\overset{1}{\sim})/\overset{0}{\sim}))$ is a minimally resolved tree that explains $G(\overset{1}{\sim})/\overset{0}{\sim}$. It is unique up to the choice of the $z_T q'_i$ in Line 11.*

*Proof.* Since every tree $T(Q_i)$ explains a connected component in $G(\overset{1}{\sim})/\overset{0}{\sim}$, from the construction of $T^*$ it is easily seen that $T^*$ explains $G(\overset{1}{\sim})/\overset{0}{\sim}$. Now we need to prove that $T^*$ is a minimally resolved tree that explains $G(\overset{1}{\sim})/\overset{0}{\sim}$.

To this end, consider an arbitrary tree $T'$ that explains $G(\overset{1}{\sim})/\overset{0}{\sim}$. Since $T'$ explains $G(\overset{1}{\sim})/\overset{0}{\sim}$, it must explain each of the connected components $Q_1, \dots Q_k$. Thus, each of the subtrees $T'_i$ of $T'$ with leaf set $V(Q_i)$ that are minimal w.r.t. inclusion must explain the connected component $Q_i$, $1 \leq i \leq k$. Note the $T'_i$ may have vertices of degree 2.

We show first that $T(Q_i)$ is obtained from $T'_i$ by contracting all interior 0-edges and all 0-edges of degree 2. If there are no vertices of degree 2, we can immediately apply Thm. 5.4.2.

If there is a vertex $v$ of degree 2, then $v$ cannot be incident to two 1-edges, as otherwise the relation explained by $T'_i$ would not be connected, contradicting the assumption that $T'_i$ explains the connected component $Q_i$. Thus, if there is a vertex $v$ of degree 2 it must be incident to a 0-edge $vw$. Contracting $vw$ preserves the property of $\stackrel{0}{\sim}$ being discrete. If $w$ is a leaf, we can contract the edge $vw$ to a new leaf vertex $w$; if $vw$ is an interior edge we simply contract it to some new inner vertex. In both cases, we can argue analogously as in the proof of Lemma 5.4.4 that the tree obtained from $T'_i$ after contracting $vw$, still explains $Q_i$. This procedure can be repeated until no degree-two vertices are in the contracted $T'_i$.

In particular, the resulting tree is a phylogenetic tree that explains $Q_i$. Now we continue to contract all remaining interior 0-edges. Thm. 5.4.2 implies that in this manner we eventually obtain $T(Q_i)$.

By Lemma 5.4.7, two distinct tree $T(Q_i)$ and $T(Q_j)$ do not have a common vertex, and moreover, any two vertices in $T(Q_i)$ and $T(Q_j)$, respectively, have distance at least two in $T'$.

This implies that the construction as in Alg. 2 yields a least resolved tree. In more detail, since the subtrees explaining $Q_i$ in any tree that explains $G(\stackrel{1}{\sim})/\stackrel{0}{\sim}$ must be vertex disjoint, the minimally resolved trees $T(Q_1), \ldots, T(Q_k)$ must be subtrees of any minimally resolved tree that explain $G(\stackrel{1}{\sim})/\stackrel{0}{\sim}$, as long as all $Q_i$ are single vertex graphs or have at least one inner vertex.

If $Q_i$ is a single edge $v_iw_i$ and thus $T(Q_i) = v_iw_i$ where $\lambda(v_iw_i) = 1$, we modify $T(Q_i)$ in Line 8 to obtain a tree isomorphic to $S_2$ with inner vertex $x_i$. This modification is necessary, since otherwise (at least one of) $v_i$ or $w_i$ would be an inner vertex in $T^*$, and we would loose the information about the leaves $v_i, w_i$. In particular, we need to add this vertex $x_i$ because we cannot attach the leaves $v_i$ (resp. $w_i$) by an edge $x_jv_i$ (resp. $x'_jw_i$) to some subtree subtree $T(Q_j)$. To see this, note that at least one of the edges $x_jv_i$ and $x'_jw_i$ must be a 0-edge. However, $x_j$ and $x'_j$ are already incident to a 0-edge $x_jv'_i$ or $x'_jw'_i$ (cf. Lemma 5.4.4), which implies that $\stackrel{0}{\sim}$ would not be discrete; a contradiction. By construction, we still have $v_i \stackrel{1}{\sim} w_i$ in Line 8.

Finally, any two distinct vertices in $T(Q_i)$ and $T(Q_j)$ have distance at least two in $T^*$, as shown above. Hence, any path connecting two subtrees $T(Q_i)$ in $T^*$ contains and least two edges and hence at least one vertex that is not contained in any of the $T(Q_i)$. Therefore, any tree explaining $Q$ has at least $1 + \sum_i |V(T(Q_i))|$ vertices.

We now show that adding a single vertex $z_T$, which we may consider as a trivial tree $(\{z_T\}, \emptyset)$, is sufficient. Indeed, we may connect the different trees to $z_T$ by insertion of an edge $z_T q_i'$, where $q_i'$ is an arbitrary inner vertex of $T(Q_i)$ and label these edges $\lambda(z_T q_i') = 1$. Thus, no two leaves $u$ and $w$ of distinct trees are either in relation $\overset{0}{\sim}$ or $\overset{1}{\sim}$, as required. The resulting trees have the minimal possible number of vertices, i.e., they are minimally resolved. □

## Binary trees

Instead of asking for least resolved trees that explain $G(\overset{1}{\sim})/\overset{0}{\sim}$, we may also consider the other extreme and ask which binary, i.e., fully resolved tree can explain $G(\overset{1}{\sim})/\overset{0}{\sim}$. Recall that an $X$-tree is called binary or fully resolved if the root has degree 2 while all other inner vertices have degree 3. From the construction of the least resolved trees we immediately obtain the following:

**Corollary 5.4.4.** A least resolved tree $T(Q)$ for a connected component $Q$ of $G(\overset{1}{\sim})/\overset{0}{\sim}$ is binary if and only if $Q$ is a path.

If a least resolved tree $T(Q)$ of $G(\overset{1}{\sim})/\overset{0}{\sim}$ is a star, we have:

**Lemma 5.4.8.** *If a least resolved tree $T(Q)$ explaining $G(\overset{1}{\sim})/\overset{0}{\sim}$ is a star with $n$ leaves, then either*

*(a) all edges in $T(Q)$ are 1-edges and $Q$ has no edge, or*

*(b) there is exactly one 0-edge in $T(Q)$ and $Q$ is a star with $n-1$ leaves.*

*Proof.* For implication in case (a) and (b) we can re-use exactly the same arguments as in the proofs of Theorem 5.4.2 and 5.4.3.

Now suppose there are at least two (incident) 0-edges in $T(Q)$, whose endpoints are the vertices $u$ and $v$. Then $u \overset{0}{\sim} v$, which is impossible in $G(\overset{1}{\sim})/\overset{0}{\sim}$. □

**Lemma 5.4.9.** *Let $(T, \lambda)$ be a least resolved tree that explains $G(\overset{1}{\sim})/\overset{0}{\sim}$. Consider an arbitrary subgraph $S_k$ that is induced by an inner vertex $v_0$ and* all *of its $k$ neighbors. Then, $S_k$ with its particular labeling $\lambda_{|E(S_k)}$ is always of type (a) or (b) as in Lemma 5.4.8.*

*Proof.* This is an immediate consequence of Lemma 5.4.4 and the fact that $\overset{0}{\sim}$ is discrete. □

To construct the binary tree explaining the star $Q = S_n$, we consider the set of all binary trees with $n$ leaves and 0/1-edge labels. If $S_n$ is of type (a) in Lemma 5.4.8, then all terminal edges are labeled 1 and all interior edges are arbitrarily labeled 0 or 1. Figure 5.3 shows an example for $S_6$. If $S_n$ is of type (b), we label the terminal edges in the same way as in $T(Q)$ and all interior edge are labeled 0. In this case, for each binary tree there is exactly one labeling.

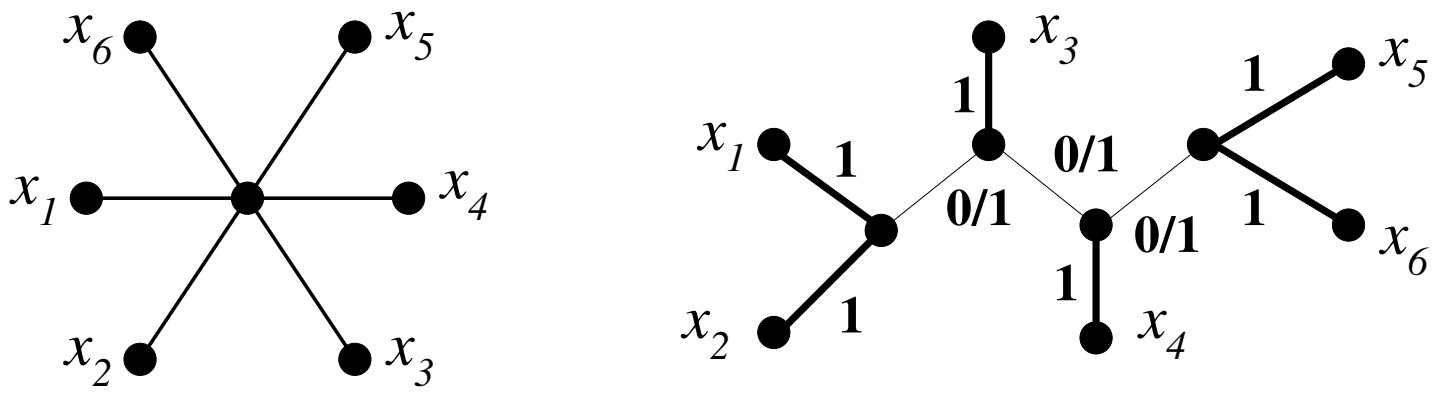

Figure 5.3: For fixed underlying tree $T(Q)$, in this case a star $S_6$ with all 1, there are in general multiple labelings $\lambda$ that explains the same relation $Q$, here the empty relation.

In order to obtain the complete set of binary trees that explain $G$ we can proceed as follows. If $G$ is connected, there is a single minimally resolved tree $T(G)$ explaining $G$. If $G$ is not connected then there are multiple minimally resolved trees $T$.

Let $\mathcal{T}_{\text{lrt}}$ be the set of all least resolved trees that explain $G$. For every such least resolved tree $T \in \mathcal{T}_{\text{lrt}}$ we iterate over all vertices $v_0$ of $T$ with degree $k > 3$ and perform the following manipulations:

1. Given a vertex $v_0$ of $T$ with degree $k > 3$, denote the set of its neighbors $v_1, v_2, \ldots, v_k$ by $N(v_0)$. Delete vertex $v_0$ and its attached edges from $T$, and rename the neighbors $v_i$ to $v_i'$ for all $1 \leq i \leq k$. Denote the resulting forest by $F(v_0)$.

2. Generate all binary trees with leaves $v_1', \ldots, v_k'$.

3. Each of these binary trees is inserted into a copy of the forest $F(v_0)$ by identifying $v_i$ and $v_i'$ for all $1 \leq i \leq k$.

4. For each of the inserted binary trees $T'$ that results from a "local" star $S_k$ in step 3. we must place an edge label.
   Put $\lambda(xv_i) = \lambda(v_0 v_i)$ for all edges $xv_i$ in $T'$ and mark $xv_i$ as *LABELED.* If $S_k$ is of type (a) (cf. Lemma 5.4.8), then choose an arbitrary 0/1-label for the interior edges of $T'$. If $S_k$ is of type (b) we need to consider the two exclusive cases for the vertex $v_j$ for which $\lambda(v_j x) = 0$:

(i) For all $y \in V(G)$ for which $v_j \overset{1}{\sim} y$, label all *interior* edges on the unique path $\mathbb{P}(v_j, y)$ that are also contained in $T'$ and are not marked as *LABELED* with 0 and mark them as *LABELED* and choose an arbitrary 0/1-label for all other un*LABELED* interior edges of $T'$.

(ii) Otherwise, choose an arbitrary 0/1-label for the interior edges of $T'$.

It is well known that each binary tree has $k-3$ interior edges [72]. Hence, for a binary tree there are $2^{k-3}$ possibilities to place a 0/1 label on its interior edges. Let $t(k)$ denote the number of binary trees with $k$ leaves and $V_a, V_b$ be a partition of the inner vertices into those where the neighborhood corresponds to a star of type (a) and (b), respectively. Note, if $T$ is minimally resolved, then $|V_a| \leq 1$. For a given least resolved tree $T \in \mathcal{T}_{\text{lrt}}$, the latter procedure yields the set of all $(\prod_{v \in V_a} t(\deg(v_0)) 2^{\deg(v_0)-3}) \cdot (\prod_{v \in V_b} t(\deg(v_0)))$ pairwise distinct binary trees that one can obtain from $T$. The union of these tree sets and its particular labeling over all $T \in \mathcal{T}_{\text{lrt}}$ is then the set of all binary trees explaining $G$. To establish the correctness of this procedure, we prove

**Lemma 5.4.10.** *The procedure outlined above generates all binary trees $(T, \lambda)$ explaining $G$.*

*Proof.* We first note that there may not be a binary tree explaining $G$. This is case whenever $T(G)$ has a vertex of degree 2, which is present in particular if $G$ is forest with two connected components.

Now consider an arbitrary binary tree $(T_B(G), \lambda)$ that is not least resolved for $G$. Then a least resolved tree $T(G)$ explaining $G$ can be obtained from $T_B(G)$ by contracting edges and retaining the the labeling of all non-contracted edges. In the following we will show that the construction above can be be used to recover $T_B(G) = (V, E)$ from $T(G)$. To this end, first observe that only interior edges can be contracted in $T_B(G)$ to obtain $T(G)$. Let $E' = \{e_1, \ldots, e_h\}$ be a maximal (w.r.t. inclusion) subset of contracted edges of $T_B(G)$ such that the subgraph $(V' = \cup_{i=1}^{h} e_i, E')$ is connected, and thus forms a subtree of $T_B(G)$. Furthermore, let $F = \{f_1, \ldots, f_k\} \subseteq E \setminus E'$ be a maximal subset of edges of $T_B(G)$ such that for all $f_i \in F$ there is an edge $e_j \in E'$ such that $f_i \cap e_j \neq \emptyset$. Moreover, set $W = \cup_{i=1}^{k} f_i$. Thus, the contracted subtree $(V', E')$ locally corresponds to the vertex $v_0$ of degree $k > 3$ and thus, to a local star $S_k$. Now, replacing $S_k$ by the tree $(V' \cup W, E' \cup F)$ (as in Step 3) yields the subtree from which we have contracted all interior edges that are contained in $E'$. Since the latter procedure can be repeated for all such maximal sets $E'$, we can recover $T_B(G)$.

It remains to show that one can also recover the labeling $\lambda$ of $T_B(G)$. Since $T(G)$ is a least resolved tree obtained from $T_B(G)$ that explains $G$, we have, by definition, $u \overset{1}{\sim} w$ in $T(G)$ if and only if $u \overset{1}{\sim} w$ in $T_B(G)$. By Lemma 5.4.9, every local star $S_k$ in $T(G)$ is either of type (a) or (b). Assume it is of type (a), i.e., $\lambda(v_0v_i) = 1$ for all $1 \leq i \leq k$. Let $u$ and $w$ be leaves of $G$ for which the unique path connecting them contains the edge $v_0v_i$ and $v_0v_j$. Thus, there are at least two edges labeled 1 along the path; hence $u \overset{1}{\nsim} w$. The edges $xv_i$ and $x'v_j$ in $T'$ are both labeled by 1; therefore $u$ and $w$ are not in relation $\overset{1}{\sim}$ after replacing $S_k$ by $T'$. Therefore all possible labelings can be used in Step 4 except for the edges that are not contained in $T'$ and $xv_i$ and $x'v_j$ which are marked as *LABELED*). Therefore, we also obtain the given labeling of the subtree $T = (V' \cup W, E' \cup F)$ as a result.

If the local star $S_k$ is of type (b), then there is exactly one edge $v_0v_j$ with $\lambda(v_0v_j) = 0$. By Lemma 5.4.4 and because $(T(G), \lambda)$ is least resolved, the vertex $v_j$ must be a leaf of $G$. For two leaves $u, w$ of $T(G)$ there are two cases: either the unique path $\mathbb{P}(u, w)$ in $T(G)$ contains $v_0$ or not.

If $\mathbb{P}(u, w)$ does not contain $v_0$, then this path and its labeling remains unchanged after replacing $S_k$ by $T'$. Hence, the relations $\overset{1}{\sim}$ and $\overset{1}{\nsim}$ are preserved for all such vertices $u, w$.

First, assume that $\mathbb{P}(u, w)$ contains $v_0$ and thus two edges incident to vertices in $N(v_0)$. If the path $\mathbb{P}(u, w)$ contains two edges $v_0v_i$ and $v_0v_l$ with $i, l \neq j$, then $\lambda(v_0v_i) = \lambda(v_0v_l) = 1$. Thus, $u \overset{1}{\nsim} w$ in $T(G)$. The extended path (after replacing $S_k$ by $T'$) still contains the two 1-edges $xv_i$ and $x'v_l$, independently from the labeling of all other edges in $T'$ that have remained un*LABELED* edges up to this point. Thus, $u \overset{1}{\nsim} w$ is preserved after replacing $S_k$ by $T'$.

Next, assume that $\mathbb{P}(u, w)$ contains $v_0$ and the 0-edge $v_0v_j$ in $T(G)$. In the latter case, $u = v_j$ of $w = v_j$. Note, there must be another edge $v_0v_l$ in $\mathbb{P}(u, w)$ with $v_l \in N(v_0)$, $l \neq j$ and therefore, with $\lambda(v_0v_l) = 1$. There are two cases, either $u \overset{1}{\sim} w$ or $u \overset{1}{\nsim} w$ in $T(G)$.

If $u \overset{1}{\sim} w$ then there is exactly one 1-edge (the edge $v_0v_l$) contained in $\mathbb{P}(u, w)$ in $T(G)$. By construction, all interior edges on the path $\mathbb{P}(u, w)$ that are contained in $T'$ are labeled with 0 and all other edge-labelings remain unchanged in $T(G)$ after replacing $S_k$ by $T'$. Thus, $u \overset{1}{\sim} w$ in $T(G)$ after replacing $S_k$ by $T'$. Analogously, if $u \overset{1}{\nsim} w$, then there are at least two edges 1-edges $\mathbb{P}(u, w)$ in $T(G)$. Since $\lambda(v_0v_j) = 0$ and $\lambda(v_0v_l) = 1$, the 1-edge different from $v_0v_l$ is not contained in $T'$ and its label 1 remains unchanged. Moreover, the edge $xv_l$ in $T'$ gets also the label 1 in Step 4. Thus, $\mathbb{P}(u, w)$ still contains at least two 1-edges in $T(G)$ after replacing $S_k$ by $T'$ independently of the labeling chosen for the other un*LABELED* interior edges of $T'$. Thus, $u \overset{1}{\not\sim} w$ in $T(G)$ after replacing $S_k$ by $T'$.

We allow all possible labelings and fix parts where necessary. In particular, we obtain the labeling of the subtree $(V' \cup W, E' \cup F)$ that coincides with the labeling of $T_B(G)$. Thus, we can repeat this procedure for all stars $S_k$ in $T(G)$ and their initial labelings. Therefore, we can recover both $T_B(G)$ and its edge-labeling $\lambda$. Clearly, every binary tree $T_B(G)$ that explains $G$ is either already least resolved or there is a least tree $T(G) \in \mathcal{T}_{\text{lrt}}$ from which $T_B(G)$ can be recovered by the construction as outlined above. □

As a consequence of the proof of Lemma 5.4.10 we immediately obtain the following Corollary that characterizes the condition that $Q$ cannot be explained by a binary tree.

**Corollary 5.4.5.** $G(\overset{1}{\sim})/\overset{0}{\sim}$ cannot be explained by a binary tree if and only if $G(\overset{1}{\sim})/\overset{0}{\sim}$ is a forest with exactly two connected components.

The fact that exactly two connected components appear as a special case is the consequence of a conceptually too strict definition of "binary tree". If we allow a single "root vertex" of degree 2 in this special case, we no longer have to exclude two-component graphs.

## 5.5 The antisymmetric single-1 relation

The antisymmetric version $x \overset{1}{\rightharpoonup} y$ of the 1-relation shares many basic properties with its symmetric cousin. We therefore will not show all formal developments in full detail. Instead, we will where possible appeal to the parallels between $x \overset{1}{\rightharpoonup} y$ and $x \overset{1}{\sim} y$. For convenience we recall the definition: $x \overset{1}{\rightharpoonup} y$ if and only if all edges along $\mathbb{P}(u, x)$ are labeled 0 and exactly one edge along $\mathbb{P}(u, y)$ is labeled 1, where $u = lca(x, y)$. As an immediate consequence we may associate with $\overset{1}{\rightharpoonup}$ a symmetrized 1-relation $x \overset{1}{\sim} y$ whenever $x \overset{1}{\rightharpoonup} y$ or $y \overset{1}{\rightharpoonup} x$. Thus we can infer (part of) the underlying unrooted tree topology by considering the symmetrized version $\overset{1}{\sim}$. On the other hand, $\overset{1}{\rightharpoonup}$ cannot convey more information on the unrooted tree from which $\overset{1}{\rightharpoonup}$ and its symmetrization $\overset{1}{\sim}$ are derived. It remains, however, to infer the position of the root from directional information. Instead of the quadruples used for the unrooted trees in the previous section, structural constraints on rooted trees are naturally expressed in terms of triples.

In the previous section we have considered $\overset{1}{\sim}$ in relation to unrooted trees only. Before we start to explore $\overset{1}{\rightharpoonup}$ we first ask whether $\overset{1}{\sim}$ contains any information about the position of the root and if it already places any constraints on $\overset{1}{\rightharpoonup}$ beyond those derived for $\overset{1}{\sim}$ in the previous section. In general the answer to this question will be negative, as suggested by the example of the tree $T_5^*$ in Figure 5.4. Any of its inner vertex can be chosen as the root, and each choice of a root vertex yields a different relation $\overset{1}{\rightharpoonup}$.

Nevertheless, at least partial information on $\overset{1}{\rightharpoonup}$ can be inferred uniquely from $\overset{1}{\sim}$ and $\overset{0}{\sim}$. Since all connected components in $G(\overset{1}{\sim})/\overset{0}{\sim}$ are trees, we observe that the underlying graphs $\underline{G(\overset{1}{\rightharpoonup})/\overset{0}{\sim}}$ of all connected components in $G(\overset{1}{\rightharpoonup})/\overset{0}{\sim}$ must be trees as well. Moreover, since $\overset{0}{\sim}$ is discrete in $G(\overset{1}{\sim})/\overset{0}{\sim}$, it is also discrete in $G(\overset{1}{\rightharpoonup})/\overset{0}{\sim}$.

Let $Q$ be a connected component in $G(\overset{1}{\rightharpoonup})/\overset{0}{\sim}$. If $Q$ is an isolated vertex or a single edge, there is only a single phylogenetic rooted tree (a single vertex and a tree with two leaves and one inner root vertex, resp.) that explains $Q$ and the position of its root is uniquely determined.

Thus we assume that $Q$ contains at least three vertices from here on. By construction, any three vertices $x, y, z$ in a connected component $Q$ in $G(\overset{1}{\rightharpoonup})/\overset{0}{\sim}$ either induce a disconnected graph, or a tree on three vertices. Let $x, y, z \in V(Q)$ induce such a tree. Then there are three possibilities (up to relabeling of the vertices) for the induced subgraph contained in $G(\overset{1}{\rightharpoonup})/\overset{0}{\sim} = (V, E)$:

(i) $xy, yz \in E$ implying that $x \overset{1}{\rightharpoonup} y \overset{1}{\rightharpoonup} z$,

(ii) $yx, yz \in E$ implying that $y \overset{1}{\rightharpoonup} x$ and $y \overset{1}{\rightharpoonup} z$,

(iii) $xy, zy \in E$ implying that $x \overset{1}{\rightharpoonup} y$ and $z \overset{1}{\rightharpoonup} y$.

Below, we will show that Cases (i) and (ii) both imply a unique tree on the three leaves $x, y, z$ together with a unique 0/1-edge labeling for the unique resolved tree $T(Q)$ that displays $Q$, see Fig. 5.5. Moreover, we shall see that Case (iii) cannot occur.

**Lemma 5.5.1.** *In Case (i), the unique triple* $\mathsf{yz|x}$ *must be displayed by any tree* $T(Q)$ *that explains* $Q$. *Moreover, the paths* $\mathbb{P}(u, v)$ *and* $\mathbb{P}(v, z)$ *in* $T(Q)$ *contain both exactly one 1-edge, while the other paths* $\mathbb{P}(u, x)$ *and* $\mathbb{P}(v, y)$ *contain only 0-edges, where* $u = lca(xy) = lca(xz) \neq lca(yz) = v$.

*Proof.* Let $x, y, z \in V(Q)$ such that $xy, yz \in E$ and thus, $x \overset{1}{\rightharpoonup} y \overset{1}{\rightharpoonup} z$. Notice first that there must be two distinct last common ancestors for pairs of the three vertices $x, y, z$; otherwise, if $u = lca(xy) = lca(xz) = lca(yz)$, then the path $\mathbb{P}(uy)$ contains a 1-edge (since $x \overset{1}{\rightharpoonup} y$) and hence $y \overset{1}{\rightharpoonup} z$ is impossible. We continue to show that $u = lca(xy) = lca(xz)$. Assume that $u = lca(xy) \neq lca(xz)$. Hence, either the triple $\mathsf{xz|y}$ or $\mathsf{xy|z}$ is displayed by $T(Q)$. In either case the path $\mathbb{P}(u, y)$ contains a 1-edge, since $x \overset{1}{\rightharpoonup} y$. This, however, implies $y \overset{1}{\not\rightharpoonup} z$, a contradiction. Thus, $u = lca(xy) = lca(xz)$.

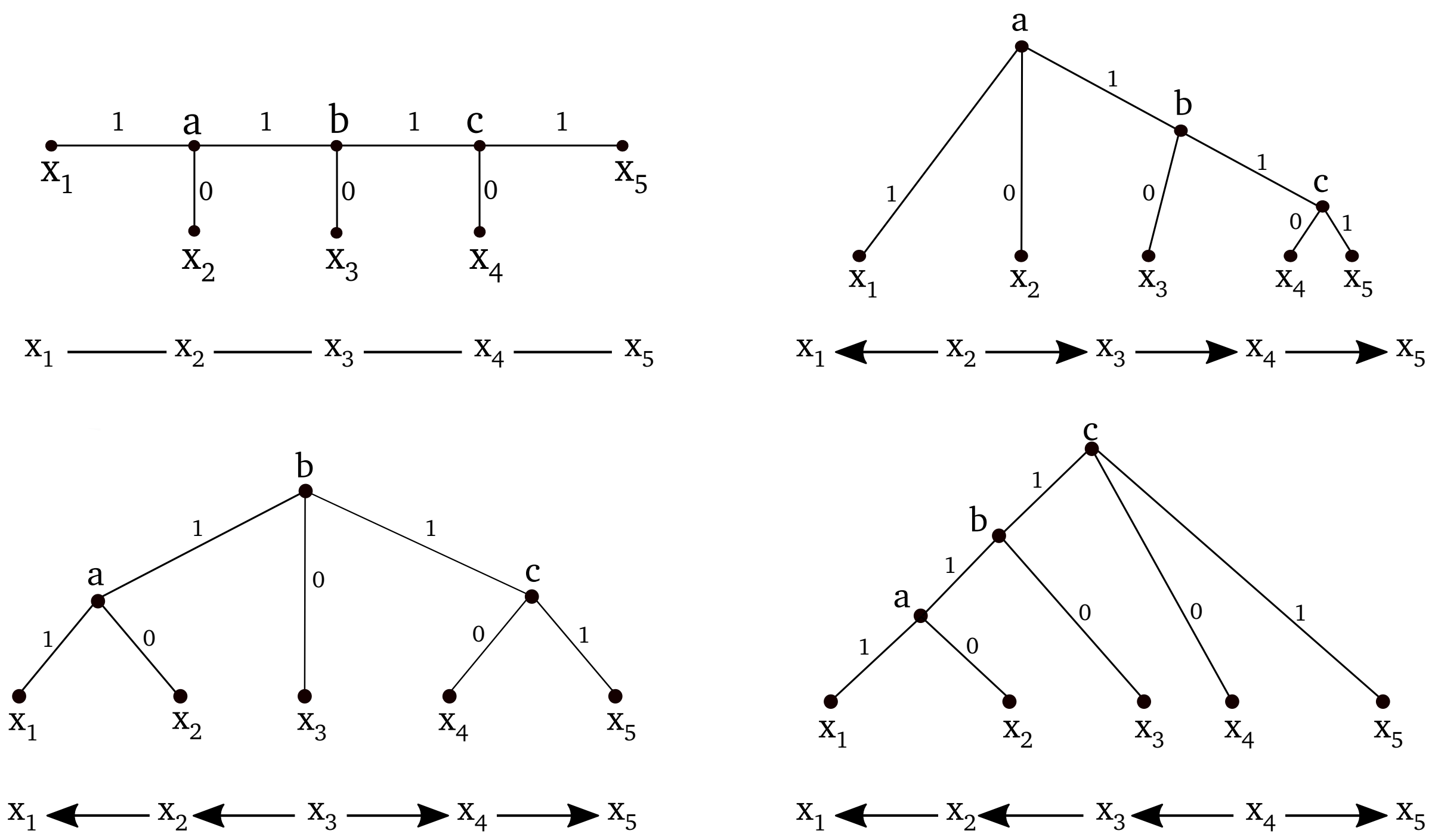

Figure 5.4: Placing the root. The tree $T$ in the upper left is the unique minimally resolved tree that explains $G(\overset{1}{\sim})/\overset{0}{\sim}$ (shown below $T$). Each of the tree inner vertices $a$, $b$, or $c$ of $T$ can be chosen as the root, giving rise to three distinct relations $\overset{1}{\rightharpoonup}$. For the "siblings" in the unrooted tree $x_1, x_2$ as well as $x_4, x_5$ it holds that $x_2 \overset{1}{\rightharpoonup} x_1$ and $x_4 \overset{1}{\rightharpoonup} x_5$ for all three distinct relations. Thus, there are uniquely determined parts of $\overset{1}{\rightharpoonup}$ conveyed by the information of $\overset{0}{\sim}$ and $\overset{1}{\sim}$ only.

Since there are two distinct last common ancestors, we have $u \neq v = lca(yz)$. Therefore, the triple $\mathsf{yz|x}$ must be displayed by $T(Q)$. From $y \overset{1}{\rightharpoonup} z$ we know that $\mathbb{P}(v, y)$ only contains 0-edges and $\mathbb{P}(v, z)$ contains exactly one 1-edge; $x \overset{1}{\rightharpoonup} y$ implies that $\mathbb{P}(x, u)$ contains only 0-edges. Moreover, since $\mathbb{P}(x, y) = \mathbb{P}(x, u) \cup \mathbb{P}(u, v) \cup \mathbb{P}(v, y)$ and $x \overset{1}{\rightharpoonup} y$, the path $\mathbb{P}(u, v)$ must contain exactly one 1-edge. □

**Lemma 5.5.2.** *In Case (ii), there is a unique tree on the three vertices $x, y, z$ with single root $\rho$ displayed by any least resolved tree $T(Q)$ that explains $Q$. Moreover, the path $\mathbb{P}(\rho, y)$ contains only 0-edges, while the other paths $\mathbb{P}(\rho, x)$ and $\mathbb{P}(\rho, z)$ must both contain exactly one 1-edge.*

*Proof.* Assume for contradiction that there is a least resolved tree $T(Q)$ that displays $\mathsf{xy|z}$, $\mathsf{yz|x}$, or $\mathsf{xz|y}$.

The choice of $\mathsf{xy|z}$ implies $u = lca(xy) \neq lca(xz) = lca(yz) = v$. Since $y \overset{1}{\rightharpoonup} x$

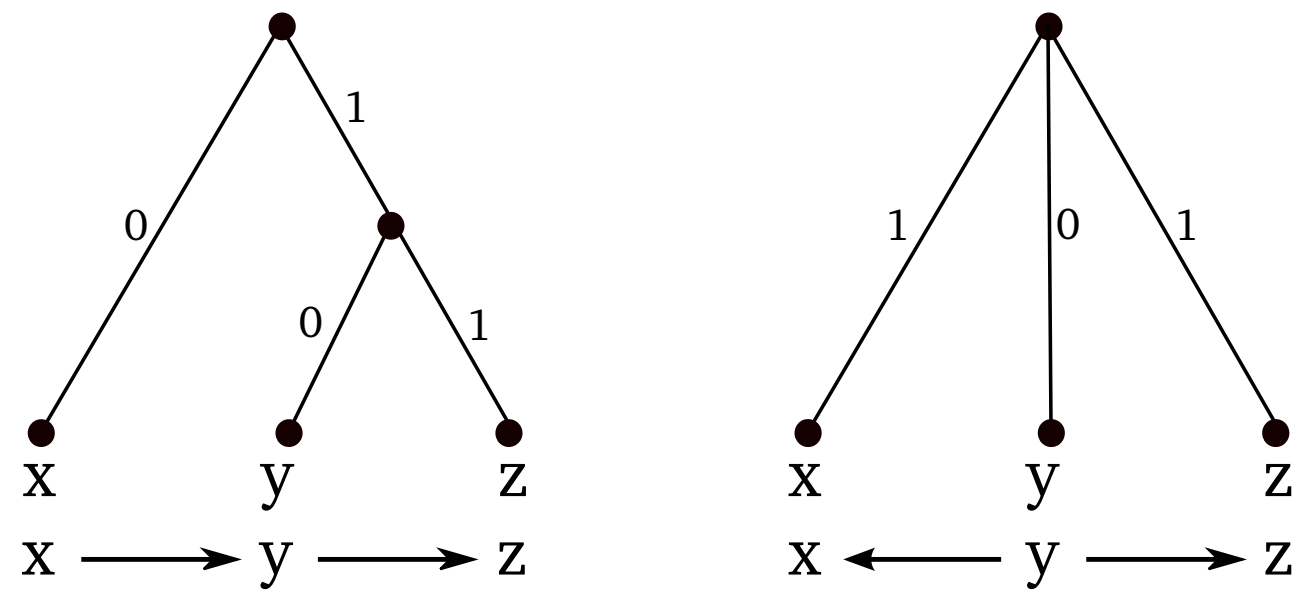

Figure 5.5: There are only two possibilities for induced connected subgraphs $H$ in $G(\overset{1}{\rightharpoonup})/\overset{0}{\sim}$ on three vertices, cf. Lemma 5.5.1, 5.5.2, and 5.5.3. Each of the distinct induced subgraphs imply a unique least resolved tree with a unique labeling.

and $y \overset{1}{\rightharpoonup} z$, $\mathbb{P}(v,y) \subsetneq \mathbb{P}(u,y)$ contain only 0-edges, while $\mathbb{P}(u,x)$ and $\mathbb{P}(v,z)$ each must contain exactly one 1-edge, respectively. This leads to a tree $T'$ that yields the correct $\overset{1}{\rightharpoonup}$-relation. However, this tree is not least resolved. By contracting the path $\mathbb{P}(u,v)$ to a single vertex $\rho$ and maintaining the labels on $\mathbb{P}(\rho,x)$, $\mathbb{P}(\rho,y)$, and $\mathbb{P}(\rho,z)$ we obtain the desired labeled least resolved tree with single root.

For the triple $\mathsf{yz|x}$ the existence of the unique, but not least resolved tree can be shown by the same argument with exchanged roles of $x$ and $y$.

For the triple $\mathsf{xz|y}$ we $u = lca(xy) = lca(yz) \neq lca(xz) = v$. From $y \overset{1}{\rightharpoonup} x$ and $y \overset{1}{\rightharpoonup} z$ we see that both paths $\mathbb{P}(u,x) = \mathbb{P}(u,v) \cup \mathbb{P}(v,x)$ and $\mathbb{P}(u,z) = \mathbb{P}(u,v) \cup \mathbb{P}(v,z)$ contain exactly one 1-edge, while all edges in $\mathbb{P}(u,y)$ are labeled 0. There are two cases: (1) The path $\mathbb{P}(u,v)$ contains this 1-edge, which implies that both paths $\mathbb{P}(v,x)$ and $\mathbb{P}(v,z)$ contain only 0-edges. But then $x \overset{0}{\sim} z$, a contradiction to $\overset{0}{\sim}$ being discrete. (2) The path $\mathbb{P}(u,v)$ contains only 0-edges, which implies that each of the paths $\mathbb{P}(v,x)$ and $\mathbb{P}(v,z)$ contain exactly one 1-edge. Again, this leads to a tree that yields the correct $\overset{1}{\rightharpoonup}$-relation, but it is not least resolved. By contracting the path $\mathbb{P}(u,v)$ to a single vertex $\rho$ and maintaining the labels on $\mathbb{P}(\rho,x)$, $\mathbb{P}(\rho,y)$, and $\mathbb{P}(\rho,z)$ we obtain the desired labeled least resolved tree with single root. □

**Lemma 5.5.3.** *Case (iii) cannot occur.*

*Proof.* Let $x, y, z \in V(Q)$ such that $xy, zy \in E$ and thus, $x \overset{1}{\rightharpoonup} y$ and $z \overset{1}{\rightharpoonup} y$. Hence, in the rooted tree that explains this relationship we have the following situation: All edges along $\mathbb{P}(u,x)$ are labeled 0; exactly one edge along $\mathbb{P}(u,y)$ is labeled 1, where $u = lca(x,y)$; all edges along $\mathbb{P}(v,z)$ are labeled 0, and exactly one edge along $\mathbb{P}(v,y)$ is labeled 1, where $v = lca(y,z)$. Clearly, $lca(x,y,z) \in \{u,v\}$. If $u = v$, then all edges in $\mathbb{P}(u,x)$ and $\mathbb{P}(u,z)$ are labeled 0, implying that $x \overset{0}{\sim} y$, contradicting that $\overset{0}{\sim}$ is discrete.

Now assume that $u = lca(x,y) \neq v = lca(y,z)$. Hence, one of the triples $\mathsf{xy|z}$ or $\mathsf{yz|x}$ must be displayed by $T(Q)$. W.l.o.g., we can assume that $\mathsf{yz|x}$ is displayed, since the case $\mathsf{xy|z}$ is shown analogously by interchanging the role of $x$ and $z$. Thus, $lca(x,y,z) = lca(x,y) = u \neq lca(yz) = v$. Hence, $\mathbb{P}(u,y) = \mathbb{P}(u,v) \cup \mathbb{P}(v,y)$. Since $z \overset{1}{\rightharpoonup} y$, the path $\mathbb{P}(v,y)$ contains a single 1-edge and $\mathbb{P}(v,z)$ contains only 0-edges. Therefore, the paths $\mathbb{P}(u,x)$ and $\mathbb{P}(u,v)$ contain only 0-edges, since $x \overset{1}{\sim} y$. Since $\mathbb{P}(x,z) = \mathbb{P}(x,u) \cup \mathbb{P}(u,v) \cup \mathbb{P}(v,z)$ and all edges along $\mathbb{P}(u,x)$, $\mathbb{P}(u,v)$ and $\mathbb{P}(v,z)$ are labeled 0, we obtain $x \overset{0}{\sim} z$, again a contradiction. □

Taken together, we obtain the following immediate implication:

**Corollary 5.5.1.** The graph $G(\overset{1}{\rightharpoonup})/\overset{0}{\sim}$ does not contain a pair of edges of the form $xv$ and $yv$.

Recall that the connected components $\underline{Q}$ in $G(\overset{1}{\rightharpoonup})/\overset{0}{\sim}$ are trees. By Cor. 5.5.1, $Q$ must be composed of distinct paths that "point away" from each other. In other words, let $P$ and $P'$ be distinct directed path in $Q$ that share a vertex $v$, then it is never the case that there is an edge $xv$ in $P$ and an edge $yv$ in $P'$, that is, both edges "pointing" to the same vertex $v$. We first consider directed paths in isolation.

**Lemma 5.5.4.** *Let $Q$ be a connected component in $G(\overset{1}{\rightharpoonup})/\overset{0}{\sim}$ that is a directed path with $n \geq 3$ vertices labeled $x_1, \dots, x_n$ such that $x_i x_{i+1} \in E(Q)$, $1 \leq i \leq n-1$. Then the tree $T(Q)$ explaining $Q$ must display all triples in $\mathcal{R}_Q = \{\mathsf{x_i x_j | x_l} \mid i,j > l \geq 1\}$. Hence, $T(Q)$ must display $\binom{n}{3}$ triples and is therefore the unique (least resolved) binary rooted tree $(\dots(x_n, x_{n-1})x_{n-2})\dots)x_2)x_1$ that explains $Q$. Moreover, all interior edges in $T(Q)$ and the edge incident to $x_n$ are labeled* 1 *while all other edges are labeled* 0.

*Proof.* Let $Q$ be a directed path as specified in the lemma. We prove the statement by induction. For $n = 3$ the statement follows from Lemma 5.5.1. Assume the statement is true for $n = k$. Let $Q$ be a directed path with vertices $x_1, \dots, x_k, x_{k+1}$ and edges $x_i x_{i+1}$, $1 \leq i \leq k$ and let $T(Q)$ be a tree that explains $Q$. For the subpath $Q'$ on the vertices $x_2, \dots, x_{k+1}$ we can apply the induction hypothesis and conclude that $T(Q')$ must display the triples $\mathsf{x_i x_j | x_l}$ with $i,j > l \geq 2$ and that all interior edges in $T(Q')$ and the edge incident to $x_{k+1}$ are labeled 1 while all other edges are labeled 0. Since $T(Q)$ must explain in particular the subpath $Q'$ and since $T(Q')$ is fully resolved, we can conclude that $T(Q')$ is displayed by $T(Q)$ and that all edges in $T(Q)$ that are also in $T(Q')$ retain the same label as in $T(Q')$ .

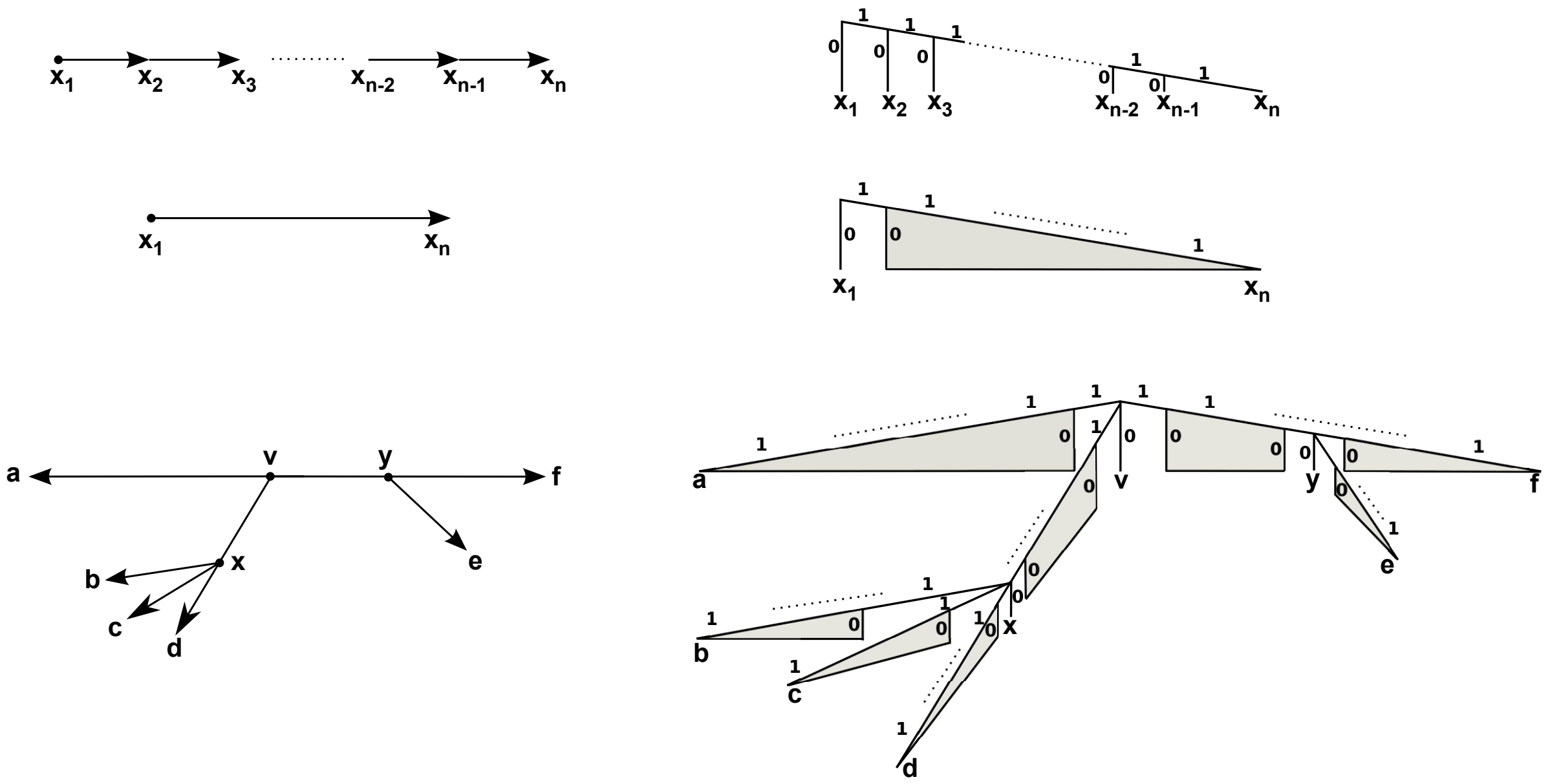

Figure 5.6: Connected components $Q$ in $G(\stackrel{1}{\rightharpoonup})/\stackrel{0}{\sim}$ are trees composed of paths pointing away from each other. On the top, a directed path $\mathbb{P}(x_1, x_n)$ together with its unique least resolved tree representation is shown. In the middle, an abstract picture of the path $\mathbb{P}(x_1, x_n)$ and its tree is given. Bottom, a larger example of a connected component $Q$ in $G(\stackrel{1}{\rightharpoonup})/\stackrel{0}{\sim}$ and its tree-representation are sketched.

Thus $T(Q)$ displays in particular the triples $\mathsf{x_i x_j | x_l}$ with $i, j > l \geq 2$. By Lemma 5.5.1, and because there are edges $x_1x_2$ and $x_2x_3$, we see that $T(Q)$ must also display $\mathsf{x_2x_3|x_1}$. Take any triple $\mathsf{x_3x_j|x_2}$, $j > 3$. Application of the triple-inference rules shows that any tree that displays $\mathsf{x_2x_3|x_1}$ and $\mathsf{x_3x_j|x_2}$ must also display $\mathsf{x_2x_j|x_1}$ and $\mathsf{x_3x_j|x_1}$. Hence, $T(Q)$ must display these triples. Now we apply the same argument to the triples $\mathsf{x_2x_j|x_1}$ and $\mathsf{x_ix_j|x_2}$, $i, j > 2$ and conclude that in particular, the triple $\mathsf{x_ix_j|x_1}$ must be displayed by $T(Q)$ and thus, the the entire set of triples $\{\mathsf{x_ix_j|x_l}\colon i, j > l \geq 1\}$. Hence, there are $\binom{n}{3}$ triples and thus, the set of triples that needs to be displayed by $T(Q)$ is strictly dense. Making use of a technical result from [43, Suppl. Material], we obtain that $T(Q)$ is the unique binary tree $(\dots(x_n, x_{n-1})x_{n-2})\dots)x_2)x_1$. Now it is an easy exercise to verify that the remaining edge containing $x_1$ must be labeled 0, while the interior edge not contained in $T(Q')$ must all be 1-edges. $\square$

If $Q$ is connected but not a simple path, it is a tree composed of the paths pointing away from each other as shown in Fig. 5.6. It remains to show how to connect the distinct trees that explain these paths to obtain a tree $T(Q)$ for $Q$. To this end, we show first that there is a unique vertex $v$ in $Q$ such that no edge ends in $v$.

**Lemma 5.5.5.** *Let $Q$ be a connected component in $G(\overset{1}{\rightharpoonup})/\overset{0}{\sim}$. Then there is a unique vertex $v$ in $Q$ such that there is no edge $xv \in E(Q)$.*

*Proof.* Corollary 5.5.1 implies that for each vertex $v$ in $Q$ there is at most one edge $xv \in E(Q)$. If for all vertices $w$ in $Q$ we would have an edge $xw \in E(Q)$, then $Q$ contains cycles, contradicting the tree structure of $Q$. Hence, there is at least one vertex $v$ so that there is no edge of the form $xv \in E(Q)$. Assume there are two vertices $v, v'$ so that there are no edges of the form $xv, yv'$, then all edges incident to $v, v'$ are of the form $vx, v'y$. However, in this case, the unique path connecting $v, v'$ in $Q$ must contain two edges of the form $aw$, $bw$; a contradiction to Corollary 5.5.1. Thus, there is exactly one vertex $v$ in $Q$ such that there is no edge $xv \in E(Q)$. □

By Lemma 5.5.5, for each connected component $Q$ of $G(\overset{1}{\rightharpoonup})/\overset{0}{\sim}$ there is a unique vertex $v_Q$ s.t. all edges incident to $v_Q$ are of the form $v_Qx$. That is, all directed paths that are maximal w.r.t. inclusion start in $v_Q$. Let $\mathcal{P}_Q$ denote the sets of all such maximal paths. Thus, for each path $P \in \mathcal{P}_Q$ there is the triple set $\mathcal{R}_P$ according to Lemma 5.5.4 that must be displayed by any tree that explains also $Q$. Therefore, $T(Q)$ must display all triples in $\mathcal{R}_Q = \cup_{P\in\mathcal{P}_Q}\mathcal{R}_P$.

The underlying undirected graph $\underline{G(\overset{1}{\rightharpoonup})/\overset{0}{\sim}}$ is isomorphic to $G(\overset{1}{\sim})/\overset{0}{\sim}$. Thus, with Algorithm 1, one can similar to the unrooted case, first construct the tree $T(\underline{Q})$ and then set the root $\rho_Q = v'_Q$ to obtain $T(Q)$. It is easy to verify that this tree $T(Q)$ displays all triples in $\mathcal{R}_Q$. Moreover, any edge-contradiction in $T(Q)$ leads to the loss of an input triple $\mathcal{R}_Q$ and in particular, to a wrong pair of vertices w.r.t. $\overset{1}{\rightharpoonup}$ or $\overset{0}{\sim}$. Thus, $T(Q)$ is a least resolved tree for $\mathcal{R}_Q$ and therefore, a least resolved tree that explains $Q$.

We summarize these arguments in

**Corollary 5.5.2.** Let $Q$ be a connected component in $G(\overset{1}{\rightharpoonup})/\overset{0}{\sim}$. Then a tree $T(Q)$ that explains $Q$ can be obtained from the unique minimally resolved tree $T(\underline{Q})$ that explains $\underline{Q}$ by choosing the unique vertex $v$ where all edges incident to $v$ are of the form $vx$ as the root $\rho_Q$.

If $G(\overset{1}{\rightharpoonup})/\overset{0}{\sim}$ is disconnected, one can apply Algorithm 2, to obtain the tree $T(G(\overset{1}{\sim})/\overset{0}{\sim})$ and then chose either one of the vertices $\rho_Q$ or the vertex $z_T$ as root to obtain $T(G(\overset{1}{\rightharpoonup})/\overset{0}{\sim})$, in which case all triples of $\mathcal{R}_{G(\overset{1}{\rightharpoonup})/\overset{0}{\sim}} = \cup_Q\mathcal{R}_Q$ are displayed. Again, it is easy to verify that any edge-contradiction leads to a wrong pair of vertices in $\overset{1}{\rightharpoonup}$ or $\overset{0}{\sim}$. Thus, $T(G(\overset{1}{\rightharpoonup})/\overset{0}{\sim})$ is a least resolved tree for $\mathcal{R}_{G(\overset{1}{\rightharpoonup})/\overset{0}{\sim}}$.

To obtain uniqueness of minimally resolved trees one can apply similar arguments as in the proofs of Theorems 5.4.2 and 5.4.3. This yields the following characterization:

**Theorem 5.5.3.** *Let $Q_1, \dots Q_k$ be the connected components in $G(\overset{1}{\rightharpoonup})/\overset{0}{\sim}$. Up to the choice of the vertices $q_i'$ in Line 11 of Alg. 2 for the construction of $T(\underline{Q_i})$ and the choice of the root $\rho \in \{\rho_{Q_1}, \dots \rho_{Q_k}, z\}$, the tree $T^* = T(G(\overset{1}{\rightharpoonup})/\overset{0}{\sim}))$ is the unique minimally resolved tree that explains $G(\overset{1}{\rightharpoonup})/\overset{0}{\sim}$.*

## 5.6 Mix of symmetric and anti-symmetric relations

In real data, e.g., in the application to mitochondrial genome arrangements, one can expect that the known relationships are in part directed and in part undirected. Such data are naturally encoded by a relation $\overset{1}{\rightharpoonup}$ with directional information and a relation $\overset{1\star}{\sim}$ comprising the set of pairs for which it is unknown whether one of $x \overset{1}{\rightharpoonup} y$ and $y \overset{1}{\rightharpoonup} x$ or $x \overset{1}{\sim} y$ are true. Here, $\overset{1}{\sim}$ is a subset of $\overset{1\star}{\sim}$. The disjoint union $\overset{1\star}{\sim} \uplus \overset{1}{\rightharpoonup}$ of these two parts can be seen as refinement of a corresponding symmetrized relation $x \overset{1}{\sim} y$. Ignoring the directional information one can still construct the tree $T(G(\overset{1}{\sim})/\overset{0}{\sim})$. In general there will be less information of the placement of the root in $T(G(\overset{1\star}{\sim} \cup \overset{1}{\rightharpoonup})/\overset{0}{\sim})$ than with a fully directed edge set.

In what follows, we will consider all edges of $G(\overset{1\star}{\sim} \cup \overset{1}{\rightharpoonup})/\overset{0}{\sim}$ to be directed, that is, for a symmetric pair $(a, b) \in \overset{1\star}{\sim}$ we assume that both directed edges $ab$ and $ba$ are contained in $G(\overset{1\star}{\sim} \cup \overset{1}{\rightharpoonup})/\overset{0}{\sim}$. Still, for any connected component $Q$ the underlying undirected graph $\underline{Q}$ is a tree. Given a component $Q$ we say that a directed edge $xy \in E(Q)$ *points away from the vertex* $v$ if the unique path in $\underline{Q}$ from $v$ to $x$ does not contain $y$. In this case the path from $v$ to $y$ must contain $x$. Note that in this way we defined "pointing away from $v$" not only for the edges incident to $v$, but for all directed edges. A vertex $v$ is a *central vertex* if, for any two distinct vertices $x, y \in V$ that form an edge in $T$, either $xy$ or $yx$ in $T$ points away from $v$.

As an example consider the tree $a \leftarrow b \rightarrow c \leftrightarrow d \rightarrow e$. There is only the edge $bc$ containing $b$ and $c$. However, $bc$ does not point away from vertex $d$, since the unique path from $d$ to $b$ contains $c$. Thus $d$ is not central. On the other hand, $b$ is a central vertex. The only possibility in this example to obtain a valid relation $\overset{1}{\rightharpoonup}$ that can be displayed by rooted 0/1-edge-labeled tree is provided by removing the edge $dc$, since otherwise Cor. 5.5.1 would be violated.

In the following, for given relations $\overset{1\star}{\sim}$ and $\overset{1}{\rightharpoonup}$ we will denote with $\overset{1\star}{\rightharpoonup}$ a relation that contains $\overset{1}{\rightharpoonup}$ and exactly one pair, either $(x, y)$ or $(y, x)$, from $\overset{1\star}{\sim}$.

**Lemma 5.6.1.** *For a given graph $G(\overset{1\star}{\sim} \cup \overset{1}{\rightharpoonup})/ \overset{0}{\sim}$ the following statements are equivalent:*

*(i) There is a relation $\overset{1\star}{\rightharpoonup}$ that is the antisymmetric single-1-relation of some 0/1-edge-labeled tree.*

*(ii) There is a central vertex in each connected component $Q$ of $G(\overset{1\star}{\sim} \cup \overset{1}{\rightharpoonup})/ \overset{0}{\sim}$.*

*Proof.* If there is a relation $\overset{1\star}{\rightharpoonup}$ that can be displayed by a rooted 0/1-edge-labeled tree, then $G(\overset{1\star}{\rightharpoonup})/ \overset{0}{\sim}$ consists of connected components $Q$ where each connected component is a tree composed of maximal directed paths that point away from each other. Hence, for each connected component $Q$ there is the unique vertex $v_Q$ such that all edges incident to $v_Q$ are of the form $v_Q x$ and, in particular, $v_Q$ is a central vertex $v_Q$ in $Q$ and thus, in $G(\overset{1\star}{\sim} \cup \overset{1}{\rightharpoonup})/ \overset{0}{\sim}$.

Conversely, assume that each connected component $Q$ has a central vertex $v_Q$. Hence, one can remove all edges that do not point away from $v_Q$ and hence obtain a connected component $Q'$ that is still a tree with $V(Q) = V(Q')$ so that all maximal directed paths point away from each other and in particular, start in $v_Q$. Thus, for the central vertex $v_Q$ all edges incident to $v_Q$ are of the form $v_Q x$. Since $Q'$ is now a connected component in $G(\overset{1\star}{\rightharpoonup})/ \overset{0}{\sim}$, we can apply Cor. 5.5.2 to obtain the tree $T(Q')$ and Thm. 5.5.3 to obtain $T(G(\overset{1\star}{\rightharpoonup})/ \overset{0}{\sim})$. □

The key consequence of this result is the following characterization of the constraints on the possible placements of the root.

**Corollary 5.6.1.** Let $Q$ be a connected component in $G(\overset{1\star}{\sim} \cup \overset{1}{\rightharpoonup})/ \overset{0}{\sim}$ and let $T(\underline{Q})$ be the unique least resolved tree that explains the underlying undirected graph $\underline{Q}$. Then each copy $v'$ of a vertex $v$ in $Q$ can be chosen to be the root in $T(\underline{Q})$ to obtain $T(Q)$ if and only if $v$ is a central vertex in $Q$.

## 5.7 Concluding Remarks

In this contribution we have introduced a class of binary relations deriving in a natural way from edge-labeled trees. This construction has been inspired by the conceptually similar class of relations induced by vertex-labeled trees [41]. The latter have co-graph structure and are closely related to orthology and paralogy [37, 52, 43, 40]. Defining $x \sim y$ whenever at least one 1-edge lies along the path from $x$ to $y$ is related to the notion of xenology: the edges labeled 1 correspond to horizontal gene transfer events, while the 0-edge encode vertical transmission. In its simplest setting, this idea can also be combined with vertex labels, leading to the directed analog of co-graphs [39]. Here, we have explored an even simpler special case: the existence of a single 1-label along the connecting path, which captures the structure of rare event data as we have discussed in the introduction. We have succeeded here in giving a complete characterization of the relationships between admissible relations, which turned out to be forests, and the underlying phylogenetic tree. Moreover, for all such cases we gave polynomial-time algorithms to compute the trees explaining the respective relation.

The characterization of single-1 relations is of immediate relevance for the use of rare events in molecular phylogenetics. In particular, it determines lower bounds on the required data: if too few events are known, many taxa remain in $\overset{0}{\sim}$ relation and thus unresolved. On the other hand, if taxa are spread too unevenly, $G(\overset{1}{\rightharpoonup})/\overset{0}{\sim}$ will be disconnected, consisting of connected components separated by multiple events. The approach discussed here, of course, is of practical use only if the available event data are very sparse. If taxa are typically separated by multiple events, classical phylogenetic approaches, i.e., maximum parsimony, maximum likelihood, or Bayesian methods, will certainly be preferable. An advantage of using the single-1 relation is that it does not require the product structure of independent characters implicit in the usual, sequence or character-based methods, nor does it require any knowledge of the algebraic properties of the underlying operations as, e.g., in phylogenetic reconstruction from (mitochondrial) gene order rearrangements. This begs the question under which conditions on the input data identical results are obtained from the direct translation of the single-1 relation and maximum parsimony on character data or break point methods for genome rearrangments.

A potentially useful practical application of our results is a new strategy to incorporate rare event data into conventional phylogenetic approaches. Trees obtained from a single-1 relation will in general be poorly resolved. Nevertheless, they determine some monophyletic groups (in the directed case) or a set of split (in the undirected case) that have to be present in the true phylogeny. Many of the software tools commonly used in molecular phylogenetic can incorporate this type of constraints. Assuming the there is good evidence that the rare events are homoplasy-free, the $T(G(\overset{1}{\rightharpoonup})/\overset{0}{\sim})$ represents the complete information from the rare events that can be used constrain the tree reconstruction process.

The analysis presented here makes extensive use of the particular properties of the single-1 relation and hence does not seem to generalize easily to other interesting cases.

Horizontal gene transfer, for example, is expressed naturally in terms of the "at-least-one-1" relation $\rightsquigarrow$. It is worth noting that $\rightsquigarrow$ also has properties (L1) and (L2) and hence behaves well w.r.t. contraction of the underlying tree and restriction to subsets of leaves. Whether this is sufficient to obtain a complete characterization remains an open question.

Several general questions arise naturally. For instance, is there a characterization of admissible relations in terms of forbidden subgraphs graphs or minors? For instance, the relation $\stackrel{1}{\rightharpoondown}$ / $\stackrel{0}{\sim}$ is characterized in terms of the forbidden subgraph $x \rightarrow v \leftarrow y$. Hence, it would be of interest, whether such characterizations can be derived for arbitrary relations $\stackrel{k}{\rightharpoondown}$ or for $\rightsquigarrow$. If so, can these forbidden substructures be inferred in a rational manner from properties of vertex and/or edge labels along the connecting paths in the explaining tree? Is this the case at least for labels and predicates satisfying (L1) and (L2)?

# Chapter 6

# Reconstructing unrooted phylogenetic trees from symbolic ternary metrics

In 1998, Böcker and Dress presented a 1-to-1 correspondence between symbolically dated rooted trees and symbolic ultrametrics. We consider the corresponding problem for unrooted trees. More precisely, given a tree $T$ with leaf set $X$ and a proper vertex coloring of its interior vertices, we can map every triple of three different leaves to the color of its median vertex. We characterize all ternary maps that can be obtained in this way in terms of 4- and 5-point conditions, and we show that the corresponding tree and its coloring can be reconstructed from a ternary map that satisfies those conditions. Further, we give an additional condition that characterizes whether the tree is binary, and we describe an algorithm that reconstructs general trees in a bottom-up fashion.

## 6.1 Introduction

A phylogenetic tree is a rooted or unrooted tree where the leaves are labeled by some objects of interest, usually taxonomic units (taxa) like species. The edges have a positive edge length, thus the tree defines a metric on the taxa set. It is a classical result in phylogenetics that the tree can be reconstructed from this metric, if it is unrooted or ultrametric. The latter means that the tree is rooted and all taxa are equally far away from the root. An ultrametric tree is realistic whenever the edge lengths are proportional to time and the taxa are species that can be observed in the present. In an ultrametric tree, the distance between two taxa is twice of the distance between each of the taxa and their last common ancestor (lca), hence pairs of taxa with the same lca must have the same distance. For three taxa $x, y, z$, it follows that there is no unique maximum within their three pairwise distances, thus we have $d(x, y) \leq \max\{d(x, z), d(y, z)\}$. This 3-point condition turns out to be sufficient for a metric to be ultrametric, too, and it is the key for reconstructing ultrametric trees from their distances. In 1995, Bandelt and Steel [4] observed that the complete ordering of the real numbers is not necessary to reconstruct trees, and they showed that the real-valued distances can be replaced by maps from the pairs of taxa into a cancellative abelian monoid. Later, Böcker and Dress [6] pushed this idea to the limit by proving that the image set of the symmetric map does not need any structure at all (see Section 6.2 for details). While this result is useful for understanding how little information it takes to reconstruct an ultrametric phylogenetic tree, it was not until recently that it turned out to have some practical applications. In 2013, Hellmuth et al. [38] found an alternative characterization of symbolic ultrametrics in terms of cographs and showed that, for perfect data, phylogenetic trees can be reconstructed from orthology information. By adding some optimization tools, this concept was then applied to analyze real data [42].

Motivated by the practical applicability of symbolic ultrametrics, we are considering their unrooted version. However, in an unrooted tree there is in general no interior vertex associated to a pair of taxa that would correspond to the last common ancestor in a rooted tree. Instead, there is a median associated to every set of three taxa that represents, for every possible rooting of the tree, a last common ancestor of at least two of the three taxa. Therefore, we consider ternary maps from the triples of taxa into an image set without any structure. We will show that an unrooted phylogenetic tree with a proper vertex coloring can be reconstructed from the function that maps every triple of taxa to the color of its median.

In order to apply our results to real data, we need some way to assign a state to every set of three taxa, with the property that 3-sets with the same median will usually have the same state. For symbolic ultrametrics, the first real application was found 15 years after the development of the theory. In addition to the hope that something similar happens with symbolic ternary metrics, we have some indication that they can be useful to construct unrooted trees from orthology relations (see Section 6.6 for details).

Consider an unrooted tree $T$ with vertex set $V$, edge set $E$, and leaf set $X$, and a dating map $t : V \to M^{\odot}$, where $M^{\odot} = M \cup \{\odot\}$ such that $t(x) = \odot$ for all $x \in X$, and $t(v_1) \neq t(v_2)$ if $v_1 v_2 \in E$.

For any $S = \{x, y\} \in \binom{V}{2}$ there is a unique path $[x, y]$ with end points $x$ and $y$, and for any 3-set $S = \{x, y, z\} \in \binom{V}{3}$ there is a unique *triple point* or *median* $\text{med}(x, y, z)$ such that $[x, y] \cap [y, z] \cap [x, z] = \{\text{med}(x, y, z)\}$. Putting $[x, x] = \{x\}$, the definition also works, if some or all of $x, y$ and $z$ equal.

Given a phylogenetic tree $T$ on $X$ and a dating map $t : V \to M^{\odot}$, we can define the symmetric symbolic ternary map $d_{(T;t)} : X \times X \times X \to M^{\odot}$ by $d_{(T;t)}(x, y, z) = t(\mathrm{med}(x, y, z))$.

In this set-up, our question can be phrased as follows: Suppose we are given an arbitrary symbolic ternary map $\delta : X \times X \times X \to M^{\odot}$, can we determine if there is a pair $(T; t)$ for which $d_{(T:t)}(x, y, z) = \delta(x, y, z)$ holds for all $x, y, z \in X$?

The rest of this paper is organized as follows. In Section 6.1.1, we present the basic and relevant concepts used in this paper. In Section 6.2 we recall the one-to-one correspondence between symbolic ultrametrics and symbolically dated trees, and introduce our main results Theorem 6.2.4 and Theorem 6.2.5.

In Section 6.3 we give the proof of Theorem 6.2.4. In order to prove our main result, we first introduce the connection between phylogenetic trees and quartet systems on $X$ in Subsection 6.3.1. Then we use a graph representation to analyze all cases of the map $\delta$ for 5-taxa subsets of $X$ in Subsection 6.3.2.

In Section 6.4 we use a similar method to prove Theorem 6.2.5, which gives a sufficient and necessary condition to reconstruct a binary phylogenetic tree on $X$.

In Section 6.5, we give a criterion to identify all pseudo-cherries of the underlying tree from a symbolic ternary metric. This result makes it possible to reconstruct the tree in a bottom-up fashion.

In the last section we discuss some open questions and future work.

### 6.1.1 Preliminaries

We introduce the relevant basic concepts and notation. Unless stated otherwise, we will follow the monographs [73] and [18].

In the remainder of this paper, $X$ denotes a finite set of size at least three.

An *(unrooted) tree* $T = (V, E)$ is an undirected connected acyclic graph with vertex set $V$ and edge set $E$. A vertex of $T$ is a *leaf* if it is of degree 1, and all vertices with degree at least two are *interior* vertices.

A *rooted tree* $T = (V, E)$ is a tree that contains a distinguished vertex $\rho_T \in V$ called the *root*. We define a partial order $\preceq_T$ on $V$ by setting $v \preceq_T w$ for any two vertices $v, w \in V$ for which $v$ is a vertex on the path from $\rho_T$ to $w$. In particular, if $v \preceq_T w$ and $v \neq w$ we call $v$ an *ancestor* of $w$.

An unrooted *phylogenetic tree* $T$ on $X$ is an unrooted tree with leaf set $X$ that does not contain any vertex of degree 2. It is *binary*, if every interior vertex has degree 3.

A rooted phylogenetic tree $T$ on $X$ is a rooted tree with leaf set $X$ that does not contain any vertices with in- and out-degree one, and whose root $\rho_T$ has in-degree zero. For a set $A \subseteq X$ with cardinality at least 2, we define *the last common ancestor* of $A$, denoted by $\mathrm{lca}_T(A)$, to be the unique vertex in $T$ that is the greatest lower bound of $A$ under the partial order $\preceq_T$. In case $A = \{x, y\}$ we put $\mathrm{lca}_T(x, y) = \mathrm{lca}_T(\{x, y\})$.

Given a set $Q$ of four taxa $\{a, b, c, d\}$, there exist always exactly three partitions into two pairs: $\{\{a, b\}, \{c, d\}\}$,$\{\{a, c\}, \{b, d\}\}$ and $\{\{a, d\}, \{b, c\}\}$. These partitions are called *quartets*, and they represent the three non-isomorphic unrooted binary trees with leaf set $Q$. These trees are usually called quartet trees, and they – as well as the corresponding quartets –are symbolized by $ab|cd, ac|bd, ad|bc$ respectively. We use $Q(X)$ to denote the set of all quartets with four taxa in $X$. A phylogenetic tree $T$ on $X$ *displays* a quartet $ab|cd \in Q(X)$, if the path from $a$ to $b$ in $T$ is vertex-disjoint with the path from $c$ to $d$. The collection of all quartets that are displayed by $T$ is denoted by $Q_T$.

Let $M$ be a non-empty finite set, $\odot$ denotes a special element not contained in $M$, and $M^{\odot} := M \cup \{\odot\}$. Note that in biology the symbol $\odot$ corresponds to a "non-event" and is introduced for purely technical reasons [38]. A *symbolic ternary map* is a mapping from $X \times X \times X$ to $M^{\odot}$. Suppose we have a symbolic ternary map $\delta : X \times X \times X \to M^{\odot}$, we say $\delta$ is *symmetric* if the value of $\delta(x, y, z)$ is only related to the set $\{x, y, z\}$ but not on the ordering of $x, y, z$, i.e., if $\delta(x, y, z) = \delta(y, x, z) = \delta(z, y, x) = \delta(x, z, y)$ for all $x, y, z \in X$. For simplicity, if a map $\delta : X \times X \times X \to M^{\odot}$ is symmetric, then we can define $\delta$ on the set $\{x, y, z\}$ to be $\delta(x, y, z)$.

Let $S$ be a set, we define $|S|$ to be the number of elements in $S$.

## 6.2 Symbolic ultrametrics and our main results

In this section, we first recall the main result concerning symbolic ultrametrics by Böcker and Dress [6].

Suppose $\delta : X \times X \to M^{\odot}$ is a map. We call $\delta$ a *symbolic ultrametric* if it satisfies the following conditions:

(U1) $\delta(x, y) = \odot$ if and only if $x = y$;

(U2) $\delta(x, y) = \delta(y, x)$ for all $x, y \in X$, i.e., $\delta$ is symmetric;

(U3) $|\{\delta(x, y), \delta(x, z), \delta(y, z)\}| \leq 2$ for all $x, y, z \in X$; and

(U4) there exists no subset $\{x, y, u, v\} \in \binom{X}{4}$ such that
$\delta(x, y) = \delta(y, u) = \delta(u, v) \neq \delta(y, v) = \delta(x, v) = \delta(x, u)$.

Now suppose that $T = (V, E)$ is a rooted phylogenetic tree on $X$ and that $t : V \to M^{\odot}$ is a map such that $t(x) = \odot$ for all $x \in X$. We call such a map $t$ a *symbolic dating map* for $T$; it is *discriminating* if $t(u) \neq t(v)$, for all edges $\{u, v\} \in E$. Given $(T, t)$, we associate the map $d_{(T;t)}$ on $X \times X$ by setting, for all $x, y \in X$, $d_{(T;t)}(x, y) = t(\text{lca}_T(x, y))$. Clearly $\delta = d_{(T;t)}$ satisfies Conditions (U1),(U2),(U3),(U4) and we say that $(T; t)$ is a *symbolic representation* of $\delta$. Böcker and Dress established in 1998 the following fundamental result which gives a 1-to-1 correspondence between symbolic ultrametrics and symbolic representations [6], i.e., the map defined by $(T, t) \mapsto d_{(T,t)}$ is a bijection from the set of symbolically dated trees into the set of symbolic ternary metrics.

**Theorem 6.2.1** (Böcker and Dress 1998 [6])**.** *Suppose $\delta : X \times X \to M^{\odot}$ is a map. Then there is a discriminating symbolic representation of $\delta$ if and only if $\delta$ is a symbolic ultrametric. Furthermore, up to isomorphism, this representation is unique.*

Similarly, we consider unrooted trees. Suppose that $T = (V, E)$ is an unrooted tree on $X$ and that $t : V \to M^{\odot}$ is a symbolic dating map, i.e., $t(x) = \odot$ for all $x \in X$, it is discriminating if $t(x) \neq t(y)$ for all $(x, y) \in E$. Given the pair $(T; t)$, we associate the map $\delta_{(T;t)}$ on $X \times X \times X$ by setting, for all $x, y, z \in X$, $\delta_{(T;t)}(x, y, z) = t(\text{med}(x, y, z))$.

Before stating our main results, we need the following definition:

**Definition 6.2.2** (*n*-*m* partitioned)**.** Suppose $\delta : X \times X \times X \to M^{\odot}$ is a symmetric map. We say that a subset $S$ of $X$ is *n-m partitioned (by $\delta$)*, if among all the 3-element subsets of $S$, there are in total 2 different values of $\delta$, and $n$ of those 3-sets are mapped to one value while all other $m$ 3-sets are mapped to the other value.

Note that $S$ can be $n$-$m$ partitioned, only when $\binom{|S|}{3} = m + n$.

**Definition 6.2.3** (symbolic ternary metrics)**.** We say $\delta : X \times X \times X \to M^{\odot}$ is a *symbolic ternary metric*, if the following conditions hold.

(1) $\delta$ is symmetric, i.e., $\delta(x, y, z) = \delta(y, x, z) = \delta(z, y, x) = \delta(x, z, y)$ for all $x, y, z \in X$.

(2) $\delta(x, y, z) = \odot$ if and only if $x = z$ or $y = z$ or $x = y$.

(3) for any distinct $x, y, z, u$ we have

$$|\{\delta(x, y, z), \delta(x, y, u), \delta(x, z, u), \delta(y, z, u)\}| \leq 2,$$

and when the equality holds then $\{x, y, z, u\}$ is 2-2 partitioned by $\delta$.

(4) there is no distinct 5-element subset $\{x, y, z, u, e\}$ of $X$ which is 5-5 partitioned by $\delta$.

We will refer to these conditions throughout the paper.

Our main result is:

**Theorem 6.2.4.** *There is a 1-to-1 correspondence between the discriminating symbolically dated phylogenetic trees and the symbolic ternary metrics on $X$.*

Let $\delta$ be a ternary symbolic ultrametric on $X$. Then we call $\delta$ *fully resolved*, if the following condition holds:

(*) If $|\{\delta(x, y, z), \delta(x, y, u), \delta(x, z, u), \delta(y, z, u)\}| = 1$, then there exists $e \in X$ such that $e$ can resolve $xyzu$. i.e., the set $\{x, y, z, u, e\}$ is 4-6 partitioned by $\delta$.

Now we can characterize ternary symbolic ultrametrics that correspond to binary phylogenetic trees:

**Theorem 6.2.5.** *There is a 1-to-1 correspondence between the discriminating symbolically dated binary phylogenetic trees and the fully resolved symbolic ternary metrics on $X$.*

## 6.3 Reconstructing a symbolically dated phylogenetic tree

The aim of this section is to prove Theorem 6.2.4.

### 6.3.1 Quartet systems

We will use quartet systems to prove Theorem 6.2.4. In 1981, Colonius and Schulze [14] found that, for a quartet system $Q$ on a finite taxa set X, there is a phylogenetic tree $T$ on $X$ such that $Q = Q_T$, if and only if certain conditions on subsets of $X$ with up to five elements hold. The following theorem (Theorem 3.7 in [18]) states their result.

A quartet system $Q$ is *thin*, if for every 4-subset $a, b, c, d \subseteq X$, at most one of the three quartets $ab|cd$, $ac|bd$ and $ad|bc$ is contained in $Q$. It is *transitive*, if for any 5 distinct elements $a, b, c, d, e \in X$, the quartet $ab|cd$ is in $Q$ whenever both of the quartets $ab|ce$ and $ab|de$ are contained in $Q$. It is *saturated*, if for any five distinct elements $a$ , $b$ , $c$ , $d$ , $e \in X$ with $ab|cd \in Q$, at least one of the two quartets $ae|cd$ and $ab|ce$ is also in $Q$.

**Theorem 6.3.1.** *A quartet system $Q \subseteq Q(X)$ is of the form $Q = Q(T)$ for some phylogenetic tree $T$ on $X$ if and only if $Q$ is thin, transitive and saturated.*

We can encode a phylogenetic tree on $X$ in terms of a quartet system by taking all the quartets displayed by the tree, as two phylogenetic trees on $X$ are isomorphic if and only if the associated quartet systems coincide [18]. Hence, a quartet system that satisfies Theorem 6.3.1 uniquely determines a phylogenetic tree.

### 6.3.2 Graph representations of a ternary map

Suppose we have a symmetric map $\delta : X \times X \times X \to M^{\odot}$. Then we can represent the restriction of $\delta$ to all 3-element subsets of any 5-element subset $\{x, y, z, u, v\}$ of $X$ by an edge-colored complete graph on the 5 vertices $x, y, z, u, v$. For any distinct $a, b, c, d \in \{x, y, z, u, v\}$, edge $ab$ and edge $cd$ have the same color if and only if the value of $\delta$ for $\{x, y, z, u, v\} \setminus \{a, b\}$ is the same as for $\{x, y, z, u, v\} \setminus \{c, d\}$. It follows from Condition (3) in the definition of a symbolic ternary metric that, for any vertex of the graph, either 2 incident edges have one color and the other 2 edges have another color, or all 4 incident edges have the same color. By symmetry, there are exactly five non-isomorphic graph representations.

**Lemma 6.3.2.** *Let the edges of a $K_5$ be colored such that for each vertex, the 4 incident edges are either colored by the same color, or 2 of them colored by one color and the other 2 by another color. Then there are exactly 5 non-isomorphic colorings, and they are depicted in Figure 6.1.*

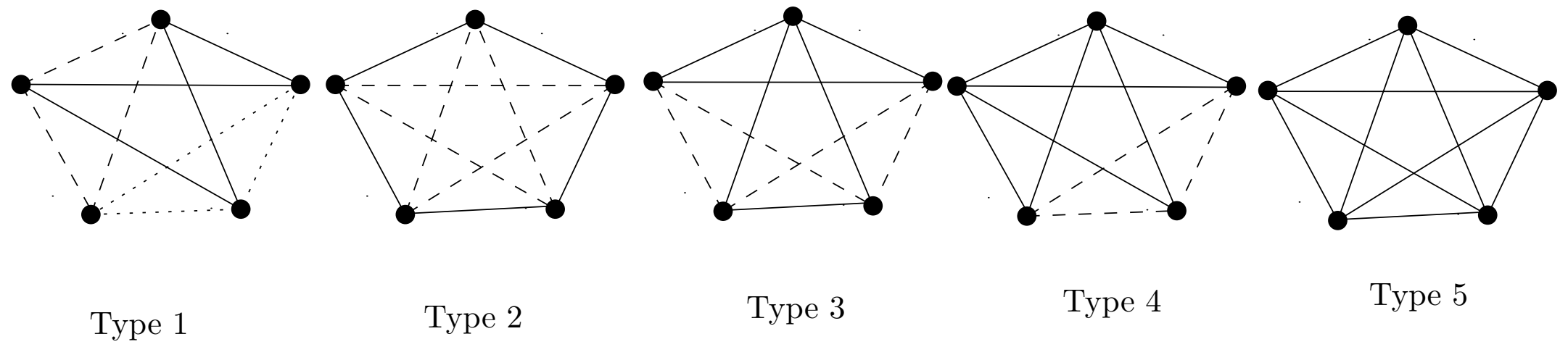

Figure 6.1: The 5 non-isomorphic colorings of $K_5$ for which every color class induces an Eulerian graph. Note that in type 1 there are 3 types of edges, solid edges, dotted edge and dashed edges

*Proof.* It follows from the condition on the coloring that every color class induces an *Eulerian* subgraph (a graph where all vertices have even degree) of $K_5$. Therefore, ignoring isolated vertices, every such induced subgraph either is a cycle or it contains a vertex of degree four. Since there are only ten edges, the only way to have three color classes is two triangles and one 4-cycle. In that case each of the triangles must contain two non-adjacent vertices of the 4-cycle and the vertex that is not in the 4-cycle, thus we get a coloring isomorphic to Type 1 in Figure 6.1. If there are exactly two color classes, then one of them has to be a cycle and the other one its complement. This yields Types 2, 3, and 4, if the length of the cycle is 5, 4, and 3, respectively. Finally, if there is only one color, we get Type 5.

□

Note that the vertices are not labeled and it does not matter which color we are using.

We will prove Theorem 6.2.4 by obtaining a quartet system from any symbolic ternary metric. More precisely, we say that the symbolic ternary metric $\delta$ on $X$ *generates* the quartet $xy|zu$ if either $\delta(x,z,u) = \delta(y,z,u) \neq \delta(x,y,z) = \delta(x,y,u)$, or $|\{\delta(x,y,z), \delta(x,y,u), \delta(x,z,u), \delta(y,z,u)\}| = 1$ and there is $e \in X$ such that

$$\delta(x,y,e) = \delta(x,y,z) = \delta(x,y,u) = \delta(x,z,u) = \delta(z,u,e) = \delta(y,z,u)$$

$$\neq \delta(x,u,e) = \delta(x,z,e) = \delta(y,z,e) = \delta(y,u,e).$$

In the latter case, we say that $e$ *resolves* $x,y,z,u$. Note that the 3-sets obtained by adding $e$ to the pairs of the generated quartet both have the same $\delta$-value as the subsets of $\{x,y,z,u\}$. The following lemma will show that the set of all quartets generated by a symbolic ternary metric is thin.

**Lemma 6.3.3.** *Let $\delta : X \times X \times X \to M^{\odot}$ be a symbolic ternary metric and let $x,y,z,u \in X$ be four different taxa with $|\{\delta(x,y,z), \delta(x,y,u), \delta(x,z,u), \delta(y,z,u)\}| = 1$. Let $e, e' \in X - \{x,y,z,u\}$ such that $\{x,y,z,u,e\}$ and $\{x,y,z,u,e'\}$ are both 4-6-partitioned, and let*

$$\delta(x,y,e) = \delta(x,y,z) = \delta(x,y,u) = \delta(x,z,u) = \delta(z,u,e) = \delta(y,z,u)$$

$$\neq \delta(x,u,e) = \delta(x,z,e) = \delta(y,z,e) = \delta(y,u,e).$$

*Then we also have*

$$\delta(x,y,e') = \delta(x,y,z) = \delta(x,y,u) = \delta(x,z,u) = \delta(z,u,e') = \delta(y,z,u)$$

$$\neq \delta(x,u,e') = \delta(x,z,e') = \delta(y,z,e') = \delta(y,u,e').$$

*Proof.* We already know that $|\{\delta(x,y,z), \delta(x,y,u), \delta(x,z,u), \delta(y,z,u)\}| = 1$ and $\{x,y,z,u,e'\}$ is 4-6-partitioned.

So there are three possible cases for the values of $\delta$ on $\{x,y,z,u,e'\}$.

(1) $\delta(x,y,z) = \delta(x,y,u) = \delta(x,z,u) = \delta(y,z,u)$ and the rest 6 are equal. Then consider $\delta$ on $\{x,y,z,e'\}$, if it is 1-3 partitioned instead of 2-2 partitioned, then it contradicts to the definition of symbolic ternary metric, thus this case would not happen.

(2) $\delta(x,y,z) = \delta(x,y,u) = \delta(x,z,u) = \delta(y,z,u) = \delta(e',a,b) = \delta(e',a,c)$ where $\{a,b,c\} \in \binom{\{x,y,z,u\}}{3}$. There are totally 12 different cases. Since $x,y,z,u$ are symmetric, w.l.o.g., we assume $\delta(x,y,z) = \delta(x,y,u) = \delta(x,z,u) = \delta(y,z,u) = \delta(e',x,y) = \delta(e',x,z)$ and the rest are equal. Then consider $\delta$ on $\{e',x,y,z\}$, if it is 1-3 partitioned

instead of 2-2 partitioned, then it contradicts to the definition of symbolic ternary metric, thus this case would not happen.

(3) $\delta(x,y,z) = \delta(x,y,u) = \delta(x,z,u) = \delta(y,z,u) = \delta(e',a,b) = \delta(e',c,d)$ where $\{a,b,c,d\} \in \binom{\{x,y,z,u\}}{4}$. There are totally 3 different cases. Suppose the statement of the lemma is wrong. Then because $x,y,z,u$ are symmetric w. l. o. g. we can assume that

$$\delta(x,z,e') = \delta(x,y,z) = \delta(x,z,u) = \delta(y,u,e') = \delta(x,y,u) = \delta(y,z,u)$$
$$\neq \delta(x,u,e') = \delta(x,y,e') = \delta(y,z,e') = \delta(z,u,e').$$

Case (3a): $\delta(x,u,e) \neq \delta(x,u,e')$. We assume that $\delta(x,u,e)$ is dashed, $\delta(x,u,e')$ is dotted, and $\delta(x,y,z)$ is solid. Since $\delta$ is a symbolic ternary metric, by Lemma 6.3.2, the graph representation of $\{y,z,u,e,e'\}$ has to be Type 1, so the color classes are one 4-cycle and two 3-cycles. The values of $\delta$ for the sets that contain at most one of $e$ and $e'$ are shown in Figure 6.2. There is a path of length 3 that is colored with $\delta(x,y,z)$ (solid), and there are paths of length 2 colored with $\delta(x,u,e)$ (dashed) and $\delta(x,u,e')$ (dotted), respectively. It follows that we only can get Type 1 by coloring the edges connecting the end vertices of each of those paths with the same color as the edges on the path. We get $\delta(u,e,e') = \delta(x,y,z)$ (solid), $\delta(z,e,e') = \delta(x,u,e')$ (dotted), and $\delta(y,e,e') = \delta(x,u,e)$ (dashed). Now doing the same analysis for $\{x,y,z,e,e'\}$ yields $\delta(z,e,e') = \delta(x,u,e)$, in contradiction to $\delta(z,e,e') = \delta(x,u,e')$.

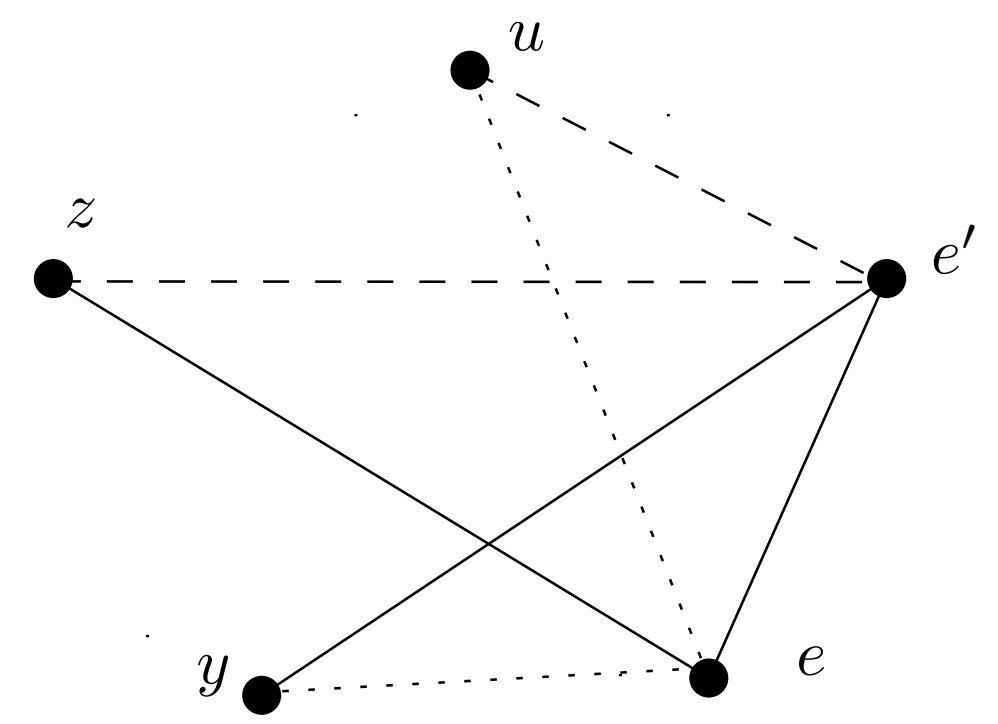

Figure 6.2: The partial coloring of $K_5$ as described in Case (a).

Case (3b): $\delta(x,u,e) = \delta(x,u,e')$. The graph representation of $\{y,z,u,e,e'\}$ can be obtained from Figure 6.2 by identifying the colors dashed and dotted. It contains a path of length 3 of edges colored with $\delta(x,y,z)$ and a path of length 4 of edges colored with $\delta(x,u,e)$. Since a path is not an Eulerian graph, both colors must be used for at least one of the remaining three edges, thus Type 1 is not possible.

Due to $\delta$ being a symbolic ternary map, $\{y,z,u,e,e'\}$ is not 5-5-partitioned, and since there are only two colors with at least 4 edges of one color and 5 edges of the other color, Lemma 6.3.2 implies that the corresponding graph representation must be Type 2 and therefore, $\{y, z, u, e, e'\}$ is 4-6-partitioned. We get $\delta(u, e, e') = \delta(x, y, z)$, $\delta(y, e, e') = \delta(z, e, e') = \delta(x, u, e)$. Now we consider the graph representation of $\{x, z, u, e, e'\}$, and we observe that the edges colored with $\delta(x, y, z)$ contain a path of length 4, and the edges colored with $\delta(x, u, e)$ contain a 5-cycle. Hence, $\{x, z, u, e, e'\}$ is 5-5-partitioned, in contradiction to Condition (4). $\square$

*Proof of Theorem 6.2.4.* By the definitions of the median, the ternary map $\delta_{(T;t)}$ associated with a discriminating symbolically dated phylogenetic tree $(T, t)$ on $X$ satisfies Conditions (1) and (2). For any distinct leaves $x, y, z, u$, the smallest subtree of $T$ connecting those four leaves has at most two vertices of degree larger than two. If there are two such vertices, then each of them has degree 3 and is the median of two 3-sets in $\{x, y, z, u\}$. Therefore, $\delta_{(T;t)}$ satisfies Condition (3).

For any 5 distinct leaves, the smallest subtree of $T$ connecting them either has three vertices of degree 3, or one vertex of degree 3 and one of degree 4, or one vertex of degree 5, while all other vertices have degree 1 or 2. The first case is depicted in Figure 6.3. There $v_1$ is the median for the 3-sets that contain $x_1$ and $x_2$, $v_3$ is the median for the 3-sets that contain $z_1$ and $z_2$, and $v_2$ is the median for the remaining four 3-sets. Hence, either $\{x_1, x_2, y, z_1, z_2\}$ is 4-6-partitioned (if $t(v_1) = t(v_3)$), or there are three different values of $\delta_{(T;t)}$ within those five taxa. For the other two cases, the set of five taxa is either 3-7-partitioned or $\delta_{(T;t)}$ is constant on all its subsets with 3 taxa. Hence, no subset of $X$ of cardinality five is 5-5-partitioned by $\delta_{(T;t)}$, thus $\delta_{(T;t)}$ satisfies Condition (4).

On the other hand, let $\delta$ be a symbolic ternary metric on $X$. By Lemma 6.3.2, taking any 5-element subset of $X$, the possible graph representations of the delta system satisfying (1), (2), and (3) are shown in Figure 6.1. Except for Type 2, all other types satisfy (4).

For the first type, the delta system is $\delta(y, z, u) = \delta(x, y, z) = \delta(w, y, z) \neq \delta(w, x, z) = \delta(w, x, u) = \delta(w, x, y) \neq \delta(x, z, u) = \delta(w, z, u) = \delta(x, y, u) = \delta(w, y, u) \neq \delta(y, z, u)$. The corresponding quartet system is $\{xw|yu, xw|zu, xw|yz, uw|yz, xu|yz\}$.

For the third type, the delta system is $\delta(y, z, u) = \delta(x, z, u) = \delta(w, z, u) = \delta(x, y, w) = \delta(y, u, w) = \delta(y, z, w) \neq \delta(x, y, z) = \delta(w, x, z) = \delta(x, y, u) = \delta(w, x, u)$. The corresponding quartet system is $\{wy|xu, wy|xz, xy|zu, wy|zu, wx|zu\}$.

For the fourth type, the delta system is $\delta(w, x, z) = \delta(w, x, u) = \delta(w, x, y) \neq \delta(y, z, u) = \delta(x, y, z) = \delta(x, z, u) = \delta(w, z, u) = \delta(w, y, z) = \delta(x, y, u) = \delta(w, y, u)$. The corresponding quartet system is $\{xw|yu, xw|zu, xw|yz\}$.

For the fifth type, the delta system is $\delta(y, z, u) = \delta(x, y, z) = \delta(w, y, z) = \delta(w, x, z) = \delta(w, x, u) = \delta(w, x, y) = \delta(x, z, u) = \delta(w, z, u) = \delta(x, y, u) = \delta(w, y, u)$. The corre-sponding quartet system is $\emptyset$.

All quartet systems are thin, transitive, and saturated. Indeed, the delta systems of Types 1 and 3 generate all quartets displayed by a binary tree, Type 4 generates all quartets displayed by a tree with exactly one interior edge, and Type 5 corresponds to the star tree with 5 leaves. Now we take the union of all quartets generated by $\delta$. The resulting quartet system is transitive and saturated, since we have verified these properties for all subsets of five taxa, and it is thin in view of Lemma 6.3.3. Therefore, by Theorem 6.3.1, every delta system satisfying Conditions (1), (2), (3) and (4) uniquely determines a phylogenetic tree $T$ on $X$.

It only remains to show that two 3-element subsets of $X$ that have the same median in $T$ must be mapped to the same value of $\delta$, since then we can define $t$ to be the dating map with $t(v) = \delta(x, y, z)$ for every interior vertex $v$ of $T$ and for all 3-sets $\{x, y, z\}$ whose median is $v$. It suffices to consider two sets which intersect in two taxa, as the general case follows by exchanging one taxon up to three times. We assume that $v$ is the median of both, $\{w, x, y\}$ and $\{w, x, z\}$. If $\delta(w, x, y) \neq \delta(w, x, z)$, then by the definition of symbolic ternary metric, $\{x, y, w, z\}$ is 2-2-partitioned by $\delta$, thus $\delta$ generates one of the quartets $wy|xz$ and $wz|xy$, thus $T$ must display that quartet. However, in both cases $\{w, x, y\}$ and $\{w, x, z\}$ do not have the same median in $T$. We also claim that the associated $t$ is discriminating. Suppose otherwise, there is an edge $uv$ in $T$ such that $t(u) = t(v)$. Thus for any four leaves $x_1, x_2, y_1, y_2$ such that $\text{med}(x_1, x_2, y_1) = \text{med}(x_1, x_2, y_2) = u$ and $\text{med}(y_1, y_2, x_1) = \text{med}(y_1, y_2, x_2) = v$, we have $\delta(x_1, x_2, y_1) = \delta(x_1, x_2, y_2) = \delta(y_1, y_2, x_1) = \delta(y_1, y_2, x_2)$. Suppose there is a leave $e$ which resolves $x_1x_2y_1y_2$, since $uv$ is an edge in $T$, then there exists $i$ and $j$ such that $\text{med}(x_i, y_j, e) \in \{u, v\}$, we have $\delta(x_i, y_j, e) = \delta(x_1, x_2, y_1)$, which means $e$ cannot resolve $x_1x_2y_1y_2$, a contradiction. Hence a thin, transitive and saturated quartet system uniquely determines a discriminating symbolic dated tree. □

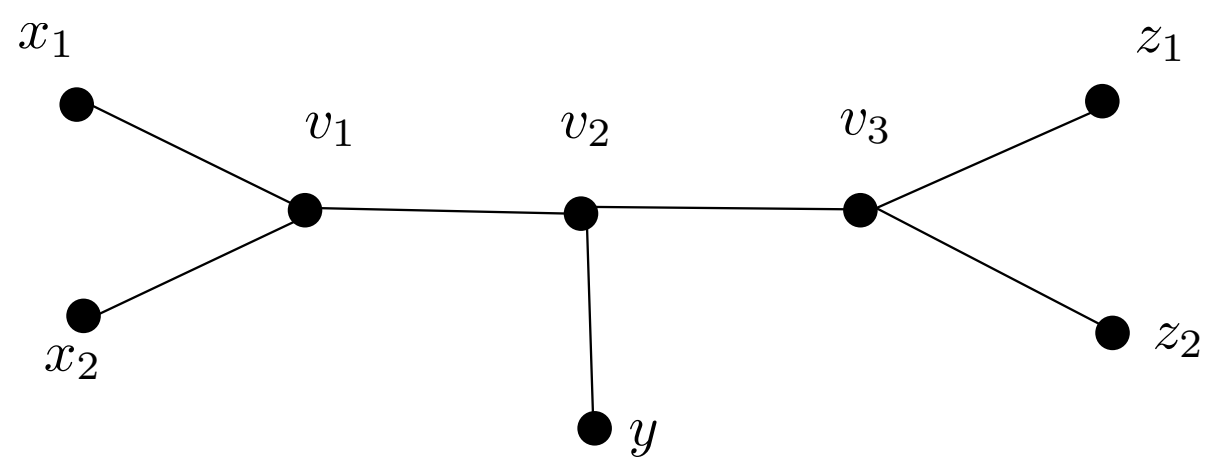

Figure 6.3: The leaves and median vertices for a 5-taxa binary tree.

The set of quartets generated by Type 2 does not satisfy the condition of being saturated from Theorem 6.3.1. Without loss of generality, label the vertices by $x, y, z, u, w$ as in Figure 6.4. Then the quartet system is $\{yw|zu, xu|yz, xz|uw, xy|zw, xw|yu\}$. In order to be saturated, the presence of $xu|yz$ would induce that we have $xw|yz$ or $xu|yw$, but we have $xy|zw$, $xw|yu$ instead.

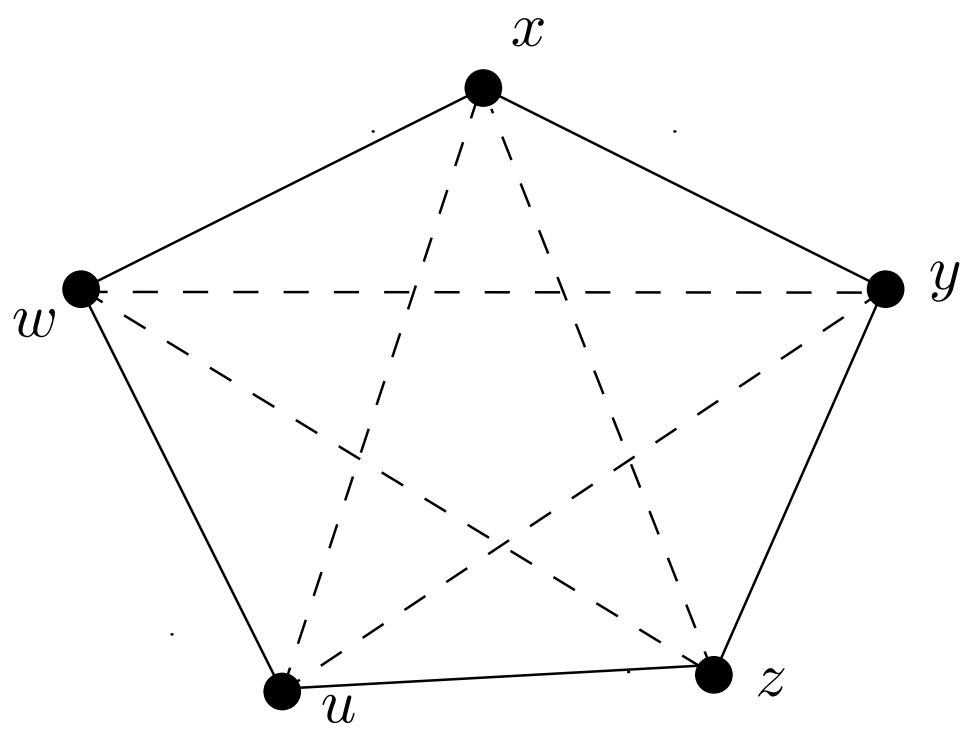

Figure 6.4: The graph representation of a ternary map satisfying Conditions (1), (2), and (3), but not (4).

## 6.4 Reconstructing a binary phylogenetic tree

The aim of this section is to prove Theorem 6.2.5.

A quartet system $Q$ on $X$ is *complete*, if

$$|\{Q \cap \{ab|cd, ac|bd, ad|bc\}\}| = 1$$

holds for all $\{a, b, c, d\} \in \binom{X}{4}$. Using the easy observation that a phylogenetic tree is binary if and only if it displays a quartet for every 4-set, the following result is a direct consequence of Theorem 6.3.1.

**Corollary 6.4.1.** *A quartet system $Q \subseteq Q(X)$ is of the form $Q = Q(T)$ for some binary phylogenetic tree $T$ on $X$ if and only if $Q$ is complete, transitive, and saturated.*

Condition (*) ensures that a ternary metric $\delta$ generates a quartet for every set of four taxa, even if $\delta$ is constant on all of its 3-taxa subsets. In view of Lemma 6.3.3 we also have that $\delta$ can not generate two different quartets for the same 4-set. Hence, we have the following corollary.

**Corollary 6.4.2.** *A symbolic ternary metric $\delta : X \times X \times X \to M^{\odot}$ that satisfies Condition (*) generates a complete quartet system on $X$.*

Now we prove Theorem 6.2.5.

*Proof.* Let $(T, t)$ be a symbolically dated binary phylogenetic tree. By Theorem 6.2.4, $\delta_{(T;t)}$ is a symbolic ternary metric. Since $T$ is binary, it displays a quartet for every 4-taxa subset $\{x, y, z, u\}$ of $X$. Assume that $T$ displays $xu|yz$, thus $\text{med}(x, y, z) \neq \text{med}(x, y, u)$. If $|\{\delta_{(T;t)}(x, y, z), \delta_{(T;t)}(x, y, u), \delta_{(T;t)}(x, z, u), \delta_{(T;t)}(y, z, u)\}| = 1$, then there is at least one vertex $v$ on the path in $T$ connecting $\text{med}(x, y, z)$ and $\text{med}(x, y, u)$ with $t(v) \neq t(\text{med}(x, y, z))$, as $t$ is discriminating. Hence, there is a leaf $e \in X$ such that $v = \text{med}(x, y, e)$. It follows that the set $\{x, y, z, u, e\}$ is 4-6 partitioned by $\delta_{(T;t)}$, thus $\delta_{(T;t)}$ satisfies Condition (*).

On the other hand, if $\delta$ is a symbolic ternary metric on $X$ and satisfies (*), then by Corollary 6.4.2, it corresponds to a unique complete quartet system, thus it encodes a binary phylogenetic tree $T$ in view of Corollary 6.4.1. As in the last paragraph of the proof of Theorem 6.2.4, we can define a dating map $t$ by $t(v) = \delta(x, y, z)$ for every interior vertex $v$ of $T$ and for all 3-sets $\{x, y, z\}$ whose median is $v$. Hence, $(T, t)$ is a symbolically dated binary phylogenetic tree. □

## 6.5 The recognition of pseudo-cherries

In Theorem 6.2.4, we have established a 1-to-1 correspondence between symbolically dated phylogenetic trees and symbolic ternary metrics on $X$, and a bijection is given by mapping $(T, t)$ to $d_{(T,t)}$. To get the inverse of this map, we can first compute the set of all quartets generated by a symbolic ternary metric, and then apply an algorithm that reconstructs a phylogenetic tree from the collection of all its displayed quartets. Finally, the dating map is defined as in our proof of Theorem 6.2.4. This approach would correspond to first extracting rooted triples from a symbolic ultrametric and then reconstruct the rooted tree (see Section 7.6 of [73]). However, a more direct way to reconstruct the corresponding tree from a symbolic ultrametric was presented in [38]. It is based on identifying maximal sets of at least two taxa that are adjacent to the same interior vertex, so-called *pseudo-cherries*. These can iteratively be identified into a single new taxon, thereby reconstructing the corresponding tree in a bottom-up fashion.

Here we show how to find the pseudo-cherries of $T$ from a symbolic ternary metric $d_{(T,t)}$. Then an algorithm to reconstruct $T$ can be designed exactly as in [38]. We point out that it is not necessary to check Condition (4) of a symbolic ternary metric, as a violation would make the algorithm recognize that the ternary map does not correspond to a tree.

Given an arbitrary symbolic ternary map $\delta : X \times X \times X \to M^{\odot}$ satisfying Conditions (1), (2), and (3). For $x, y \in X$ and $m \in M^{\odot}$, we say that $x$ and $y$ are *m-equivalent*, if there is $z \in X$ such that $\delta(x, y, z) = m$, and for $u, v \in X - x - y$, $\delta(x, u, v) = m$ if and only if $\delta(y, u, v) = m$.

**Lemma 6.5.1.** *If $x$ and $y$ are $m$-equivalent and $y$ and $z$ are $m'$-equivalent, then $m = m'$ and $x$ and $z$ are $m$-equivalent.*

*Proof.* Assume $\delta(x, y, z) \neq m$. Then let $u \in X$ with $\delta(x, y, u) = m$. Since $\delta$ is not constant on $\{x, y, z, u\}$, this 4-set must be 2-2-partitioned, thus exactly one of $\delta(x, u, z) = m$ and $\delta(y, u, z) = m$ must hold, in contradiction to $x$ and $y$ being $m$-equivalent. Hence, we have $\delta(x, y, z) = m$, and by symmetry we also have $\delta(x, y, z) = m'$, thus $m = m'$. In order to verify that $x$ and $z$ must be $m$-equivalent, we have already shown $\delta(x, y, z) = m$. For $w, w' \in X - \{x, y, z\}$, we have $\delta(x, w, w') = m$ if and only if $\delta(y, w, w') = m$, since $x$ and $y$ are $m$-equivalent, and we have $\delta(y, w, w') = m$ if and only if $\delta(z, w, w') = m$, since $y$ and $z$ are $m$-equivalent. Finally, we have $\delta(x, y, w) = m$ if and only if $\delta(x, z, w) = m$ if and only if $\delta(y, z, w) = m$. Hence $x$ and $z$ are $m$-equivalent. □

We say $x, y \in X$ are *$\delta$-equivalent*, denoted by $x \sim_\delta y$, if there exists $m \in M^\odot$ such that $x$ and $y$ are $m$-equivalent.

**Lemma 6.5.2.** *The relation of being $\delta$-equivalent is an equivalence relation.*

*Proof.* For any $x \in X$, since $\delta(x, x, y) = \odot$ for any $y \in X$, by definition $x$ and $x$ are $\odot$-equivalent, hence $x$ and $x$ are $\delta$-equivalent. Hence $\sim_\delta$ is reflexive. For any $x \sim_\delta y$, we know that there exists an $m \in M^\odot$ such that $x$ and $y$ are $m$-equivalent. Since $\delta$ is symmetric, by the definition of $m$-equivalent, $y$ and $x$ are also $m$-equivalent, thus $y \sim_\delta x$. Hence, $\sim_\delta$ is symmetric. To prove the transitivity of $\sim_\delta$, assume $x \sim_\delta y$ and $y \sim_\delta z$, by Lemma 6.5.1 we know that $x \sim_\delta z$. Therefore, $\delta$-equivalent is an equivalence relation. □

Suppose that $T$ is a phylogenetic tree on $X$. Let $C \subseteq X$ be a subset of $X$ with $|C| \geq 2$. We call $C$ a *pseudo-cherry* of $T$, if there is an interior vertex $v$ of $T$ such that $C$ is the set of all leaves adjacent to $v$ .

**Theorem 6.5.3.** *If $(T; t)$ is a symbolically dated phylogenetic tree, then a non-empty subset $C$ of $X$ is a non-trivial equivalence class of $\sim_{\delta_{(T,t)}}$ if and only if $C$ is a pseudo-cherry of $T$.*

*Proof.* For the ease of notation, we let $\delta = \delta_{(T,t)}$.

Since $t$ is discriminating, the definition of a pseudo-cherry immediately implies that any pseudo-cherry of $T$ must be contained in a non-trivial equivalence class of $\sim_\delta$.

Conversely, if a non-trivial equivalence class $C$ of $\sim_\delta$ is not a pseudo-cherry, then there are $x_1, x_2 \in C$ such that the path in $T$ that contains $x_1$ and $x_2$ has length at least 3, and since $t$ is discriminating, it has at least 2 interior vertices labeled by two different elements of $M$. Suppose that all elements of $C$ are $m$-equivalent, and that $v$ is an interior vertex on the path from $x_1$ to $x_2$ such that $t(v) = m'$ and $m' \neq m$. Further, let $y \in X$ such that $v = \text{med}(x_1, x_2, y)$. Since $x_1$ and $x_2$ are $m$-equivalent, there is $z \in X$ such that $\delta(x_1, x_2, z) = m$. Then the median $u$ of $x_1, x_2, z$ is also on the path from $x_1$ and $x_2$, and we assume without loss of generality that $u$ is on the path from $x_1$ to $v$. It follows that $\delta(x_1, y, z) = m$ but $\delta(x_2, y, z) = m'$ in contradiction to $x_1 \sim_m x_2$. □

## 6.6 Discussions and open questions

The proofs of our main results heavily rely on extracting the corresponding quartet set from a symbolic ternary metric and then checking that our Conditions (3) and (4) guarantee the quartet system to be thin, transitive, and saturated, and adding (*) makes the quartet system complete. The conditions look like (3) corresponds to thin and transitive, (4) to saturated, and (*) to complete. However, this is not true, and removing (4) from Theorem 6.2.5 does not necessarily yield a transitive complete quartet system. While for five taxa, a 5-5-partition yields the only non-saturated complete transitive quartet system, Lemma 6.3.3 does not hold without Condition (4). Indeed the ternary map that is visualized in Figure 6.5 suffices Conditions (1), (2), (3), and (*), but it generates two quartets on each of $\{a_1, a_2, b_1, b_2\}$ and $\{a_1, a_2, c_1, c_2\}$. It can be shown by checking the remaining 5-sets in Case 2 of our proof of Lemma 6.3.3 that every ternary map on 6 taxa satisfying Conditions (1), (2), (3), (*) that does not yield a thin quartet system is isomorphic to this example. This raises the question whether ternary maps satisfying these four conditions can be completely characterized. The hope is to observe something similar to the Clebsch trees that were observed by Jan Weyer-Menkhoff [82]. As a result, a phylogenetic tree with all interior vertices of degree 3, 5, or 6 can be reconstructed from every transitive complete quartet set.

Another direction to follow up this work would be to consider more general graphs than trees. A *median graph* is a graph for which every three vertices have a unique median. Given a vertex-colored median graph and a subset $X$ of its vertex set, we can get a symmetric ternary map on $X \times X \times X$ by associating the color of the median to every 3-subset of $X$. It would be interesting to see whether this map can be used to reconstruct the underlying graph for other classes of median graphs than phylogenetic trees. In phylogenetics, median graphs are used to represent non-treelike data. Since the interior vertices of those so-called *splits graphs* do in general not correspond to any ancestor of some of the taxa, reconstructing a collection of splits from the ternary map induced by a vertex-colored splits graph is probably limited to split systems that are almost compatible with a tree.

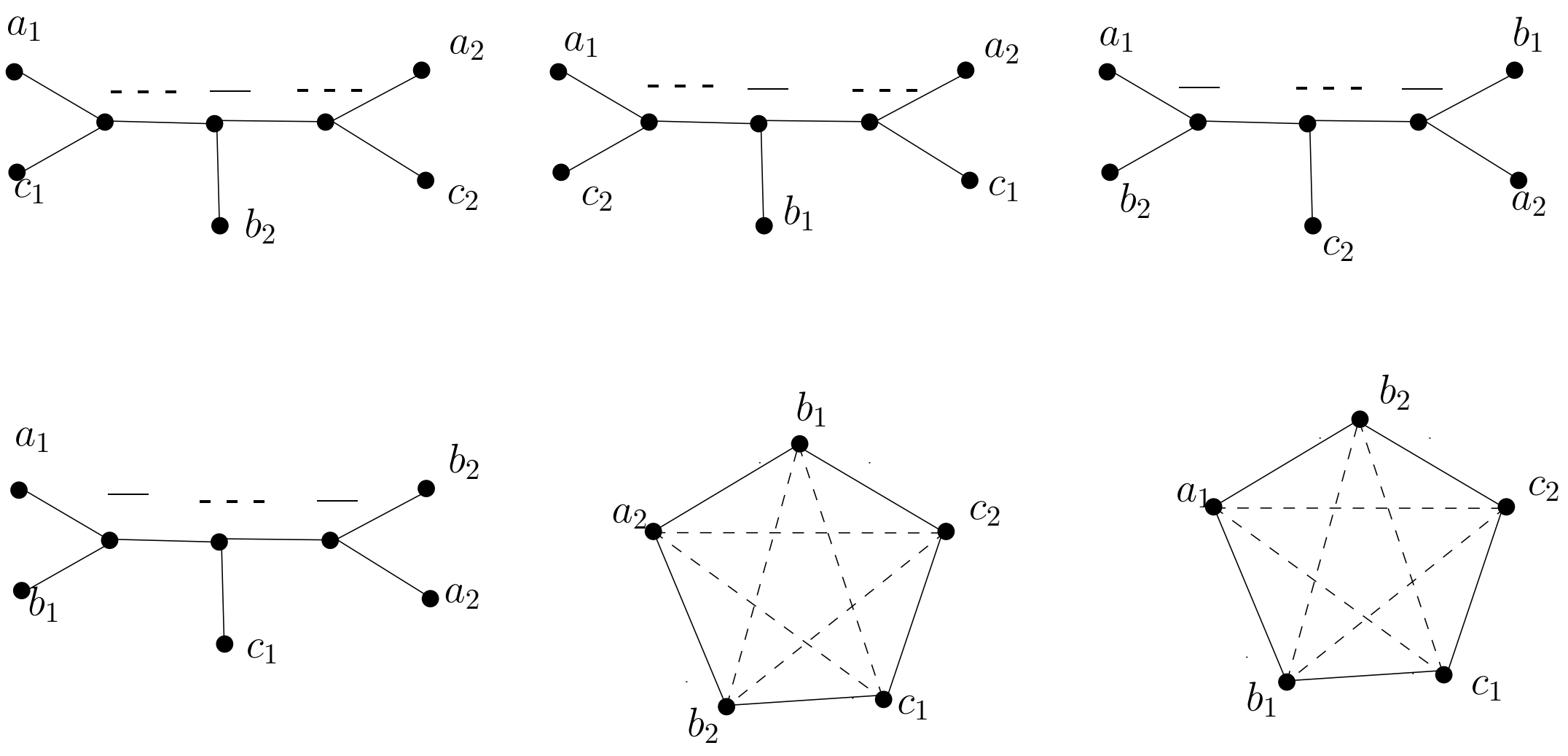

Figure 6.5: The 5-taxa trees respectively graph representations generated by a ternary map satisfying Conditions (1), (2), (3), (*) but not (4).

It is one of the main observations of [38] that the 4-point condition for symbolic ultrametrics can be formulated in terms of *cographs* which are graphs that do not contain an induced path of length 3. For the special case that $\delta : X \times X \to M^{\odot}$ with $|M| = 2$ and $m \in M$, consider the graph with vertex set $X$ where two vertices $x, y$ are adjacent, if and only if $\delta(x, y) = m$. Then deciding whether $\delta$ is a symbolic ultrametric can be reduced to checking whether this graph (as well as its complement) is a cograph. This is useful for analyzing real data which will usually not provide a perfect symbolic ultrametric, thus some approximation is required. For ternary maps and unrooted trees, the 5-taxa case looks promising, as Condition (4) translates to a forbidden graph representation that splits the edges of a $K_5$ into two 5-cycles, thus we have a self-complementary forbidden induced subgraph. However, for more taxa, the 3-sets that are mapped to the same value of a ternary map $\delta$ define a 3-uniform hypergraph on $X$ and formulating Condition (4) in terms of this hypergraph does not seem to be promising. In addition, even if there are only two values of $\delta$ for 3-sets, Condition (3) does not become obsolete. We leave it as an open question, whether an alternative characterization of symbolic ternary metrics exists that makes it easier to solve the corresponding approximation problem.

# 博士后在站期间研究成果

一、博士后在站期间发表和在审的论文列表

1. Inferring Phylogenetic Trees from the Knowledge of Rare Evolutionary Events, joint with M. Hellmuth, M. Hernandez-Rosales, P. F. Stadler, to appear in *Journal of Mathematical Biology*, available at *https://arxiv.org/abs/1612.09093*

2. Reconstructing unrooted phylogenetic trees from symbolic ternary metrics, joint with S. Grünewald, Y. Wu, minor revision in *Bulletin of Mathematical Biology*, available at *https://arxiv.org/abs/1702.00190*

3. Fractal property of homomorphism order, joint with J. Fiala, J. Hubička, J. Nešetřil, *European Journal of Combinatorics*, in press, *http://www.sciencedirect.com/science/article/pii/S0195669817300914*

4. Gaps in full homomorphism order, joint with J. Fiala, J. Hubička, *Electronic Notes in Discrete Mathematics*, (2017), **61**, 429–435

5. An universality argument for graph homomorphisms, joint with J. Fiala, J. Hubička, *Electronic Notes in Discrete Mathematics*, (2015), **49**, 643–649.

6. Word-representability of split graphs, joint with Sergey Kitaev, Jun Ma, Hehui Wu, submitted.

二、博士后在站期间主持和参与的科研项目

◇ 2016 – 2017, **Principal Investigator** of the project "Some problems on graph homomorphisms", No. 154463
**50000 RMB**, National Natural Science Fund of Chinese Postdoc, CN.

◇ 2013 – 2016 **Participant** of the project "Several problems on discrete process: structure, simulation, and algorithm" (11271255)
**600000 RMB**, National Natural Science Fund, CN.

◇ 2017 – 2020 **Main Participant** of the project "Discrete Dynamics and Beyond" (11671258)
**480000 RMB**, National Natural Science Fund, CN.

# 致 谢

首先我要深深地感谢我的合作导师吴耀琨老师，吴老师对科研有激情，投入，懂的方向很多，在他的带领下我学习和接触了很多新的领域。吴老师组里聚集着一群聪明，有想法有活力的本科生和研究生们，组里讨论的氛围非常好，使我受益匪浅。虽然我比这些学生们学数学的时间长，但是从他们那里我学到了很多，特别是一些新的想法和思路，给我不少启发。

感谢我的合作者们，感谢我的博士导师Peter F. Stadler,通过和他们合作我学到了很多东西。

感谢一些其他关心和指导过我的老师，武同锁老师，王维克老师，张晓东老师，厦门大学金贤安老师，湖南大学彭岳建，复旦大学吴河辉老师。感谢数学院行政部门的老师，郭恒亮老师，尚建辉老师，陈青老师等等，他们对待工作专业和负责，给我提供了很多工作上的协助。

感谢我的办公室同事们，吴德军老师和王晓东博士，他们给予了我不少帮助。

感谢我在上海的好友们，是他们让上海的两年成为我近几年生活上最开心的日子。感谢父母和家人一直以来的支持。

龙旸靖

二零一七年九月于上海